\documentclass[a4paper,11pt]{article}
\usepackage{amsmath, amssymb, amsthm}
\usepackage{graphicx} 

\begin{document}

\def\PP{{\mathcal P}} \def\HH{{\mathcal H}} \def\MM{{\mathcal M}}
\def\S{{\mathbb S}} \def\C{{\mathbb C}} \def\CC{{\mathcal C}}
\def\Aut{{\rm Aut}} \def\Im{{\rm Im}} \def\Ker{{\rm Ker}}
\def\D{{\mathbb D}} \def\QQ{{\mathcal Q}} \def\FF{{\mathcal F}}
\def\SS{{\mathcal S}} \def\LL{{\mathcal L}} \def\Ind{{\rm Ind}}
\def\Norm{{\mathcal G}} \def\R{{\mathbb R}} \def\Z{{\mathbb Z}}
\def\N{{\mathbb N}} \def\Sym{{\rm Sym}} \def\Diff{{\rm Diff}}
\def\T{{\mathbb T}} \def\Isom{{\rm Isom}} \def\BB{{\mathcal B}}
\def\Id{{\rm Id}} \def\mod{{\rm mod}} \def\PDiff{{\rm PDiff}}
\def\Stab{{\rm Stab}} \def\exp{{\rm exp}} \def\A{{\mathbb A}}
\def\UU{{\mathcal U}} \def\RR{{\mathcal R}} \def\source{{\rm source}}
\def\target{{\rm target}} \def\DD{{\mathcal D}} \def\Homeo{{\rm Homeo}}
\def\Norm{{\rm Norm}} \def\NNorm{{\rm \widetilde{N}orm}}


\title{\bf{Roots in the mapping class groups}}
 
\author{\textsc{Christian Bonatti and Luis Paris}}

\date{\today}

\maketitle

\begin{abstract}
\noindent
The purpose of this paper is the study of the roots in the mapping class groups. Let $\Sigma$ be a compact 
oriented surface, possibly with boundary, let $\PP$ be a finite set of punctures in the interior of 
$\Sigma$, and let $\MM (\Sigma, \PP)$ denote the mapping class group of $(\Sigma, \PP)$. We prove that, 
if $\Sigma$ is of genus $0$, then each $f \in \MM (\Sigma)$ has at most one $m$-root for all $m \ge 1$. 
We prove that, if $\Sigma$ is of genus $1$ and has non-empty boundary, then each $f \in \MM (\Sigma)$ has at most one 
$m$-root up to conjugation for all $m \ge 1$. We prove that, however, if $\Sigma$ is of genus $\ge 2$, then 
there exist $f,g \in \MM (\Sigma, \PP)$ such that $f^2=g^2$, $f$ is not conjugate to $g$, and none 
of the conjugates of $f$ commutes with $g$. Afterwards, we focus our study on 
the roots of the pseudo-Anosov elements. We prove that, if $\partial \Sigma \neq \emptyset$, then each 
pseudo-Anosov element $f \in \MM(\Sigma, \PP)$ has at most one $m$-root for all $m \ge 1$. We prove that, however, if 
$\partial \Sigma = \emptyset$ and the genus of $\Sigma$ is $\ge 2$, then there exist two pseudo-Anosov 
elements $f,g \in \MM (\Sigma)$ (explicitely constructed) such that $f^m=g^m$ for some $m\ge 2$, $f$ is 
not conjugate to $g$, and none of the conjugates of $f$ commutes with $g$. Furthermore, if the genus of 
$\Sigma$ is $\equiv 0\ (\mod\,4)$, then we can take $m=2$. Finally, we show that, if $\Gamma$ is a pure 
subgroup of $\MM (\Sigma, \PP)$ and $f \in \Gamma$, then $f$ has at most one $m$-root in $\Gamma$ for 
all $m \ge 1$. Note that there are finite index pure subgroups in $\MM (\Sigma, \PP)$.
\end{abstract}

\noindent
{\bf AMS Subject Classification:} Primary 57M99. Secondary 57N05, 57R30.  

\section{Introduction}

Throughout the paper $\Sigma$ will denote a compact oriented surface, possibly with boundary, and $\PP= 
\{ P_1, \dots, P_n \}$ a finite collection of points, called the {\it punctures}, in the interior of 
$\Sigma$. We denote by $\Diff (\Sigma, \PP)$ the group of orientation preserving diffeomorphisms $F: 
\Sigma \to \Sigma$ which are the identity on a neighborhood of the boundary of $\Sigma$ and such 
that $F(\PP)= \PP$. The {\it mapping class group} of $(\Sigma, \PP)$ is defined to be the group $\MM 
(\Sigma, \PP)= \pi_0 (\Diff (\Sigma, \PP))$ of isotopy classes of elements of $\Diff (\Sigma, \PP)$. 
Note that, in the above definition, one may replace $\Diff (\Sigma, \PP)$ by $\Homeo (\Sigma, \PP)$, 
the group of orientation preserving homeomorphisms $F: \Sigma \to \Sigma$ which are the identity on 
$\partial \Sigma$ and such that $F (\PP)= \PP$ (namely, $\MM (\Sigma, \PP)= \pi_0 (\Diff (\Sigma, 
\PP))= \pi_0( \Homeo (\Sigma, \PP))$. Note also that the hypothesis that the elements of $\Diff 
(\Sigma, \PP)$ restrict to the identity on a neighborhood of $\partial \Sigma$ 
(and not only on $\partial \Sigma$)
is especially needed when 
considering subsurfaces. Indeed, if $\Sigma'$ is a subsurface of $\Sigma$ such that $\PP \cap \partial 
\Sigma' = \emptyset$, then the embedding $\Sigma' \subset \Sigma$ determines an embedding $\Diff 
(\Sigma', \PP \cap \Sigma') \to \Diff (\Sigma, \PP)$ by extending each diffeomorphism $F \in \Diff 
(\Sigma', \PP \cap \Sigma')$ with the identity map outside $\Sigma'$, and this monomorphism determines a 
homomorphism $\MM (\Sigma', \PP \cap \Sigma') \to \MM (\Sigma, \PP)$ which is injective in most of the 
cases but not always (see \cite{ParRol1}).

\bigskip\noindent
If $\Sigma= \D$ is the standard disk and $| \PP |=n$, then $\MM (\Sigma, \PP)$ is the braid group 
$\BB_n$ introduced by Artin \cite{Artin1}, \cite{Artin2}.

\bigskip\noindent
Let $\Gamma$ be a group, $g \in \Gamma$, and $m \ge 1$. Define a {\it $m$-root} of $g$ to be an element 
$f \in \Gamma$ such that $f^m=g$. If $\Gamma = \Z^q$, then each element $g \in \Gamma$ has at most one 
$m$-root. The same result is true if $\Gamma$ is free or, more generally, if $\Gamma$ is biorderable 
(see Subsection 3.1).

\bigskip\noindent
The purpose of the present paper is the study of the roots (or, more precisely, the uniqueness of 
the roots) in the mapping class groups. Our starting point is a recent result due to Gonz\'alez-Meneses 
\cite{Gonza1} which asserts that, in the braid group $\BB_n$, a $m$-root is unique up to conjugation. 
This result, also known as the {\it Makanin conjecture}, has been for some time a classical problem in 
the field (see \cite{BaMySh1}).

\bigskip\noindent
The braid group has different definitions. One of them is as the mapping class group 
of the punctured disk, but this is not the only one. So, it is not a surprise that many results on the 
braid groups (the linearity, for example) are not known for the mapping class groups because their 
proofs are based on different approaches. But, in what concerns the result of Gonz\'alez-Meneses, the 
proof strongly uses the so-called {\it Nielsen-Thurston classification} which is one of the main tools in the 
theory of the mapping class groups.

\bigskip\noindent
We briefly recall this classification and refer to Section 2 for detailed definitions and properties. 
An {\it essential curve} is defined to be an embedding $a: \S^1 \hookrightarrow \Sigma \setminus \PP$ 
of the standard circle $\S^1$ which is not parallel to a boundary component of $\Sigma$ and which does not 
bound a disk embedded in $\Sigma$ containing $0$ or $1$ puncture. The mapping class group $\MM (\Sigma, 
\PP)$ acts on the set $\CC_0 (\Sigma, \PP)$ of isotopy classes of essential curves. We call an element $f \in 
\MM (\Sigma, \PP)$ {\it periodic} if some non-trivial power of $f$ acts trivially on $\CC_0 
(\Sigma, \PP)$, we call $f$ {\it pseudo-Anosov} if it has no finite orbit in $\CC_0 (\Sigma, \PP)$, and 
we call $f$ {\it reducible} otherwise. Note that a Dehn twist along a boundary component acts trivially 
on $\CC_0 (\Sigma, \PP)$, so, according to the above definition, it is periodic. In particular, a 
periodic element does not need to be of finite order. Actually, if $\partial \Sigma \neq \emptyset$, 
then $\MM (\Sigma, \PP)$ is torsion free (see \cite{RouWie1}). 
When $\partial \Sigma = 
\emptyset$, the elements of all these three classes (periodic, pseudo-Anosov, reducible) can be 
represented by diffeomorphisms having some special properties, and these special representatives play a 
crucial role in the understanding of their corresponding mapping classes.

\bigskip\noindent
Now, as pointed out by Birman \cite{Birma2}, the result of Gonz\'alez-Meneses cannot be extended to the 
mapping class groups. The most obvious reason is because, if $\partial \Sigma = \emptyset$, then $\MM 
(\Sigma, \PP)$ has torsion. For example, if $f,g \in \MM (\Sigma, \PP)$ are such that $f$ is of 
order $2$ and $g$ is of order $3$, then $f^6=g^6=\Id$, but $f$ and $g$ are not conjugate. More 
elementary, if $f \in \MM(\Sigma, \PP)$ is of order $m \ge 2$, then $f^mù= \Id^m = \Id$ but $f$ and 
$\Id$ are not conjugate. However, if $\partial \Sigma \neq \emptyset$, the mapping class group $\MM 
(\Sigma, \PP)$ is torsion free. Besides, our first positive result is the following.

\bigskip\noindent
{\bf Proposition 3.2.} {\it Let $\Sigma$ be a surface of genus $0$. Then $\MM (\Sigma)$ is biorderable. 
In particular, if $f,g \in \MM (\Sigma)$ are such that $f^m=g^m$ for some $m \ge 1$, then $f=g$.}

\bigskip\noindent
The case of the surfaces of genus $1$ is not so nice, but we still have a positive answer:

\bigskip\noindent
{\bf Theorem 3.6.} {\it Let $\Sigma$ be a surface of genus $1$ with non-empty boundary. Let $f,g \in 
\MM (\Sigma)$. If $f^m=g^m$ for some $m \ge 1$, then $f$ and $g$ are conjugate.}

\bigskip\noindent
The proof of Theorem 3.6 is similar to the one of the theorem of Gonz\'alez-Meneses in \cite{Gonza1} 
for the braid groups in the sense that it is based on a deep analysis of the different possible reductions 
for a mapping class. However, we also have to determine the conjugacy classes of the periodic elements 
(and of some special reducible elements), while such a step is not needed in the proof of Gonz\'alez-Meneses 
because the conjugacy classes of the periodic elements in the braid groups where previously 
determined (see \cite{ConKol1}, \cite{Eilen1}, \cite{Kerek1}).

\bigskip\noindent
Unfortunately, the good news stop at the genus 1. Indeed, we prove the following.

\bigskip\noindent
{\bf Corollary 4.3.} {\it Let $\Sigma$ be a surface of genus $\rho \ge 2$ and $q$ boundary components, 
where $q \ge 1$, and let $\PP$ be a finite set of punctures in the interior of $\Sigma$. Then there 
exist $f,g \in \MM (\Sigma, \PP)$ such that $f ^2=g^2$, $f$ and $g$ are not conjugate, and none of the 
conjugates of $f$ commutes with $g$.}

\bigskip\noindent
Our next step is to focus our study on the roots of the pseudo-Anosov elements. In the case where 
$\partial \Sigma \neq \emptyset$ we prove the following.

\bigskip\noindent
{\bf Theorem 4.5.} {\it Let $\Sigma$ be a surface of genus $\rho \ge 0$ and $q$ boundary components, 
where $q \ge 1$, and let $\PP$ be a finite set of punctures in the interior of $\Sigma$. Let $f,g \in 
\MM (\Sigma, \PP)$ be two pseudo-Anosov elements. If $f^m=g^m$ for some $m \ge 1$, then $f=g$.}

\bigskip\noindent
However, if $\partial \Sigma = \emptyset$, we may have pretty bad examples:

\bigskip\noindent
{\bf Theorem 5.2.} {\it 
\begin{enumerate}
\item
Let $\Sigma$ be a closed surface of genus $\rho \ge 2$. Then there exist two pseudo-Anosov elements $f, 
g \in \MM (\Sigma, \PP)$ such that $f^{2 (\rho+1)} = g^{2 (\rho+1)}$, $f$ is not conjugate to $g$, and 
none of the conjugates of $f$ commutes with $g$.
\item
Let $\Sigma$ be a closed surface of genus $\rho \ge 4$, with $\rho \equiv 0\ (\mod\, 4)$. Then there 
exist two pseudo-Anosov elements $f,g \in \MM (\Sigma)$ such that $f^2=g^2$, $f$ is not conjugate to 
$g$, and none of the conjugates of $f$ commutes with $g$.
\end{enumerate}}

\bigskip\noindent
In particular, Theorem 5.2 contradicts the ``popular idea'' that two pseudo-Anosov elements either 
commute or generate a free group: the examples of Theorem 5.2 ``highly'' do not commute and ``highly'' 
do not generate a free group. Nevertheless, the set of $g \in \MM (\Sigma, \PP)$ having a common power 
with a given pseudo-Anosov element $f$ is finite, namely:

\bigskip\noindent
{\bf Proposition 5.1.} {\it Let $\Sigma$ be a closed surface of genus $\rho \ge 2$, and let $\PP$ be a 
finite set of punctures in $\Sigma$. Let $f \in \MM (\Sigma, \PP)$ be a pseudo-Anosov element. Then the 
set of $g \in \MM (\Sigma, \PP)$ satisfying $g^m=f^m$ for some $m \ge 1$ is finite.}

\bigskip\noindent
Section 6 of the paper is dedicated to the problem of the roots in some special (finite index) 
subgroups of $\MM (\Sigma, \PP)$. More precisely, we prove the following.

\bigskip\noindent
{\bf Theorem 6.1.} {\it Let $\Sigma$ be a surface of genus $\rho \ge 0$ and $q$ boundary components, 
where $q \ge 0$, and let $\PP$ be a finite set of punctures in the interior of $\Sigma$. Let $\Gamma$ 
be a pure subgroup of $\MM (\Sigma, \PP)$. If $f, g \in \Gamma$ are such that $f^m=g^m$ for some $m \ge 
1$, then $f=g$.}

\bigskip\noindent
We refer to Section 2 for the definition of a pure subgroup. 
However, we insist on the fact that $f,g$ must belong to the same pure subgroup in the above statement. 
The hypothesis that both, $f$ and $g$, are pure does not suffice to get the result. For instance, the 
pseudo-Anosov elements are pure while Theorem 5.2 shows that there exist two different 
(and, even worse, non-conjugate)
pseudo-Anosov elements $f,g$ such that $f^2=g^2$.
The standard example of a pure subgroup is the so-called {\it Ivanov subgroup} 
$\Gamma(k)$ which is defined to be the kernel of the natural homomorphism $\MM (\Sigma, \PP) \to \Aut 
(H_1 (\Sigma \setminus \PP, \Z/ k\Z))$, where $k \ge 3$ (see \cite{Ivano2}). Note that $\Gamma(k)$ has 
finite index in $\MM (\Sigma, \PP)$, thus there are finite index subgroups in $\MM (\Sigma, \PP)$ which 
have the property that each of its elements has at most one $m$-root
inside the subgroup, for all $m \ge 1$. This property is 
a typical property of the biorderable groups, and it is an open question to know whether $\Gamma(k)$ is 
biorderable (see \cite{DDRW1}, \cite{Paris1}).

\bigskip\noindent
In order to achieve our main constructions we develop in Section 2 some more or less original 
techniques which may be interesting by themselves for the reader.

\bigskip\noindent
As the reader may notice, we tried to make a treatment of the question of the roots as complete as 
possible for the mapping class groups of the unpunctured surfaces, and we considered the punctured 
surfaces only if needed or if the proof do not require any extra argument. So, some problems which 
essentially concern the mapping class groups of the punctured surfaces are left. Here are some.

\bigskip\noindent
{\bf Question 1.} Let $\Sigma$ be a surface of genus $0$ and $q$ boundary components, where $q \ge 1$, 
and let $\PP$ be a finite set of punctures in the interior of $\Sigma$. We suspect that, if $f,g \in 
\MM (\Sigma, \PP)$ are such that $f^m=g^m$ for some $m \ge 1$, then $f$ and $g$ are conjugate. This is 
true if either $q=1$ (by \cite{Gonza1}) or $\PP=\emptyset$ (by Proposition 3.2).

\bigskip\noindent
{\bf Question 2.} Let $\Sigma$ be a surface of genus 1 with $q$ boundary components, where $q \ge 1$, 
and let $\PP$ be a non-empty set of punctures in the interior of $\Sigma$. We do not know whether, if 
$f,g \in \MM (\Sigma, \PP)$ are such that $f^m=g^m$ for some $m \ge 1$, then $f$ and $g$ are conjugate.

\bigskip\noindent
{\bf Question 3.} Let $\Sigma$ be a closed surface of genus $\rho \ge 2$, and let $\PP$ be a non-empty 
finite set of punctures in $\Sigma$. We do not know 
any necessary and sufficient condition on $(\Sigma, \PP)$ so that
there exist two pseudo-Anosov elements $f,g 
\in \MM (\Sigma, \PP)$ such that $f^m=g^m$ for some $m \ge 1$, $f$ and $g$ are not conjugate, and none 
of the conjugates of $f$ commutes with $g$. If $|\PP|=1$, then two pseudo-Anosov elements $f,g \in \MM 
(\Sigma, \PP)$ which satisfy $f^m=g^m$ for some $m \ge 1$ must commute. (This can be easily proved using 
the techniques of the present paper.)

\section{Preliminaries}

\subsection{The complex of curves and the Nielsen-Thurston classification}

Let $\S^1= \{ z \in \C; |z|=1 \}$ denote the standard circle. A {\it simple closed curve} of $(\Sigma, 
\PP)$ is an embedding $c: \S^1 \to \Sigma \setminus \PP$. Here we assume that the circle 
$\S^1$ as well as the simple closed curves are oriented. The curve with the same image as $c$ but 
reverse orientation is denoted by $c^{-1}$ (that is, $c^{-1} (z) = c(\bar z)$ for all $z \in \S^1$). 
We will also often 
identify a curve $c$ with its oriented image $c(\S^1)$. The orientation of the surface $\Sigma$ determines an 
orientation of each component of $\partial \Sigma$ which, by this way, is viewed as a simple closed 
curve. Two simple closed curves $a$ and $b$ are {\it isotopic} if there is a continuous family $\{ 
a_t\}_{ t \in [0,1]}$ of simple closed curves such that $a_0=a$ and $a_1=b$.

\bigskip\noindent
Let $a : \S^1 \to \Sigma \setminus \PP$ be a simple closed curve such that $a \cap \partial \Sigma = 
\emptyset$. A {\it tubular neighborhood} of $a$ is an embedding $A: [0,1] \times \S^1 \to \Sigma 
\setminus \PP$ such that $A({1 \over 2}, z)= a(z)$ for all $z \in \S^1$. Choose such a tubular 
neighborhood and 
a smooth function $\varphi: [0,1] \to \R$ such that $\varphi(0)=0$, $\varphi(1) = 1$, $\varphi'(0)= 
\varphi'(1)= 0$, and $\varphi'(t) \ge 0$ for all $t \in [0,1]$, and
consider the diffeomorphism $T: \Sigma \to \Sigma$ which is the identity outside the 
image of $A$ and which is defined by
\[
(T \circ A)(t, z)=A(t, \exp(2i\pi \varphi(t)) \cdot z)
\]
inside the image of $A$. The {\it Dehn twist} along $a$ is defined to be the isotopy class $\tau_a \in 
\MM (\Sigma, \PP)$ of $T \in \Diff (\Sigma, \PP)$. Note that $\tau_a$ does not depend neither on the choice of 
the tubular neighborhood, 
nor on the choice of the function $\varphi$. Moreover,
we have $\tau_a = \tau_b$ if $a$ is isotopic to either $b$ or $b^{-1}$. 
If $a$ bounds a disk embedded in $\Sigma$ which contains $0$ or $1$ puncture, then $\tau_a=1$. 
Otherwise, $\tau_a$ has infinite order. If $a: \S^1 \to \Sigma \setminus \PP$ is a simple closed curve 
which intersects $\partial \Sigma$ (for example, if $a$ is a boundary component), then the {\it Dehn 
twist} along $a$, denoted by $\tau_a$, is defined to be the Dehn twist along some simple closed curve 
$b: \S^1 \to \Sigma \setminus \PP$ isotopic to $a$ and satisfying $b \cap \partial \Sigma = \emptyset$.

\bigskip\noindent
{\bf The complex of curves.}
A simple closed curve $a: \S^1 \to \Sigma \setminus \PP$ is called {\it essential} if neither $a$, nor 
$a^{-1}$, is isotopic to a boundary component of $\Sigma$, and if it does not bound a disk embedded in 
$\Sigma$ containing $0$ or $1$ puncture. Two simple closed curves $a$ and $b$ are {\it equivalent} if 
either $a$ or $a^{-1}$ is isotopic to $b$. The equivalence class of a curve $a$ is denoted by $\langle 
a \rangle$. The vertices of the {\it complex of curves} $\CC (\Sigma, \PP)$ are the equivalence classes 
of essential curves. A subset $\{ \gamma_0, \gamma_1, \dots, \gamma_p \}$ of vertices is a $p$-simplex 
of $\CC (\Sigma, \PP)$ if there exist essential curves $c_0, c_1, \dots, c_p$ such that $\gamma_i= 
\langle c_i \rangle$ for all $0 \le i\le p$, $c_i \cap c_j = \emptyset$ for all $0 \le i\neq j\le p$, 
and $c_i$ is not isotopic to $c_j^{\pm 1}$ (namely, $\gamma_i \neq \gamma_j$) for all $0 \le i\neq j\le 
p$. It is a (non-trivial) fact that $\CC (\Sigma, \PP)$ is a well-defined simplicial complex. 
Furthermore, the mapping class group $\MM (\Sigma, \PP)$ acts naturally on $\CC (\Sigma, \PP)$, and this 
action is faithful if $\partial \Sigma = \emptyset$.

\bigskip\noindent
The complexes of curves have been introduced by Harvey \cite{Harve1}, \cite{Harve2} and are a 
fundamental tool in the study of the mapping class groups. They are simply connected in general (see 
\cite{Harer1}) and, if $\partial \Sigma = \emptyset$, then $\Aut (\CC (\Sigma, \PP))$ coincides with 
the so-called extended mapping class group (which is defined in the same way as the mapping class group 
but with diffeomorphisms $F: \Sigma \to \Sigma$ which do not need to preserve the orientation) (see 
\cite{Ivano1}, \cite{Korkm1}, \cite{Luo1}). They may also be used to defined the so-called {\it 
Nielsen-Thurston classification} of the elements of $\MM (\Sigma, \PP)$ as we turn now to do.

\bigskip\noindent
{\bf The Nielsen-Thurston classification.} Let $f \in \MM (\Sigma, \PP)$. We call $f$ {\it pseudo-Anosov} 
if it has no finite orbit in $\CC (\Sigma, \PP)$, we call $f$ {\it periodic} if $f^m$ acts trivially on 
$\CC (\Sigma, \PP)$ for some $m \ge 1$, and we call $f$ {\it reducible} otherwise.

\bigskip\noindent
As pointed out in the introduction, this classification is of importance in the theory because, when
$\partial \Sigma = \emptyset$, then the elements of all these three classes (periodic, pseudo-Anosov,
reducible) can be represented by diffeomorphisms having some special properties:
a pseudo-Anosov element is represented by a pseudo-Anosov diffeomorphism 
(see Subsection 2.2), a periodic element is represented by a finite order diffeomorphism (see 
Subsection 2.8), and a reducible element is represented by a diffeomorphism which globally leaves 
invariant a family of pairwise disjoint essential curves. Cutting the surface along these curves, this 
diffeomorphism can be viewed (in some sense) as a diffeomorphism of a ``simpler'' surface (but not 
necessarily connected). Some reducible elements, called the pure elements, will be considered in 
Subsection 2.7 and in Section 6.

\bigskip\noindent
Now, observe that the Dehn twist $\tau_c$ along a boundary component $c$ acts trivially on 
$\CC(\Sigma, \PP)$ In particular, $\tau_c$ is periodic. More generally, the kernel of the homomorphism 
$\MM (\Sigma, \PP) \to \Aut( \CC (\Sigma, \PP))$ has the following geometric interpretation.

\bigskip\noindent
Let $c_1, \dots, c_q$ be the boundary components of $\Sigma$. Let $\D= \{ z \in \C; |z| \le 1 \}$ 
denote the standard disk. For each $1 \le i\le q$ we take a copy $\D_i$ of $\D$ and we denote by $d_i: 
\S^1 \to \D_i$ the boundary component of $\D_i$. Set $\Sigma_0= (\Sigma \sqcup (\sqcup_{1 \le i\le q} 
\D_i)) / \sim$, where $\sim$ is the equivalence relation which identifies $c_i(z)$ with $d_i(\bar z) 
= d_i^{-1}(z)$ 
for all $z \in \S^1$. In other words, $\Sigma_0$ is the surface obtained from $\Sigma$ gluing a disk 
along each boundary component (see Figure 2.1). Choose a point $Q_i$ in the interior of $\D_i$ (say 
$Q_i =0$) for all $1 \le i\le q$ and set $\QQ= \{ Q_1, \dots, Q_q\}$. Then the embedding $(\Sigma, \PP) 
\hookrightarrow (\Sigma_0, \PP \sqcup \QQ)$ determines a homomorphism $\theta: \MM( \Sigma, \PP) \to 
\MM( \Sigma_0, \PP \sqcup \QQ)$ that we call the {\it corking} of $(\Sigma, \PP)$. Note that $\theta$ 
is not surjective in general. Its image, $\Im \theta$, is the subgroup of $f \in \MM( \Sigma_0, \PP 
\sqcup \QQ)$ such that $f(\PP)= \PP$ and $f(Q_i)=Q_i$ for all $1 \le i\le q$. In particular, $\Im 
\theta$ has finite index in $\MM (\Sigma_0, \PP \sqcup \QQ)$.

\begin{figure}[htbp]
\centerline{
\setlength{\unitlength}{.4cm}
\begin{picture}(18,8)
\put(1,1){\includegraphics[width=6cm]{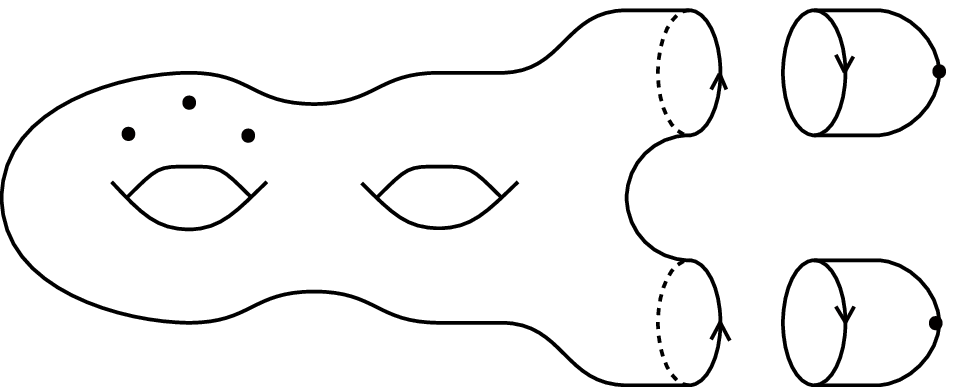}}
\put(2.5,6.1){\small $P_1$}
\put(3.5,6.5){\small $P_2$}
\put(4.5,6){\small $P_3$}
\put(11.5,2){\small $c_2$}
\put(11.5,6){\small $c_1$}
\put(13.5,2){\small $d_2$}
\put(13.5,6){\small $d_1$}
\put(16,1.5){\small $Q_2$}
\put(16,5.5){\small $Q_1$}
\end{picture}}
\centerline{{\bf Figure 2.1.} A corking.}
\end{figure}

\bigskip\noindent
The following result can be found for instance in \cite{ParRol1}.

\bigskip\noindent
{\bf Proposition 2.1.} {\it Assume $|\PP| \ge 2$ if $\Sigma= \D$ is the disk and $|\PP| \ge 1$ if 
$\Sigma$ is the annulus. Let $K$ be the subgroup of $\MM (\Sigma, \PP)$ generated by 
the Dhen twists $\tau_{c_1}, \dots, \tau_{c_q}$ along the boundary components of $\Sigma$.
Then $K= \Ker \theta$, it is a free abelian group freely generated by $\{ \tau_{c_1}, \dots, 
\tau_{c_q}\}$, and it is contained in the center of $\MM (\Sigma, \PP)$.}

\bigskip\noindent
Now, observe that $\CC (\Sigma, \PP)= \CC (\Sigma_0, \PP \sqcup \QQ)$ and that every $f \in \MM( 
\Sigma, \PP)$ acts on $\CC (\Sigma, \PP)$ in the same way as $\theta(f)$. Moreover, $\MM( \Sigma_0, \PP 
\sqcup \QQ)$ acts faithfully on $\CC (\Sigma, \PP) = \CC(\Sigma_0,\PP \sqcup \QQ)$, 
thus the kernel of $\MM (\Sigma, \PP) \to \Aut( 
\CC( \Sigma, \PP))$ is precisely $K$. Observe also that $f$ is pseudo-Anosov if and only if $\theta 
(f)$ is pseudo-Anosov, $f$ is periodic if and only if $\theta(f)$ is periodic, and $f$ is reducible if 
and only if $\theta(f)$ is reducible.


\subsection{Transverse singular foliations and pseudo-Anosov diffeomorphisms}

Throughout this subsection the surface $\Sigma$ is assumed to be closed (namely, $\partial \Sigma = 
\emptyset$). We shall also assume $|\PP| \ge 3$ if $\Sigma$ is the sphere and $|\PP| \ge 1$ if $\Sigma$ 
is the torus, so that the Euler characteristic of $\Sigma \setminus \PP$ is $<0$. Note that this last 
restriction is not important since $\MM (\S^2)= \MM (\S^2, 1)= \MM( \S^2,2)= \{ 1\}$ and $\MM 
(\T^2)= \MM (\T^2,1)$ (see \cite{Birma1}, Chapter 4). The basic reference for the material explained in 
this subsection is \cite{FaLaPo1}. Other references are \cite{BleCas1}, \cite{Thurs1}, and 
\cite{Ivano2}.

\bigskip\noindent
{\bf Transverse singular foliations.}
A pair of {\it transverse singular foliations} $(\FF_1, \FF_2)$ on $\Sigma$ is a finite set $\SS$ of 
{\it singularities} together with a covering of $\Sigma \setminus \SS$ by an atlas of charts $\{ 
\varphi_i: U_i \to ]-1,1[^2 \}_{i \in I}$ such that the changes of charts preserve the product 
structure, that is, for all $i,j \in I$ and all $x \in U_i \cap U_j$, there exists an open neighborhood 
$U \subset U_i \cap U_j$ of $x$ on which the change of coordinates $\psi_{i\,j}= \varphi_j \circ 
\varphi_i^{-1} : \varphi_i(U) \to \varphi_j(U)$ is of the form $\psi_{i\,j} (t_1,t_2)= 
(\psi_{i\,j}^{(1)}(t_1), \psi_{i\, j}^{(2)}(t_2))$. The {\it leaves} of $\FF_1$ (resp. 
$\FF_2$) are the classes of the smallest equivalence relation $\approx_1$ (resp. $\approx_2$) on 
$\Sigma \setminus \SS$ such that $x \approx_1 y$ (resp. $x \approx_2 y$) whenever $x,y \in U_i$ for 
some $i \in I$, and $\varphi_i(x)$ and $\varphi_i(y)$ belong to the same horizontal (resp. vertical) 
segment (namely, $\varphi_i(x)$ and $\varphi_i(y)$ have the same second (resp. first) coordinate).

\bigskip\noindent
From now on, we shall always assume the singularities to be of the following sort.

\bigskip\noindent
Let $P \in \SS$ be a singularity and let $k$ be an even positive number (say $k=2k'$). We say that $P$ 
is a {\it singularity with $k$ prongs} (see Figure 2.2) if there exist an open disk $U$ embedded in 
$\Sigma$ and centered at $P$ and a $k'$-folds ramified cover $\pi: U \to ]-1, 1[^2$ such that
\begin{itemize}
\item
$\pi(P)= (0,0)$, and $(0,0)$ is the unique ramification of $\pi$,
\item
if $\LL_1$ (resp. $\LL_2$) is a leaf of $\FF_1$ (resp. $\FF_2$), then the image of a connected component 
of $\LL_1 \cap U$ (resp. $\LL_2 \cap U$) is contained in a horizontal (resp. vertical) segment of $]-
1,1[^2$.
\end{itemize}
Note that, if $k=2$, then $\pi$ is a regular chart and $P$ may be viewed as a regular point. The {\it 
prongs} of $\FF_1$ (resp. $\FF_2$) at $P$ are defined to be the connected components of $\pi^{-1} ((]-
1,0[ \cup ]0,1[) \times \{0\})$ (resp. $\pi^{-1} (\{ 0\} \times (]-1,0[ \cup ]0,1[))$.

\bigskip\noindent
Let $P \in \SS$ be a singularity and let $k$ be an odd positive number. We say that $P$ is a {\it 
singularity with $k$ prongs} (see Figure 2.2) if there exists an open disk $U$ embedded in $\Sigma$ and 
centered at $P$ such that: if $\tilde \pi: \tilde U \to U$ is the 2-folds ramified cover of $U$ with 
unique ramification at $P$, $\tilde \pi^{-1}(P)= \tilde P$, and $\tilde \FF_1$ (resp. $\tilde \FF_2$) 
is the lift of $\FF_1$ (resp. $\FF_2$) in $\tilde U$, then $\tilde P$ is a singularity with $2k$ prongs 
of $(\tilde \FF_1, \tilde \FF_2)$. The {\it prongs} of $\FF_1$ (resp. $\FF_2$) at $P$ are defined to be 
the images under $\tilde \pi$ of the prongs of $\tilde \FF_1$ (resp. $\tilde \FF_2$) at $\tilde P$.

\begin{figure}[htbp]
\centerline{
\setlength{\unitlength}{.4cm}
\begin{picture}(34,12)
\put(1,3){\includegraphics[width=12.8cm]{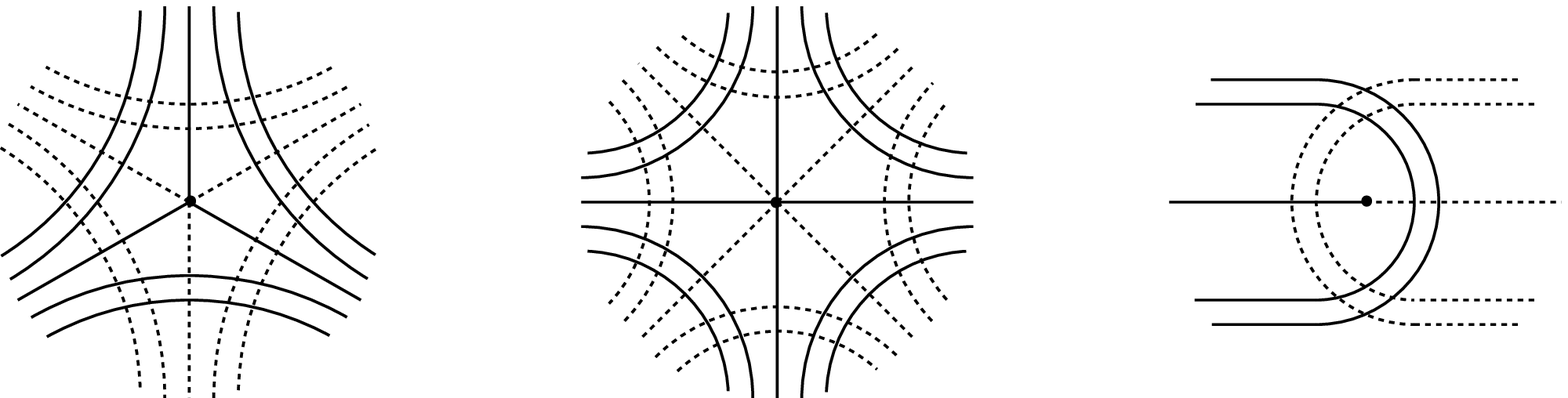}}
\put(3,1){\small $3$ prongs}
\put(15,1){\small $4$ prongs}
\put(27,1){\small $1$ prong}
\end{picture}}
\centerline{{\bf Figure 2.2.} Singularities with 1, 3, and 4 prongs.}
\end{figure}

\bigskip\noindent
A leaf of $\FF_i$ containing a prong at $P \in \SS$ is called a {\it separatrix} of $\FF_i$ at $P$.

\bigskip\noindent
From now on, we always assume the set $\PP$ of punctures to be included in $\SS$. However, a puncture 
$P \in \PP$ may be a singularity with 2 prongs.

\bigskip\noindent
The {\it index} of $(\FF_1, \FF_2)$ at a singularity $P$ is defined to be  
the number of prongs of $\FF_i$ at $P$ and is denoted by $\Ind (\FF_1, \FF_2: P)$. The 
indices and the Euler characteristic of $\Sigma$, denoted by $\chi(\Sigma)$, are related by the following 
formula (see \cite{FaLaPo1}, Expos\'e 5).

\bigskip\noindent
{\bf Proposition 2.2.}
\[
\chi(\Sigma)= \sum_{P \in \SS} \left( 1-{ \Ind( \FF_1, \FF_2: P) \over 2} \right)\,.
\]
\qed

\bigskip\noindent
{\bf Transverse invariant measures.}
A {\it transverse holonomy invariant measure} $\mu_1$ of $\FF_1$ is a function from the set of 
(non-necessarily connected) compact curves $\gamma \subset \Sigma$, with values on $[0,+\infty[$, and such 
that:
\begin{itemize}
\item
$\mu_1(\gamma_1 \cup \gamma_2)= \mu_1(\gamma_1) + \mu_1(\gamma_2)$ for all compact curves $\gamma_1, 
\gamma_2 \subset \Sigma$,
\item
if $\gamma_1, \gamma_2$ are both segments transverse to $\FF_1$, contained in the source $U$ of a 
regular chart $\varphi: U \to ]-1,1[^2$, and having their extremities on the same two horizontal 
segments of the chart, then $\mu_1(\gamma_1)= \mu_1(\gamma_2)$.
\end{itemize}
We say that $\mu_1$ has {\it total support} if $\mu_1(\gamma) >0$ for every segment $\gamma: [0,1] \to 
\Sigma\setminus \SS$ embedded in $\Sigma$ and (topologically) transverse to $\FF_1$. We say that the 
pair of transverse singular foliations $(\FF_1, \FF_2)$ is {\it measured} if each $\FF_i$ is equipped 
with a transverse holonomy invariant measure $\mu_i$ with total support.

\bigskip\noindent
If $\mu$ is a transverse holonomy invariant measure of a foliation $\FF$ and $\lambda >0$, then 
$\lambda \mu$ is also a transverse holonomy invariant measure of $\FF$. Clearly, a transverse holonomy 
invariant measure is not bounded, so there is no way to normalize it, but one may speak about {\it 
classes of measures} which are defined to be the orbits of the transverse measures under the natural 
action of $]0,+\infty[$.
Any foliation of a compact surface admits a transverse holonomy invariant 
measure (which may be with non-total support). A foliation is called {\it uniquely ergodic} if it admits 
exactly one transverse holonomy invariant measure with total support, up to a positive multiplier.

\bigskip\noindent
A diffeomorphism $F: \Sigma \to \Sigma$ acts naturally on the set of pairs of transverse foliations: the 
pair $(F(\FF_1), F(\FF_2))$ is defined by the charts $\varphi_i \circ F^{-1}: F(U_i) \to ]-1,1[^2$, $i \in 
I$, where $\{ \varphi_i: U_i \to ]-1,1[^2\}_{i\in I}$ is the defining atlas of $(\FF_1, \FF_2)$. 
Moreover, if $(\FF_1, \FF_2)$ is measured, then $(F(\FF_1), F(\FF_2))$ is also measured: the transverse 
holonomy invariant measure $F_\ast (\mu_i)$ is defined by $F_\ast(\mu_i)(\gamma)= \mu_i(F^{-1} 
(\gamma))$.

\bigskip\noindent
{\bf Pseudo-Anosov diffeomorphisms.} A diffeomorphism $F \in \Diff (\Sigma, \PP)$ is called {\it pseudo-
Anosov} if there exist a pair $(\FF^s, \FF^u)$ of transverse measured singular foliations 
with measures $(\mu^s, \mu^u)$ and set of singularities $\SS$, and a real number $\lambda >1$, such that:
\begin{itemize}
\item
$F(\SS)=\SS$, $F(\FF^s)= \FF^s$, $F(\FF^u)= \FF^u$, $F_\ast(\mu^s)= {1 \over \lambda} \mu^s$, and 
$F_\ast(\mu^u)= \lambda \mu^u$,
\item
the 1-prongs singularities of $(\FF^s, \FF^u)$ belong to $\PP$.
\end{itemize}
We emphasize that a $2$-prongs singularity is not
a regular point and must be sent to a 2-prongs singularity through a pseudo-Anosov
diffeomorphism.

\bigskip\noindent
{\bf Remark.} As pointed out in \cite{FaLaPo1}, Expos\'e 9, a pseudo-Anosov diffeomorphism is in fact a 
diffeomorphism on $\Sigma \setminus \SS$ but is not $C^1$ at the singularities. However, there are 
canonical local models for the actions of the pseudo-Anosov diffeomorphisms at the singularities, and 
this abuse, which is common in the theory, do not alterate the results stated in the remainder.

\bigskip\noindent
The following proposition collects some classical results on pseudo-Anosov diffeomorphisms.

\bigskip\noindent
{\bf Proposition 2.3.} {\it Let $F \in \Diff (\Sigma, \PP)$ be a pseudo-Anosov diffeomorphism, and let 
$\FF^s$, $\FF^u$, $\mu^s$, $\mu^u$, and $\lambda$ as above.
\begin{enumerate}
\item
$\FF^s$ and $\FF^u$ are the unique singular foliations that are invariant 
under $F$.
\item
{\rm (See \cite{FaLaPo1}, Expos\'e 12)} $\FF^s$ and $\FF^u$ are uniquely ergodic.
\item
{\rm (See \cite{FaLaPo1}, Expos\'e 9)} Each separatrix of $\FF^s$ (resp. $\FF^u$) contains exactly one 
prong.
\item
{\rm (See \cite{FaLaPo1}, Expos\'e 9)} Each leaf of either $\FF^s$ or $\FF^u$ is dense in $\Sigma$. 
More precisely, each half-leaf (connected component of a leaf minus a point) without prong is dense.
\end{enumerate}}
\qed

\bigskip\noindent
The above foliations $\FF^s$ and $\FF^u$
associated to a pseudo-Anosov diffeomorphism $F$ are called, respectively, the {\it stable} and {\it unstable 
foliations} of $F$. The number $\lambda>1$ is uniquely determined by $F$ and is called its {\it 
dilatation coefficient}.

\bigskip\noindent
Now, the coherence of the terminologies introduced in the present subsection and in the previous one is 
ensured by the following theorem (see \cite{FaLaPo1}, Expos\'e 9).

\bigskip\noindent
{\bf Theorem 2.4.} {\it The isotopy class of a pseudo-Anosov diffeomorphism is a pseudo-Anosov element 
of $\MM (\Sigma, \PP)$. Conversely, any pseudo-Anosov element of $\MM(\Sigma, \PP)$ is the 
isotopy class of a pseudo-Anosov diffeomorphism.}

\subsection{Normalizers and symmetries of a pair of measured foliations}

We assume again in this subsection that $\Sigma$ is closed, $|\PP| \ge 3$ if $\Sigma$ is the sphere, and 
$|\PP|\ge 1$ if $\Sigma$ is the torus.

\bigskip\noindent
We fix a pseudo-Anosov diffeomorphism $F \in \Diff (\Sigma, \PP)$ and we denote by $\FF^s$ (resp. $\FF^u$) 
the stable (resp. unstable) foliation of $F$, by $\SS$ the set of singularities, by 
$(\mu^s, \mu^u)$ the measures of $(\FF^s, \FF^u)$, and by $\lambda$ the dilatation coefficient of $F$. 

\bigskip\noindent
{\bf Normalizer of $(\FF^s, \FF^u)$.}
Define the {\it normalizer} of the pair $(\FF^s, \FF^u)$ to be the group 
\linebreak
$\Norm(\FF^s, \FF^u)$ of $G 
\in \Diff (\Sigma, \PP)$ such that $G (\FF^s) = \FF^s$ and $G (\FF^u) = \FF^u$. Define the {\it 
extended normalizer} of $(\FF^s, \FF^u)$ to be the group $\NNorm (\FF^s, \FF^u)$ of $G \in \Diff 
(\Sigma, \PP)$ such that $G (\{ \FF^s, \FF^u \}) = \{ \FF^s, \FF^u \}$.

\bigskip\noindent
The centralizer of $F$ is contained in $\Norm(\FF^s, \FF^u)$. Indeed, if $G$ commutes with $F$, then 
$G(\FF^s)$ and $G(\FF^u)$ are, respectively, the stable and unstable foliations of $G \circ F \circ G^{-
1} = F$, thus $G(\FF^s)= \FF^s$ and $G(\FF^u)= \FF^u$. Note also that either $\NNorm (\FF^s, \FF^u) 
= \Norm (\FF^s, \FF^u)$ or $\Norm( \FF^s, \FF^u)$ is an index 2 subgroup of $\NNorm (\FF^s, \FF^u)$.

\bigskip\noindent
Let $G \in \Norm( \FF^s, \FF^u)$. Consider the measure $G_\ast(\mu^u)$ on $\FF^u$. This is a transverse 
invariant measure on $\FF^u$, thus, by the unique ergodicity of $\FF^u$, there is some number 
$\lambda^u(G)>0$ such that $G_\ast (\mu^u)= \lambda^u(G)\mu^u$. Similarly, there is some number 
$\lambda^s(G)>0$ such that $G_\ast(\mu^s)= \lambda^s(G)\mu^s$. Note also that the map $L: \Norm( \FF^s, 
\FF^u) \to \R$, $G \mapsto \log (\lambda^u(G))$ is a group homomorphism. The proof of the following 
lemma can be found for instance in \cite{FaLaPo1}, Expos\'e 9.

\bigskip\noindent
{\bf Lemma 2.5.} {\it Let $G \in \Norm( \FF^s, \FF^u)$. Then $\lambda^s(G) \cdot \lambda^u(G)=1$.}
\qed

\bigskip\noindent
In particular, if $\lambda^u(G)>1$, then $G$ is a pseudo-Anosov diffeomorphism, $\FF^s$ (resp. $\FF^u$) 
is the stable (resp. unstable) foliation of $G$, and $\lambda^u(G)$ is its dilatation coefficient. If 
$\lambda^u(G)<1$, then $G$ is a pseudo-Anosov diffeomorphism, $\FF^u$ (resp. $\FF^s$) is the stable 
(resp. unstable) foliation of $G$, and $\lambda^s(G)= \lambda^u(G)^{-1}$ is its dilatation coefficient.

\bigskip\noindent
We refer to \cite{ArnYoc1} for the proof of the next lemma.

\bigskip\noindent
{\bf Lemma 2.6.} {\it The image $L( \Norm( \FF^s, \FF^u))$ is a discrete subgroup of $\R$, hence it is of 
the form ${\log(\lambda) \over k} \cdot \Z$, where $k \in \N \setminus \{ 0\}$.}
\qed

\bigskip\noindent
{\bf Symmetries of $(\FF^s, \FF^u)$.}
Define the {\it symmetries} of $(\FF^s, \FF^u)$ to be the elements of $\Sym( \FF^s, \FF^u)$ 
$=\Ker (L)$, and define the {\it antisymmetries} to be the elements of $\Sym_- (\FF^s, \FF^u)= 
\linebreak
\NNorm(\FF^s, \FF^u) \setminus \Norm (\FF^s, \FF^u)$. 

\bigskip\noindent
Note that, being a kernel, $\Sym (\FF^s, \FF^u)$ is a normal subgroup of $\Norm (\FF^s, \FF^u)$. More 
generally, it is easily checked that:

\bigskip\noindent
{\bf Lemma 2.7.} {\it $\Sym(\FF^s, \FF^u)$ is a normal subgroup of $\NNorm(\FF^s, \FF^u)$.}
\qed

\bigskip\noindent
Now, the coherence of the above terminology is given by the following.

\bigskip\noindent
{\bf Lemma 2.8.} {\it Let $G \in \Sym_-( \FF^s, \FF^u)$. Then $G^2 \in \Sym(\FF^s, \FF^u)$.}

\bigskip\noindent
{\bf Proof.} Let $\lambda_1, \lambda_2 >0$ be such that $G_\ast (\mu^s) = \lambda_1 \mu^u$ and 
$G_\ast(\mu^u)= \lambda_2 \mu^s$. Hence $G_\ast^2 (\mu^s) = \lambda_1 \lambda_2\mu^s$ and $G_\ast^2 
(\mu^u) = \lambda_1 \lambda_2 \mu^u$, that is, $\lambda^s(G^2) = \lambda^u(G^2) = \lambda_1 
\lambda_2$. Since $\lambda^s(G^2) \cdot \lambda^u(G^2) = 1$ (see Lemma 2.5), we conclude that 
$\lambda^s(G^2)= \lambda^u(G^2)=1$.
\qed

\bigskip\noindent
Note that we do not have necessarily $G_1 \circ G_2 \in \Sym( \FF^s, \FF^u)$ if $G_1, G_2 \in \Sym_-
(\FF^s, \FF^u)$, thus $\Sym(\FF^s, \FF^u) \sqcup \Sym_-(\FF^s, \FF^u)$ is not a group in general.

\bigskip\noindent
{\bf Lemma 2.9.} {\it Let $G$ be an element of $\Sym(\FF^s, \FF^u)$ which fixes some separatrix. Then 
$G$ is the identity map.}

\bigskip\noindent
{\bf Proof.} Let $\LL$ be a stable separatrix of a singular point $P$ which we assume to be fixed under 
$G$. In particular, we have $G(P)=P$. Observe that $x \mapsto \int_p^x d\mu^u$ induces a strictly 
increasing function from $\LL$ (oriented by the origin at $P$) to $[0, +\infty[$. Since $\lambda^u(G)=1$, 
the function $G$ preserves this parametrization of $\LL$, thus the restriction of $G$ to $\LL$ is the 
identity map. We conclude using the fact that $\LL$ is dense in $\Sigma$ (see Proposition 2.3).
\qed

\bigskip\noindent
{\bf Corollary 2.10.} {\it The group $\Sym(\FF^s, 
\FF^u)$ acts freely on the set of separatrices of $\FF^s$. In particular, $\Sym(\FF^s, \FF^u)$ is 
finite.}

\bigskip\noindent
{\bf Proof.} By the Euler-Poincar\'e formula given in Proposition 2.2, the set $\SS$ cannot be empty if 
$\Sigma$ is not a torus. On the other hand, we have assumed that $\PP \subset \SS$, and $\PP \neq 
\emptyset$ if $\Sigma$ is a torus. So, the set of separatrices is a non-empty finite set and 
$\Sym(\FF^s, \FF^u)$ acts freely on it by Lemma 2.9, thus $\Sym(\FF^s, \FF^u)$ is finite.
\qed 

\bigskip\noindent
{\bf Lemma 2.11.} {\it Let $G$ be an element of $\Norm (\FF^s, \FF^u)$ 
which is isotopic to the identity map. Then $G$ is equal to the identity map.}

\bigskip\noindent
{\bf Proof.} First, observe that $\lambda^s(G)= \lambda^u(G)=1$ (that is, $G \in \Sym(\FF^s, 
\FF^u)$), otherwise $G$ would be a pseudo-Anosov diffeomorphism which cannot be isotopic to the 
identity by Theorem 2.4. By Corollary 2.10, $G$ has finite order, and, by \cite{FaLaPo1}, Expos\'e 12, 
any finite order diffeomorphism isotopic to the identity must be the identity map.
\qed

\bigskip\noindent
{\bf Corollary 2.12.} {\it The restriction of the natural homomorphism 
$\Diff (\Sigma, \PP) \to \MM (\Sigma, \PP)$ to $\Norm (\FF^s, \FF^u)$ is injective.}
\qed


\subsection{From pseudo-Anosov mapping classes to Pseudo-Anosov diffeomorphisms}

We still assume in this section that $\Sigma$ is a closed surface, that $|\PP| \ge 3$ if $\Sigma$ is 
the sphere, and that $|\PP| \ge 1$ if $\Sigma$ is the torus. The following theorem (and its 
corollaries) is the key of many proofs in the present paper. Its proof
can be found in \cite{FaLaPo1}, Expos\'e 12 (see also \cite{Hande1}).

\bigskip\noindent
{\bf Theorem 2.13.} {\it Let $F, G \in \Diff (\Sigma, \PP)$ be two pseudo-Anosov diffeomorphisms. If $F$ 
is isotopic to $G$, then there exists $H \in \Diff (\Sigma, \PP)$ isotopic to the identity map and such 
that $G= H \circ F \circ H^{-1}$.}
\qed

\bigskip\noindent
{\bf Corollary 2.14.} {\it Let $F \in \Diff (\Sigma, \PP)$ be a pseudo-Anosov diffeomorphism and $f \in 
\MM( \Sigma, \PP)$ be its isotopy class. Let $g \in \MM (\Sigma, \PP)$ and $m \ge 1$ such that 
$f^m=g^m$. Then there is a pseudo-Anosov representative $G \in \Diff (\Sigma, \PP)$ of $g$ such that 
$F^m=G^m$.}

\bigskip\noindent
{\bf Proof.} Clearly, the equality $f^m=g^m$ implies that $g$ has no finite orbit in $\CC(\Sigma, 
\PP)$, that is, $g$ is a pseudo-Anosov element of $\MM (\Sigma, \PP)$. Let $G' \in \Diff(\Sigma, \PP)$ be 
a pseudo-Anosov representative of $g$. The diffeomorphisms $F^m$ and $G'^m$ are both pseudo-Anosov, and they 
are isotopic, thus, by Theorem 2.13, there exists $H \in \Diff (\Sigma, \PP)$ isotopic to the identity 
map such that $H \circ G'^m \circ H^{-1} = F^m$. Set $G= H \circ G' \circ H^{-1}$. Then $G$ is a 
pseudo-Anosov diffeomorphism, its isotopy class is $g$, and $G^m=F^m$.
\qed

\bigskip\noindent
{\bf Corollary 2.15.} {\it Let $F \in \Diff (\Sigma, \PP)$ be a pseudo-Anosov diffeomorphism and $f \in 
\MM(\Sigma, \PP)$ be its isotopy class. Let $g \in \MM(\Sigma, \PP)$ be a pseudo-Anosov element which 
commutes with $f$. Then there is a pseudo-Anosov representative $G \in \Diff(\Sigma, \PP)$ of $g$ such 
that $F \circ G = G \circ F$.}

\bigskip\noindent
{\bf Proof.} Let $G' \in \Diff (\Sigma, \PP)$ be a pseudo-Anosov representative of $g$. The diffeomorphism 
$G' \circ F \circ G'^{-1}$ is a pseudo-Anosov diffeomorphism isotopic to $F$, thus, by Theorem 2.13, 
there exists $H \in \Diff(\Sigma, \PP)$ isotopic to the identity map such that $H \circ G' \circ F 
\circ G'^{-1} \circ H^{-1} = F$. Set $G= H \circ G'$. Then $G \circ F = F \circ G$ and $G$ is isotopic 
to $G'$, namely, $G$ is a representative of $g$. So, it remains to prove that $G$ is a pseudo-Anosov 
diffeomorphism.

\bigskip\noindent
Let $\FF^s$ (resp. $\FF^u$) be the stable (resp. unstable) foliation of $F$. As observed before, the 
fact that $G$ commutes with $F$ implies that it belongs to $\Norm (\FF^s, \FF^u)$. Moreover, the isotopy 
class of $G$ is $g$ which is assumed to be pseudo-Anosov, thus $\lambda^s(G) \neq 1$. So, $G$ is a 
pseudo-Anosov diffeomorphism.
\qed 

\subsection{Coverings and pseudo-Anosov diffeomorphisms}

We turn now in this subsection to explain a way for constructing pseudo-Anosov diffeomorphisms having 
precise groups of symmetries by means of ramified coverings. Then the proof of Theorem 5.2 will be a 
direct application of the present subsection together with the next one.

\bigskip\noindent
{\bf Lemma 2.16.} {\it Let $M$ be a compact connected manifold, let $\pi: N \to M$ be a finite covering, 
possibly with a finite set $\RR$ of ramification points, and let $F: M \to M$ be a homeomorphism such 
that $F (\RR) = \RR$. Then there exists $k \ge 1$ such that $F^k$ has a lift on $N$.}

\bigskip\noindent
{\bf Proof.} We choose a base point $P_0 \in M \setminus \RR$ and a path $\gamma: [0,1] \to M \setminus 
\RR$ joining $P_0$ to $F(P_0)$, and we denote by $\phi_\gamma : \pi_1( M \setminus \RR, F(P_0)) \to 
\pi_1( M \setminus \RR, P_0)$ the isomorphism defined by $\alpha \mapsto \gamma^{-1} \alpha \gamma$. 
The general theory of coverings asserts that $F|_{M \setminus \RR}$ has a lift on $N \setminus 
\pi^{-1} (\RR)$ if and only if 
\[
(\phi_\gamma \circ F_\ast) (\pi_1 (N \setminus \pi^{-1} (\RR), Q_0)) = \pi_1 (N \setminus \pi^{-1} 
(\RR), Q_0)\,,
\]
where $Q_0$ is a chosen element of the preimage of $P_0$,
and $\pi_1 (N \setminus \pi^{-1} (\RR), Q_0)$ is viewed as a subgroup of $\pi_1(M \setminus \RR, P_0)$. 
Let $m$ be the number of sheets of the 
covering. The group $\pi_1( N \setminus \pi^{-1} (\RR), Q_0)$ is an index $m$ subgroup of $\pi_1 (M 
\setminus \RR, P_0)$ and there are finitely many index $m$ subgroups in $\pi_1 (M \setminus \RR, P_0)$.
Thus, there exists some $k \ge 1$ such that
\[
(\phi_\gamma \circ F_\ast)^k (\pi_1( N \setminus \pi^{-1} (\RR), Q_0)) = \pi_1(N \setminus \pi^{-1}
(\RR), Q_0) \,.
\]
On the other hand, we have $(\phi_\gamma \circ F_\ast)^k = \phi_{\gamma^{(k)}} \circ (F^k)_\ast$, where
\[
\gamma^{(k)} = (F^{k-1} \circ \gamma) \cdot \dots \cdot (F \circ \gamma) \cdot \gamma\,,
\]
which is a path from $P_0$ to $F^k(P_0)$.
So, $F^k|_{M \setminus \RR}$ has a lift $\tilde F: N \setminus \pi^{-
1}(\RR) \to N \setminus \pi^{-1} (\RR)$. Clearly, $\tilde F$ can be extended to a homeomorphism 
$\tilde F : N \to N$, and $\tilde F$ is a lift of $F^k: M \to M$.
\qed

\bigskip\noindent
We consider a regular ramified covering $\pi: \Sigma \to \Sigma_0$ of a closed surface $\Sigma_0$. That 
is, there is a finite group $\Gamma$, called the {\it Galois group} of $\pi$, acting on $\Sigma$, and such 
that $\Sigma_0= \Sigma / \Gamma$ and $\pi: \Sigma \to \Sigma_0$ is the natural quotient. Let $\RR_0$ be 
the set of ramification points of $\pi$. For each $Q \in \RR_0$, we denote by $r_\pi(Q)$ the {\it 
ramification order} of $\pi$ at $Q_0$. If $o(Q)= | \pi^{-1} (Q)|$, then this number is defined by 
$r_\pi(Q)= {m \over o(Q)}$, where $m$ is the number of sheets of $\pi$. Here we prescribe that $r_\pi(Q) \ge 
2$ for all $Q \in \RR_0$, so that there is no regular point hidden in $\RR_0$. If $Q$ is regular, that 
is, if $Q \in \Sigma_0 \setminus \PP_0$, then the {\it ramification order} of $\pi$ at $Q$ is $r_\pi(Q) 
=1$.

\bigskip\noindent
Let $F_0 \in \Diff(\Sigma_0, \RR_0)$ be a pseudo-Anosov diffeomorphism, let $\FF_0^s$ (resp. $\FF_0^u$) 
be the stable (resp. unstable) foliation of $F_0$, and let $\lambda$ be its dilatation coefficient. 
Then the following lemma can be easily proved with standard arguments.

\bigskip\noindent
{\bf Lemma 2.17.} {\it The foliation $\FF_0^s$ (resp. $\FF_0^u$) lifts to a foliation $\FF^s$ (resp. 
$\FF^u$) on $\Sigma$, and $(\FF^s, \FF^u)$ is a pair of transverse measured singular foliations of $\Sigma$. 
Moreover, we have
\[
\Ind( \FF^s, \FF^u : P) = \Ind (\FF_0^s, \FF_0^u : \pi(P)) \cdot r_\pi (\pi(P))
\]
for all $P \in \Sigma$, and $\Gamma$ is a subgroup of $\Sym( \FF^s, \FF^u)$.}
\qed

\bigskip\noindent
We keep the notations and we assume, furthermore, that $F_0$ lifts to a diffeomorphism $F: \Sigma \to 
\Sigma$. (By Lemma 2.16, this is always possible up to some power.) Then $F$ is a pseudo-Anosov 
diffeomorphism, $\FF^s$ (resp. $\FF^u$) is its stable (resp. unstable) foliation, and $\lambda$ is its 
dilatation coefficient.

\bigskip\noindent
{\bf Pivot.}
A point $Q_0 \in \RR_0$ is called a {\it pivot} of $(F_0, \pi)$ if 
it is a singular point with one prong (i.e.
$\Ind (\FF_0^s, \FF_0^u:Q_0)=1$), and 
if $\Ind( \FF_0^s, \FF_0^u: Q) \cdot r_\pi (Q) \neq r_\pi(Q_0)$ for all $Q \in \Sigma_0 \setminus \{ Q_0 
\}$.

\bigskip\noindent
{\bf Proposition 2.18.} {\it If $(F_0, \pi)$ has a pivot, then $\Sym (\FF^s, \FF^u)$ is equal to
the Galois group $\Gamma$.}

\bigskip\noindent
{\bf Proof.} We already know that $\Gamma \subset \Sym (\FF^s, \FF^u)$ (see Lemma 2.17). For $k \ge 
3$, we denote by $\SS_k$ the set of singularities of $(\FF^s, \FF^u)$ with $k$ prongs, and by 
$\SS\SS_k$ the set of separatrices of $\FF^s$ at the elements of $\SS_k$. (So, by Proposition 2.3, we 
have $|\SS \SS_k| = k \cdot | \SS_k|$.) Take a pivot $Q_0 \in \RR_0$, and set $k_0= r_\pi(Q_0)$. By 
Lemma 2.17, the fact that $Q_0$ is a pivot implies that $\SS_{k_0} = \pi^{-1} (Q_0)$. Moreover, the 
group $\Gamma$ acts transitively on $\SS\SS_{k_0}$ because the elements of $\SS\SS_{k_0}$ are precisely 
the preimages under $\pi$ of the unique separatrix of $\FF_0^s$ at $Q_0$. We conclude by Lemma 2.9 that 
$\Gamma= \Sym (\FF^s, \FF^u)$.
\qed

\subsection{Pseudo-Anosov diffeomorphisms which fix their separatrices}

Throughout this subsection $\Sigma$ is a closed surface and $F: \Sigma \to \Sigma$ is a pseudo-Anosov 
diffeomorphism. We denote by $\FF^s$ (resp. $\FF^u$) the stable (resp. unstable) foliation of $F$.

\bigskip\noindent
{\bf Fixed separatrices.} We say that $F$ {\it fixes its separatrices} if it acts trivially on the set 
of separatrices of $\FF^s$ or, equivalently, if it acts trivially on the set of separatrices of 
$\FF^u$.

\bigskip\noindent
{\bf Lemma 2.19.} {\it Suppose that $F$ fixes its separatrices.
\begin{enumerate}
\item
The centralizer of $F$ in $\Diff (\Sigma)$ is equal to $\Norm (\FF^s, \FF^u)$.
\item
Let $G \in \NNorm (\FF^s, \FF^u) \setminus \Norm (\FF^s, \FF^u) = \Sym_- (\FF^s, \FF^u)$. Then $G 
\circ F \circ G = F^{-1}$.
\end{enumerate}}

\bigskip\noindent
{\bf Proof.} As pointed out in Subsection 2.3, we already know that the centralizer of $F$ is contained in $\Norm 
(\FF^s, \FF^u)$. So, it remains to show that $\Norm (\FF^s, \FF^u)$ is contained in the centralizer of 
$F$. Let $G \in \Norm (\FF^s, \FF^u)$. Then $G \circ F \circ G^{-1}$ is a pseudo-Anosov diffeomorphism, 
its stable (resp. unstable) foliation is $\FF^s$ (resp. $\FF^u$), its dilatation coefficient is 
$\lambda$, and it acts trivially on the set of separatrices of $\FF^s$. We conclude by Lemma 2.9 that 
$F^{-1} \circ G \circ F \circ G^{-1}$ is the identity map, thus $G \circ F \circ G^{-1} = F$.

\bigskip\noindent
Let $G \in \Sym_-( \FF^s, \FF^u)$. Then $G \circ F \circ G^{-1}$ is a pseudo-Anosov diffeomorphism, its 
stable (resp. unstable) foliation is $\FF^u$ (resp. $\FF^s$), its dilatation coefficient is $\lambda$, 
and it acts trivially on the set of separatrices of $\FF^s$. We conclude by Lemma 2.9 that $F \circ G 
\circ F \circ G^{-1}$ is the identity map, thus $G \circ F \circ G^{-1} = F^{-1}$.
\qed

\bigskip\noindent
{\bf Proposition 2.20.} {\it Suppose that $F$ fixes its separatrices.
Let $R_1, R_2 \in \Sym( \FF^s, \FF^u)$. Write $F_1= F \circ R_1$ (resp. $F_2= F 
\circ R_2$) and denote by $f_1 \in \MM(\Sigma)$ (resp. $f_2 \in \MM(\Sigma)$) the isotopy class of $F_1$ 
(resp. $F_2$). There exists $h \in \MM(\Sigma)$ such that $hf_1h^{-1} = f_2$ if and only if there 
exists $H \in \Norm (\FF^s, \FF^u)$ such that $H \circ R_1 \circ H^{-1} = R_2$.}

\bigskip\noindent
{\bf Proof.} Assume that there exists $H \in \Norm (\FF^s, \FF^u)$ such that $H \circ R_1 \circ H^{-1} = 
R_2$. By Lemma 2.19, we have $H \circ F \circ H^{-1} = F$, thus $H \circ F_1 \circ H^{-1} = F \circ (H 
\circ R_1 \circ H^{-1}) = F \circ R_2 = F_2$, therefore $h f_1 h^{-1} = f_2$, where $h \in \MM (\Sigma)$ 
is the isotopy class of $H$.

\bigskip\noindent
Now, assume that there exists $h \in \MM (\Sigma)$ such that $h f_1 h^{-1} = f_2$. Let $H_1 \in \Diff 
(\Sigma)$ which represents $h$. The diffeomorphisms $H_1 \circ F_1 \circ H_1^{-1}$ and $F_2$ are both 
pseudo-Anosov and contained in the same isotopy class, thus, by Theorem 2.13, there exists $H_2 \in 
\Diff (\Sigma)$ isotopic to the identity map such that $H_2 \circ H_1 \circ F_1 \circ H_1^{-1} \circ 
H_2^{-1} = F_2$. Set $H= H_2 \circ H_1$. Then $F_2= H \circ F_1 \circ H^{-1}$.

\bigskip\noindent
The foliation $\FF^s$ (resp. $\FF^u$) is the stable (resp. unstable) foliation of $F_2$, and $H(\FF^s)$ 
(resp. $H(\FF^u)$) is the stable (resp. unstable) foliation of $H \circ F_1 \circ H^{-1} = F_2$, thus 
$H(\FF^s)= \FF^s$ and $H (\FF^u)= \FF^u$, that is, $H \in \Norm (\FF^s, \FF^u)$. By Lemma 2.19, we have $H 
\circ F \circ H^{-1} = F$, thus $F \circ (H \circ R_1 \circ H^{-1}) = H \circ F_1 \circ H^{-1} = F_2 = 
F \circ R_2$, therefore $H \circ R_1 \circ H^{-1} = R_2$.
\qed

\bigskip\noindent
{\bf Proposition 2.21.} {\it Suppose that $F$ fixes its separatrices.
Let $R_1, R_2 \in \Sym (\FF^s, \FF^u)$. Write $F_1= F \circ R_1$ (resp. $F_2= F 
\circ R_2$) and denote by $f_1 \in \MM (\Sigma)$ (resp. $f_2 \in \MM (\Sigma)$) the isotopy class of 
$F_1$ (resp. $F_2$). There exists $h \in \MM (\Sigma)$ such that $hf_1h^{-1}$ commutes with $f_2$ if 
and only if there exists $H \in \NNorm (\FF^s, \FF^u)$ such that $H \circ R_1 \circ H^{-1}$ 
commutes with $R_2$.}

\bigskip\noindent
{\bf Proof.} Assume that there exists $H \in \NNorm (\FF^s, \FF^u)$ such that $H \circ R_1 \circ 
H^{-1}$ and $R_2$ commute. By Lemma 2.19, there exists $\varepsilon \in \{ \pm 1 \}$ such that $H \circ 
F \circ H^{-1} = F^\varepsilon$. Moreover, again by Lemma 2.19, $H \circ R_1 \circ H^{-1}$ and $R_2$ 
commute with $F$ (since they belong to $\Sym (\FF^s, \FF^u) \subset \Norm (\FF^s, \FF^u)$). We conclude 
that $H \circ F_1 \circ H^{-1} = F^\varepsilon \circ (H \circ R_1 \circ H^{-1})$ commutes with $F_2 = F 
\circ R_2$, thus $hf_1h^{-1}$ commutes with $f_2$, where $h \in \MM (\Sigma)$ is the isotopy class of 
$H$.

\bigskip\noindent
Now, assume that there exists $h \in \MM (\Sigma)$ such that $h f_1 h^{-1}$ and $f_2$ commute. Let $H_1 
\in \Diff (\Sigma)$ which represents $h$. By Corollary 2.15, there exists a pseudo-Anosov 
diffeomorphism $F_1' \in \Diff (\Sigma)$ which represents $h f_1 h^{-1}$ and which comùmutes with 
$F_2$. The diffeomorphisms $H_1 \circ F_1 \circ H_1^{-1}$ and $F_1'$ are both pseudo-Anosov and 
contained in the same isotopy class, thus, by Theorem 2.13, there exists $H_2 \in \Diff(\Sigma)$ 
isotopic to the identity map such that $F_1'= H_2 \circ H_1 \circ F_1 \circ H_1^{-1} \circ H_2^{-
1}$. Set $H = H_2 \circ H_1$. Then $H \circ F_1 \circ H^{-1}$ and $F_2$ commute.

\bigskip\noindent
The foliation $\FF^s$ (resp. $\FF^u$) is the stable (resp. unstable) foliation of $F_2$, and $H \circ 
F_1 \circ H^{-1}$ belongs to the centralizer of $F_2$, thus $(H \circ F_1 \circ H^{-1}) \in \Norm (\FF^s, 
\FF^u)$. On the other hand, we have $(H \circ F_1 \circ H^{-1}) \in \Norm (H (\FF^s), H (\FF^u))$, 
because $\FF^s$ (resp. $\FF^u$) is the stable (resp. unstable) foliation of $F_1$. By Proposition 2.3, 
it follows that 
$\{ H (\FF^s), H(\FF^u) \} = \{ \FF^s, \FF^u \}$, thus
$H \in \NNorm (\FF^s, \FF^u)$. By Lemma 2.19, there exists $\varepsilon \in \{ \pm 1 
\}$ such that $H \circ F \circ H^{-1}= F^\varepsilon$. Moreover, $H \circ R_1 \circ H^{-1}$ and $R_2$ 
commute with $F$ (by Lemma 2.19). Since $H \circ F_1 \circ H^{-1} = F^\varepsilon \circ (H \circ R_1 
\circ H^{-1})$ commutes with $F_2= F \circ R_2$, we conclude that $H \circ R_1 \circ H^{-1}$ commutes with 
$R_2$.
\qed

\subsection{Canonical reduction systems and pure elements}

\bigskip\noindent
Now, the surface $\Sigma$ is not assumed to be closed anymore. To each reducible element $f \in \MM 
(\Sigma, \PP)$ one can associate a well-defined simplex $\Delta (f)$ of $\CC(\Sigma, \PP)$ called the 
{\it canonical reduction system} for $f$. It happens that it is easier to first define $\Delta(f)$ for 
a certain kind of elements, called the {\it pure elements}, as we turn now to do.
On the other hand, the pure subgroups are the main 
object of Section 6 and here is a good place to introduce them. A standard reference for the material 
explained in this subsection is \cite{Ivano2}.

\bigskip\noindent
{\bf Pure elements and pure subgroups.}
An element $f \in \MM (\Sigma, \PP)$ is called {\it pure} if, for all essential curves $c: \S^1 \to 
\Sigma \setminus \PP$, either $f(c)$ is isotopic to $c$, or $f^m(c)$ is not isotopic to $c$ for any $m 
\in \Z \setminus \{0\}$. 
We emphasize that, by $f(c)$ is isotopic to $c$, we also mean that $f$ preserves the orientation 
of $c$ up to isotopy.
Note that a pure element has no finite orbit of order $\ge 2$ in $\CC (\Sigma, 
\PP)$, and, if it is not the identity, it has infinite order. The pseudo-Anosov elements and the Dehn 
twists are examples of pure elements. A subgroup $\Gamma \subset \MM (\Sigma, \PP)$ is called {\it pure} 
if all its elements are pure. Note that a pure subgroup is torsion-free.

\bigskip\noindent
{\bf Example.} Let $m \ge 3$, and let $\Gamma(m)$ denote the kernel of the natural homomorphism 
$\MM (\Sigma, \PP) \to \Aut( H_1( \Sigma \setminus \PP, \Z/m\Z))$. Then, by \cite{Ivano2},  
$\Gamma(m)$ is a pure subgroup. Note that $\Gamma(m)$ has finite index in $\MM(\Sigma, \PP)$ since 
$H_1( \Sigma \setminus \PP, \Z/m\Z)$ is finite. So, there are finite index pure subgroups and, in 
particular, there are finite index torsion-free subgroups in $\MM (\Sigma, \PP)$.

\bigskip\noindent
Let $\Delta= \{ \delta_0, \delta_1, \dots, \delta_p \}$ be a simplex of $\CC (\Sigma, \PP)$. We choose 
essential curves $d_0, d_1, \dots, d_p$ such that $\langle d_i \rangle = \delta_i$ for all $0 \le i \le 
p$, and $d_i \cap d_j = \emptyset$ for all $0 \le i \neq j \le p$. We denote by $\Sigma_\Delta$ the 
surface obtained from $\Sigma$ cutting along the $d_i$'s. Let $c_1, \dots, c_q$ be the boundary 
components of $\Sigma$. Then $\Sigma_\Delta$ is a non-necessarily connected compact surface, whose 
boundary components are $c_1, \dots, c_q$, $d_0^{(1)}, d_1^{(1)}, \dots, d_p^{(1)}$, $d_0^{(2)}, d_1^{(2)}, 
\dots, d_p^{(2)}$, and $\Sigma= \Sigma_\Delta/\sim$, where $\sim$ is the equivalence relation which 
identifies $d_i^{(1)} (z)$ with $d_i^{(2)}(\bar z) = (d_i^{(2)})^{-1}(z)$ 
for all $0 \le i\le p$ and all $z \in \S^1$. Let 
$\pi_\Delta: \Sigma_\Delta \to \Sigma$ denote the natural quotient. Then $d_i(z)= \pi_\Delta \circ 
d_i^{(1)}(z) = \pi_\Delta \circ (d_i^{(2)})^{-1}(z)$ for all $0\le i\le p$ and all $z \in \S^1$. Let 
$\Sigma_{\Delta\,1}, \dots, \Sigma_{\Delta\,l}$ be the connected components of $\Sigma_\Delta$, let 
$\PP_\Delta= \pi_\Delta^{-1} (\PP)$, and let $\PP_{\Delta\,k} = \PP_\Delta \cap \Sigma_{\Delta\,k}$ for 
all $1 \le k\le l$. Set
\[
\MM(\Sigma_\Delta, \PP_\Delta) = \MM( \Sigma_{\Delta\,1}, \PP_{\Delta\,1}) \times \dots \times \MM 
(\Sigma_{\Delta\,l}, \PP_{\Delta\,l})\,.
\]
Then the map $\pi_\Delta: \Sigma_\Delta \to \Sigma$ induces a homomorphism $\theta_\Delta: \MM
(\Sigma_\Delta, \PP_\Delta) \to \MM (\Sigma, \PP)$.

\bigskip\noindent
The following result is well-known and its proof is implicit in \cite{ParRol1}.

\bigskip\noindent
{\bf Lemma 2.22.} {\it For $0 \le i \le p$, 
set $\varphi_i = \tau_{d_i^{(1)}} \cdot {\tau_{d_i^{(2)}}}^{-1}$.
Let $K_\Delta$ be the subgroup of $\MM (\Sigma_\Delta, \PP_\Delta)$ generated by 
$\varphi_0, \varphi_1, \dots, \varphi_p$. Then $K_\Delta= \Ker \theta_\Delta$, it is a free abelian 
group freely generated by $\{ \varphi_0, \varphi_1, \dots, \varphi_p\}$, and it is contained in the 
center of $\MM (\Sigma_\Delta, \PP_\Delta)$.}
\qed

\bigskip\noindent
Let $\tilde K$ be the subgroup of $\MM (\Sigma_\Delta, \PP_\Delta)$ generated by $\{ \tau_{c_1}, \dots, 
\tau_{c_q}, \tau_{d_0^{(1)}}, \tau_{d_1^{(1)}}, \dots, \tau_{d_p^{(1)}}, 
\linebreak
\tau_{d_0^{(2)}}, 
\tau_{d_1^{(2)}}, \dots, \tau_{d_p^{(2)}} \}$, and, for $1 \le j\le l$, let $\tilde K_j = \tilde K \cap 
\MM (\Sigma_{\Delta\, j}, \PP_{\Delta\, j})$. Then $\tilde K_j$ is the subgroup of $\MM 
(\Sigma_{\Delta\, j}, \PP_{\Delta\, j})$ generated by the Dehn twists along the boundary components of 
$\Sigma_{\Delta\, j}$, and we have $\tilde K = \tilde K_1 \times \dots \times \tilde K_l$.

\bigskip\noindent
Define the {\it stabilizer} of $\Delta$ to be $\Stab (\Delta)= \{ f \in \MM (\Sigma, \PP); f (\Delta) = 
\Delta \}$. The image of $\theta_\Delta$ is contained in $\Stab( \Delta)$, but, in general, is not 
equal to $\Stab (\Delta)$. However, if $f \in \Stab(\Delta)$ is pure, 
then $f$ fixes each $d_i$ up to isotopy,
thus it belongs to $\Im 
(\theta_\Delta)$. Moreover, if $(f_1, \dots, f_l) \in \theta_\Delta^{-1} (f)$, then each $f_j$ is 
a pure element of $\MM (\Sigma_{\Delta\, j}, \PP_{\Delta\, j})$.

\bigskip\noindent
{\bf Reduction systems.}
Let $f \in \MM (\Sigma, \PP)$ be a pure element. A simplex $\Delta$ of $\CC (\Sigma, \PP)$ is called a 
{\it reduction system} for $f$ if $f(\Delta)= \Delta$ (that is, if $f \in \Stab(\Delta)$). Let $(f_1, 
\dots, f_l) \in \theta_\Delta^{-1} (f)$. If, moreover, for all $1 \le j \le l$, either $f_j$ is a 
pseudo-Anosov element of $\MM (\Sigma_{\Delta\, j}, \PP_{\Delta\, j})$, or $f_j$ is an element of 
$\tilde K_j$, then the reduction system $\Delta$ is called a {\it complete reduction system} for $f$. 
If $f$ is reducible, then the intersection of all the complete reduction systems is a complete 
reduction system itself, called the {\it canonical reduction system} for $f$ and denoted by 
$\Delta(f)$. 
If $f$ is 
a non-reducible pure element (for example, if $f$ is either a pseudo-Anosov element or a Dehn twist 
along a boundary component of $\Sigma$), then the {\it canonical reduction system} for $f$ is defined 
to be $\Delta(f)=\emptyset$. If $f \in \MM(\Sigma,\PP)$ is a pure element and $m \in \Z \setminus 
\{0\}$, then $f^m$ is also pure and $\Delta (f^m)= \Delta(f)$. In particular, if $f$ and $g$ are both 
pure elements and $f^m=g^m$ for some $m \in \Z \setminus \{0\}$, then $\Delta(f) = \Delta(g)$.

\bigskip\noindent
The notion of a canonical reduction system can be extended to all the elements of $\MM (\Sigma, \PP)$ 
using the following result which can be found for instance in \cite{Ivano2} and \cite{BiLuMc1}.

\bigskip\noindent
{\bf Proposition 2.23.} {\it Let $f \in \MM (\Sigma, \PP)$ be a reducible element. Then there exists a 
unique simplex $\Delta= \Delta(f)$ of $\CC (\Sigma, \PP)$ such that:
\begin{itemize}
\item
$f(\Delta) = \Delta$;
\item
there exists $m \ge 1$ such that $f^m$ is pure and $\Delta= \Delta(f^m)$.
\end{itemize}}
\qed

\bigskip\noindent
{\bf Canonical reduction systems.}
The above simplex $\Delta(f)$ is called the {\it canonical reduction system} for $f$, and, 
if $f$ is non-reducible (namely, is either pseudo-Anosov or periodic), then the {\it canonical reduction 
system} for $f$ is defined to be $\Delta(f)= \emptyset$.
 

\subsection{The Kerckhoff realization theorem}

Throughout this subsection the surface $\Sigma$ is again assumed to be closed, $|\PP| \ge 3$ if 
$\Sigma$ is the sphere, and $|\PP| \ge 1$ if $\Sigma$ is the torus.

\bigskip\noindent
{\bf Theorem 2.24 (Kerckhoff \cite{Kerck1}).} {\it Assume $\PP = \emptyset$. Let $\Gamma$ be a finite 
subgroup of $\MM (\Sigma)$. Then there exist a metric $g$ on $\Sigma$ with constant curvature $-1$ and 
a finite subgroup $\tilde \Gamma \subset \Isom (\Sigma, g)$ such that $\tilde \Gamma$ is sent isomorphically 
onto $\Gamma$ under the natural homomorphism $\Isom (\Sigma, g) \to \MM (\Sigma)$.}
\qed

\bigskip\noindent
The above result, known as the Kerckhoff realization theorem, can be extended to the mapping class 
groups of the punctured surfaces as follows.

\bigskip\noindent
{\bf Theorem 2.25.} {\it Let $\Gamma$ be a finite subgroup of $\MM (\Sigma, \PP)$. Then there exists a 
finite subgroup $\tilde \Gamma \subset \Diff (\Sigma, \PP)$ which is sent isomorphically onto $\Gamma$ under 
the natural homomorphism $\Diff (\Sigma, \PP) \to \MM (\Sigma, \PP)$.}

\bigskip\noindent
A proof of Theorem 2.25 which extends Kerckhoff's arguments is sketched in \cite{Ivano3}. 
However, one can also use Theorem 2.24 together with the following construction (which may be useful 
for other purposes) in order to achieve the result.

\bigskip\noindent
Write $\PP= \{ P_1, \dots, P_n\}$. For each $i \in \{ 1, \dots, n\}$, we choose a small disk $\D_i$ 
embedded in $\Sigma$ and containing $P_i$ in its interior. Moreover, we choose the $\D_i$'s so that 
$\D_i \cap \D_j = \emptyset$ for all $1 \le i \neq j \le n$. We denote by 
$\stackrel{{\small o}}{{\mathbb D}}_i$ the interior 
of $\D_i$ and by $c_i: \S^1 \to \D_i$ its boundary component. Set $\Sigma_1= \Sigma \setminus ( 
\cup_{i=1}^n \stackrel{{\small o}}{{\mathbb D}}_i)$. Then $\Sigma_1$ is a compact surface whose boundary 
components are $c_1^{-1}, \dots, c_n^{-1}$. Now, take a copy $\Sigma_1'$ of $\Sigma_1$ and define the 
closed surface $\hat{\Sigma} = (\Sigma_1 \sqcup \Sigma_1') / \sim$, where $\sim$ is the equivalence 
relation which identifies $c_i(z) \in \partial \Sigma_1$ with $c_i(\bar{z})
= c_i^{-1}(z) \in \partial \Sigma_1'$ 
for all $1 \le i\le n$ and all $z \in \S^1$ (see Figure 2.3).

\begin{figure}[htbp]
\centerline{
\setlength{\unitlength}{.4cm}
\begin{picture}(26,9)
\put(1,1){\includegraphics[width=9.6cm]{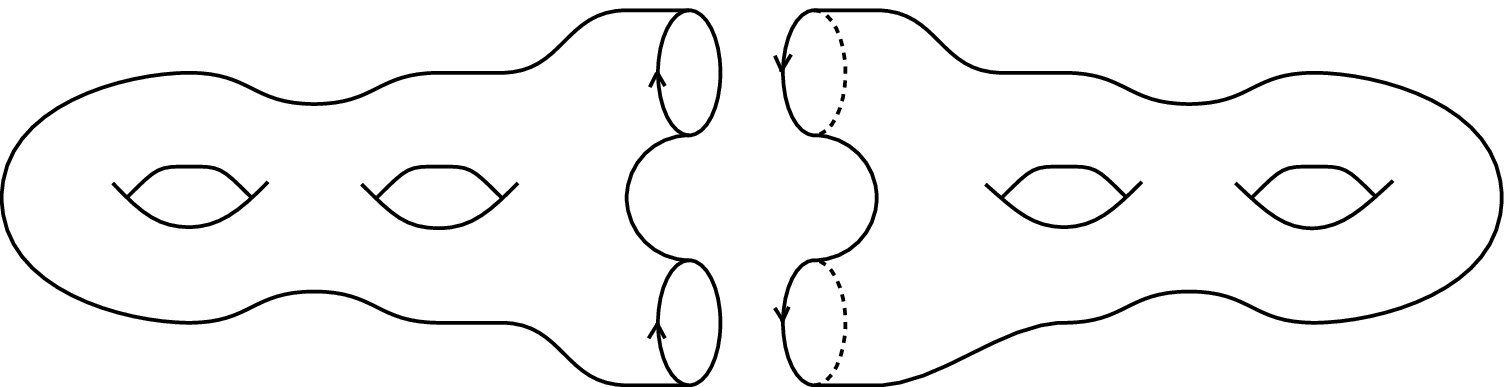}}
\put(10.5,2){\small $c_2$}
\put(10.5,6){\small $c_1$}
\put(6,6.5){$\Sigma_1$}
\put(19,6.5){$\Sigma_1'$}
\end{picture}}
\centerline{{\bf Figure 2.3.} The surface $\hat \Sigma$.}
\end{figure}

\bigskip\noindent
We take a diffeomorphism $\eta: \Sigma_1 \to \Sigma_1$ of order 2, which reverses the orientation, and 
such that $(\eta \circ c_i) (z)= c_i(\bar{z})$ for all $1 \le i \le n$ and all $z \in \S^1$. Let $f 
\in \MM(\Sigma, \PP)$, and let $\sigma \in \Sym_n$ such that $f(P_i)= P_{\sigma(i)}$ for all $1 \le i 
\le n$. We choose a representative $F: \Sigma \to \Sigma$ of $f$ such that $F \circ c_i= c_{\sigma(i)}$ 
and $F(\D_i)= \D_{\sigma(i)}$ for all $1 \le i\le n$. Then $F$ restricts to a diffeomorphism $F: 
\Sigma_1 \to \Sigma_1$ (which does not necessarily belong to $\Diff(\Sigma_1)$ since it does not need 
to be the identity on $\partial \Sigma_1$). Let $\hat{F}: \hat{\Sigma} \to \hat{\Sigma}$ be defined by 
$\hat{F}|_{\Sigma_1}= F$ and $\hat{F}|_{\Sigma_1'}= \eta \circ F \circ \eta$. Then $\hat{F}$ is a 
well-defined homeomorphism, and it is a diffeomorphism on $\hat{\Sigma} \setminus (\cup_{i=1}^n c_i)$. 
Moreover, we can and do choose $F$ such that $\hat{F}$ is a diffeomorphism on the whole surface $\hat{\Sigma}$.

\bigskip\noindent
The proof of the following lemma is left to the reader.

\bigskip\noindent
{\bf Lemma 2.26.} {\it The mapping $\MM (\Sigma, \PP) \to \Diff(\hat{\Sigma})$, $f \mapsto \hat{F}$ 
determines a well-defined injective homomorphism $\iota: \MM (\Sigma, \PP) \to \MM (\hat{\Sigma})$.}
\qed

\bigskip\noindent
{\bf Proof of Theorem 2.25.} Let $\Gamma$ be a finite subgroup of $\MM (\Sigma, \PP)$. Consider the image 
$\iota (\Gamma)$ of $\Gamma$ under the monomorphism $\iota: \MM (\Sigma, \PP) \to \MM (\hat{\Sigma})$. By 
Theorem 2.24, there exist a metric $\hat{g}$ on $\hat{\Sigma}$ of constant curvature $-1$ and a finite 
subgroup $\hat{\Gamma}$ of $\Isom (\hat{\Sigma}, \hat{g})$ such that $\hat{\Gamma}$ is sent isomorphically 
onto $\iota (\Gamma)$ under the natural homomorphism $\Isom (\hat{\Sigma}, \hat{g}) \to \MM 
(\hat{\Sigma})$. We can and do assume that each $c_i$ is the unique geodesic circle in its isotopy 
class. Let $\hat{F} \in \hat{\Gamma}$. Then $\hat{F} (\Sigma_1)= \Sigma_1$, $\hat{F} (\Sigma_1')= 
\Sigma_1'$, and there exists a permutation $\sigma \in \Sym_n$ such that $\hat{F}(c_i)= c_{\sigma(i)}$ 
for all $1 \le i\le n$. We endow each $\D_i$ with the flat metric of radius ${\rm length}(c_i)/2 \pi$ 
and assume that $P_i$ is the center of $\D_i$. Then $\hat{F}|_{\Sigma_1}$ extends in a unique way to an 
isometry $F: \Sigma \to \Sigma$ which satisfies $F(\PP)=\PP$. Now, the mapping $\hat{\Gamma} \to 
\Isom(\Sigma, \PP) \subset \Diff (\Sigma, \PP)$, $\hat{F} \mapsto F$, is a well-defined homomorphism. 
Its image, $\tilde{\Gamma}$, is a finite subgroup of $\Diff(\Sigma, \PP)$ which is sent isomorphically 
onto $\Gamma$ under the natural homomorphism $\Diff(\Sigma, \PP) \to \MM(\Sigma, \PP)$.
\qed

\bigskip\noindent
We finish this subsection with a simple and well-known result whose proof is left to the reader.

\bigskip\noindent
{\bf Lemma 2.27.} {\it Let $\Sigma$ be a closed surface, and let $\Gamma \subset \Diff (\Sigma)$ be a 
finite subgroup of order $m \ge 1$. Then $\Sigma/ \Gamma$ is a surface 
and the natural quotient $\pi: \Sigma \to \Sigma / \Gamma$ is a ramified covering. Let $Q_1, \dots, 
Q_l$ be the ramification points of $\pi$, and, for $1 \le i \le l$, let $o(Q_i)= | \pi^{-1} (Q_i)|$
be the cardinality of $\pi^{-1} (Q_i)$. 
Then the Euler characteristics of $\Sigma$ and $\Sigma/\Gamma$ are related by the formula
\[
\chi (\Sigma) + \sum_{ i=1}^l (m -o(Q_i)) = m \cdot \chi (\Sigma / \Gamma)\,.
\]}
\qed
 
\section{Low genus surfaces}

\subsection{Genus 0 surfaces}

A group $\Gamma$ is {\it biorderable} if there exists a linear ordering $<$ on $\Gamma$ invariant by 
left and right multiplications (namely, if $f <g$, then $h_1fh_2 < h_1gh_2$, for all $f,g,h_1, h_2 \in 
\Gamma$).

\bigskip\noindent
{\bf Example.} Let $\pi: \BB_n \to \Sym_n$ be the natural epimorphism from the braid group $\BB_n$ into 
the symmetric group $\Sym_n$. Then the kernel of $\pi$, denoted by $\PP\BB_n$ and called the {\it pure 
braid group on $n$ strands}, is biorderable (see \cite{KimRol1}).

\bigskip\noindent
The biorderable groups have the following property.

\bigskip\noindent
{\bf Lemma 3.1.} {\it Let $\Gamma$ be a biorderable group, and let $f,g \in \Gamma$. If $f^m=g^m$ for 
some $m \ge 1$, then $f=g$.}

\bigskip\noindent
{\bf Proof.} The fact that $<$ is invariant by left and right multiplications implies that $f^m<g^m$ if 
$f<g$, and $f^m>g^m$ if $f>g$. So, if $f^m=g^m$, then $f=g$.
\qed

\bigskip\noindent
Now, the aim of the present subsection is to prove the following.

\bigskip\noindent
{\bf Proposition 3.2.} {\it Let $\Sigma$ be a surface of genus $0$. Then $\MM(\Sigma)$ is biorderable. 
In particular, if $f,g \in \MM(\Sigma)$ are such that $f^m=g^m$ for some $m \ge 1$, then $f=g$.}

\bigskip\noindent
{\bf Remark.} A group $\Gamma$ is {\it left-orderable} if there exists a linear ordering $<$ on 
$\Gamma$ invariant by left multiplication. By \cite{RouWie1}, the mapping class group $\MM (\Sigma, 
\PP)$ is left-orderable if $\partial \Sigma \neq \emptyset$. If $\partial \Sigma = \emptyset$, then 
$\MM (\Sigma, \PP)$ is not left-orderable since it has torsion, and a left-orderable group cannot
have torsion.

\bigskip\noindent
The following lemma is a preliminary to the proof of Proposition 3.2.

\bigskip\noindent
{\bf Lemma 3.3.} {\it Let $1 \to A \to B \xrightarrow[]{\varphi} C \to 1$ be an exact sequence of 
groups such that $A \simeq \Z^q$ for some $q\ge 0$, $A$ is contained in the center of $B$, and 
$C$ is biorderable. Then $B$ is biorderable, too.}

\bigskip\noindent
{\bf Proof.} Let $<_C$ be a linear ordering on $C$ invariant by left and right multiplications, and let 
$<_A$ be a linear ordering on $A$ invariant by (left and right) multiplications (take the lexicographic 
ordering, for example). Define the linear ordering $<_B$ on $B$ by:
\begin{description}
\item[]
$\beta_1 <_B \beta_2$ if either $\varphi(\beta_1) <_C \varphi(\beta_2)$, or $\varphi(\beta_1)= 
\varphi(\beta_2)$ and $\beta_1 \beta_2^{-1} <_A 1$.
\end{description}
It is easily checked that $<_B$ is invariant by left and right multiplications.
\qed

\bigskip\noindent
{\bf Proof of Proposition 3.2.} Let $\Sigma$ be a surface of genus $0$ and $q$ boundary components,  
$c_1, \dots, c_q$. If $q \le 1$, then $\MM(\Sigma) = \{1\}$, thus we may assume $q 
\ge 2$. Let $\D= \{ z \in \C; |z| \le 1\}$ denote the standard disk. For all $1 \le i \le q-1$, we take 
a copy $\D_i$ of $\D$ and we denote by $d_i: \S^1 \to \D_i$ its boundary component. Define the surface 
$\Sigma_0= (\Sigma \sqcup (\sqcup_{1 \le i\le q-1} \D_i))/ \sim$, where $\sim$ is the equivalence 
relation which identifies $c_i(z)$ with $d_i(\bar z) = d_i^{-1} (z)$ for all $1 \le i \le q-1$ and all $z \in \S^1$. 
Note that $\Sigma_0$ is (topologically) a disk and its unique boundary component is $c_q$. We choose a 
point $Q_i$ in the interior of $\D_i$ for all $1 \le i\le q-1$, and we set $\QQ= \{Q_1, \dots, Q_{q-
1}\}$. The group $\MM (\Sigma_0, \QQ)$ is (isomorphic to) the braid group $\BB_{q-1}$ on $q-1$ strands. 
Moreover, the embedding $\Sigma \to \Sigma_0$ determines a homomorphism $\theta: \MM(\Sigma) \to 
\MM(\Sigma_0, \QQ)= \BB_{q-1}$ whose image is the pure braid group $\PP\BB_{q-1}$, a biorderable group 
(see \cite{KimRol1}). Let 
$K$ be the subgroup of $\MM(\Sigma)$ generated by $\{\tau_{c_1}, \dots, \tau_{c_{q-1}}\}$. Then, 
by Proposition 2.1, the kernel of 
$\theta$ is $K$, it is a free abelian group of rank $q-1$, and it is contained in the center of 
$\MM(\Sigma)$ (see \cite{ParRol1}). We conclude by Lemma 3.3 that $\MM(\Sigma)$ is biorderable.
\qed


\subsection{Genus 1 surfaces}

Let $\T^2$ denote the 2-dimensional torus. Let $a,b$ be the closed curves pictured in Figure 3.1, and 
let $\alpha= \tau_a \tau_b$, $\beta= \tau_a \tau_b \tau_a$, and $\delta= (\tau_a \tau_b)^3$. Then 
$\MM(\T^2)$ is well-known to be isomorphic to $SL_2(\Z)$ and has a presentation with generators 
$\alpha, \beta, \delta$ and relations
\[
\alpha^3= \beta^2 = \delta\,, \quad \delta^2=1\,.
\]
The center $Z(\MM (\T^2))$ of $\MM(\T^2)$ is the cyclic group $\{1, \delta \}$ of order 2, and $\MM 
(\T^2) / Z(\MM(\T^2)) = \CC_3 \ast \CC_2$ is the free product of $\CC_3= \{ 1, \bar \alpha, \bar 
\alpha^2 \}$, a cyclic group of order 3, and $\CC_2= \{ 1, \bar \beta \}$, a cyclic group of order 2.

\begin{figure}[htbp]
\centerline{
\setlength{\unitlength}{.4cm}
\begin{picture}(8,6)
\put(1,1){\includegraphics[width=2.4cm]{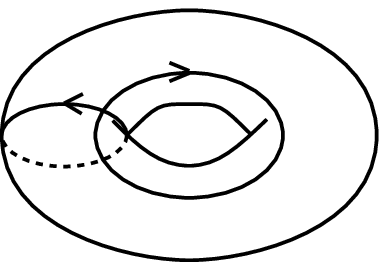}}
\put(2,3.8){\small $a$}
\put(4,4.1){\small $b$}
\end{picture}}
\centerline{{\bf Figure 3.1.} The torus $\T^2$.}
\end{figure}

\bigskip\noindent
Using the fact that any finite order element of $\CC_3 \ast \CC_2$ is conjugate to an element of either 
$\CC_3$ or $\CC_2$ (see \cite{Serre1}), one can easily prove the following.

\bigskip\noindent
{\bf Lemma 3.4.} {\it There is a unique element of order 2 in $\MM(\T^2)$: $\delta$. There are two 
conjugacy classes of elements of order 4 represented by $\beta$ and $\beta^{-1}$. There are two 
conjugacy classes of elements of order 3 represented by $\alpha^2$ and $\alpha^{-2}$. There are two 
conjugacy classes of elements of order 6 represented by $\alpha$ and $\alpha^{-1}$. These are all the 
conjugacy classes of elements of finite order in $\MM(\T^2)$.}
\qed

\bigskip\noindent
A $m$-root of a given element of infinite order in $\MM(\T^2)$ is not necessarily unique up to 
conjugation. However, such a $m$-root is unique up to the center as we turn now to show.

\bigskip\noindent
{\bf Proposition 3.5.} {\it Let $f,g \in \MM (\T^2)$ be two elements of infinite order. 
If $f^m=g^m$ for some $m \ge 1$, then either $f=g$ or 
$f=\delta g$.}

\bigskip\noindent
{\bf Proof.}
First, observe that and element $f \in SL_2 (\Z) \simeq \MM( \T^2)$ has infinite order if and only if it has real 
eigenvalues and $f \not\in \{ \Id, -\Id\}$. (The only non-trivial part of the above statement is that 
an element $f \in SL_2(\Z)$ with non-real eigenvalues has finite order. In order to prove that, notice 
that the norm $\| f^m\|$ of any power of such a $f$ is bounded by a constant $K$, and there are finitely 
many elements $g \in SL_2 (\Z)$ such that $\| g \| \le K$.)

\bigskip\noindent
Now, take two elements $f,g \in SL_2(\Z)$ of infinite order such that $f^m = g^m$ for some $m \ge 1$. 
In particular, $f$ and $g$ have real eigenvalues. The eigenspaces of $f$ (resp. $g$) coincide with the 
eigenspaces of $f^m = g^m$, hence $f$ and $g$ have the same eigenspaces.

\bigskip\noindent
{\bf Case 1 :} $f$ has two distinct eigenvalues $\lambda$ and $1/ \lambda$. Then $g$ has also two 
distinct eigenvalues $\mu$ and $1/\mu$. The equality $f^m=g^m$ implies that $\lambda^m = \mu^m$, thus 
$\lambda = \pm \mu$ and $f = \pm g$.

\bigskip\noindent
{\bf Case 2 :}
$f$ has a unique eigenvalue (which is equal to $1$ or $-1$). Then there is a basis of $\R^2$ (indeed of 
$\Z^2$) relatively to which $f$ and $g$ have the forms
$\left(\begin{array}{cc} \varepsilon_f&\nu_f\\ 0&\varepsilon_f \end{array}\right)$ and
$\left(\begin{array}{cc} \varepsilon_g&\nu_g\\ 0&\varepsilon_g \end{array}\right)$,
respectively, where $\varepsilon_f, \varepsilon_g \in \{ \pm 1\}$. Now, it is easily deduced from the 
equality $f^m=g^m$ that $f= \pm g$.
\qed

\bigskip\noindent
We turn now 
to study the mapping class groups of the genus 1 surfaces with non-empty boundary.
Our main result is the following.

\bigskip\noindent
{\bf Theorem 3.6.} {\it Let $\Sigma$ be a surface of genus 1 with non-empty boundary. Let $f,g \in \MM 
(\Sigma)$. If $f^m=g^m$ for some $m \ge 1$, then $f$ and $g$ are conjugate.}

\bigskip\noindent
The proof of Theorem 3.6 is long and tedious, so we put it separately in the last section and go
directly to the next section and the study of the mapping class groups of the high genus surfaces.


\section{Surfaces with non-empty boundary}

We start with a simple lemma which will be used later.

\bigskip\noindent
{\bf Lemma 4.1.} {\it Let $\Sigma$ be a closed surface (i.e. $\partial \Sigma = \emptyset$), and let 
$\PP$ be a non-empty finite set of punctures in $\Sigma$. Let $P_0 \in \PP$, and let $\Gamma$ be a 
finite subgroup of $\MM (\Sigma, \PP)$ such that $h(P_0)=P_0$ for all $h \in \Gamma$. Then $\Gamma$ is 
cyclic.}

\bigskip\noindent
{\bf Proof.} By Theorem 2.25, there exists a finite subgroup $\tilde \Gamma \subset \Diff (\Sigma, 
\PP)$ which is sent isomorphically onto $\Gamma$ under the natural homomorphism $\Diff (\Sigma, \PP) 
\to \MM (\Sigma, \PP)$. Moreover, we have $H(P_0)=P_0$ for all $H \in \tilde \Gamma$. Since $\tilde 
\Gamma$ is finite, there exists a Riemannian metric $g$ on $\Sigma$ which is invariant under the action 
of $\tilde \Gamma$, that is, $\tilde \Gamma \subset \Isom_+(\Sigma, g)$. Now, the fact that 
$H(P_0)=P_0$ for all $H \in \tilde \Gamma$ implies that $\tilde \Gamma$ is a finite subgroup of $O_+( 
T_{P_0} \Sigma) \simeq \S^1$, hence $\tilde \Gamma$ is cyclic.
\qed

\bigskip\noindent
Now, the bad news of the section are the following.

\bigskip\noindent
{\bf Proposition 4.2.} {\it Let $\Sigma$ be a surface of genus $\rho \ge 2$ and 1 boundary component, 
$c$. Then there exist two periodic elements $f,g \in \MM (\Sigma)$ such that $f^2=g^2= \tau_c$, $f$ and 
$g$ are not conjugate, and none of the conjugates of $f$ commutes with $g$.}

\bigskip\noindent
{\bf Proof.} Let $\A = [0,1] \times \S^1$ denote the standard annulus. Let $a_k, a_k': [0,1] \to 
\partial \A$, $1 \le k\le 2 \rho$, be the arcs defined by
\[
a_k(t)= (1, \exp({i \pi \over 2 \rho} (k-1+t)))\,, \quad a_k'(t)= (1, -\exp( {i\pi \over 2 \rho} (k-
t)))\,.
\]
(See Figure 4.1.) Then $\Sigma$ can be presented as the quotient $\Sigma= \A /\sim_f$, where $\sim_f$ 
is the equivalence relation which identifies $a_k(t)$ with $a_k'(t)$ for all $1 \le k\le 2 \rho$ and 
all $t \in [0,1]$. Let $\pi_f: \A \to \Sigma$ be the natural quotient, and, for $1 \le k\le 2 \rho$, 
let $\bar a_k= \pi_f \circ a_k = \pi_f \circ a_k'$. Then the $\bar a_k$'s are loops in $\Sigma$ and $\{ [ 
\bar a_k]; 1 \le k\le 2 \rho \}$ is a basis for $H_1 (\Sigma, \R) \simeq \R^{2 \rho}$, where $[ \bar 
a_k]$ denotes the homology class of $\bar a_k$. We choose a smooth map $\varphi: [0,1] \to [0,1]$ such 
that $\varphi(0)=0$, $\varphi(1)=1$, $\varphi'(0)= \varphi'(1)=0$, and $\varphi'(t) \ge 0$ for all $t 
\in [0,1]$, and we define $F_\A: \A \to \A$ by
\[
F_\A (t,z)= (t, \exp( i \pi \varphi(t)) \cdot z)\,.
\]
Then $F_\A$ induces a diffeomorphism $F \in \Diff (\Sigma)$ via $\pi_f$, and we define $f \in \MM 
(\Sigma)$ as the isotopy class of $F$. Note that $F \circ \bar a_k = \bar a_k^{-1}$ for all $1 \le k 
\le 2 \rho$, thus the spectrum of the action of $f$ on $H_1(\Sigma, \R)$ is $(-1, \dots, -1)$. 
Moreover, we have $f^2 = \tau_c$.

\begin{figure}[htbp]
\centerline{
\setlength{\unitlength}{.4cm}
\begin{picture}(12,12)
\put(2,2){\includegraphics[width=3.2cm]{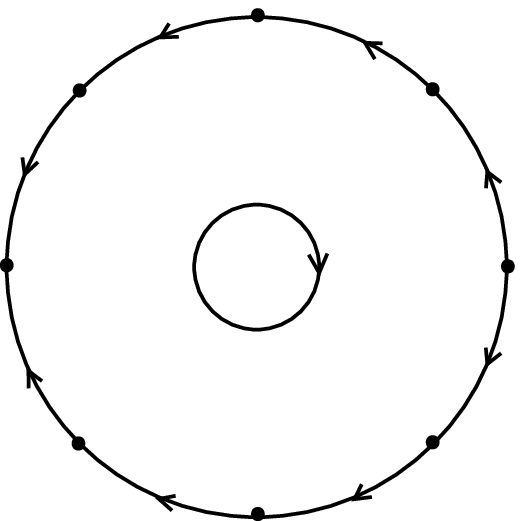}}
\put(10,7.5){\small $a_1$}
\put(7.5,10){\small $a_2$}
\put(3.5,10.1){\small $a_3$}
\put(1,7.5){\small $a_4$}
\put(1.2,3.5){\small $a_1'$}
\put(4,1.2){\small $a_2'$}
\put(7.5,1.3){\small $a_3'$}
\put(10,4){\small $a_4'$}
\put(6.3,5.7){\small $c$}
\end{picture}}
\centerline{{\bf Figure 4.1.} A surface of genus 2 with 1 boundary component.}
\end{figure}

\bigskip\noindent
Let $b_j, c_j, b_j' c_j': [0,1] \to \partial \A$, $j=1,2$, be the arcs defined by
\begin{gather*}
b_j(t)= (1, \exp( {i \pi \over 2 \rho} (2j \rho -4 +t)))\,, \quad c_j(t)= (1, \exp( {i \pi \over 2 
\rho} (2j \rho -3 +t)))\,,\\
b_j'(t)= (1, \exp( {i \pi \over 2 \rho} (2 j \rho -1 -t)))\,, \quad c_j'(t) = (1, \exp( {i \pi \over 2 
\rho} (2 j \rho -t)))\,.
\end{gather*}
(See Figure 4.2.) The surface $\Sigma$ can also be presented as the quotient $\Sigma= \A / \sim_g$, 
where $\sim_g$ is the equivalence relation which identifies $a_k(t)$ with $a_k'(t)$ for all $1 \le k\le 
2 \rho -4$ and all $t \in [0,1]$, and which identifies $b_j(t)$ with $b_j'(t)$, and $c_j(t)$ with 
$c_j'(t)$, for all $j \in \{ 1,2\}$ and all $t \in [0,1]$. Let $\pi_g: \A \to \Sigma$ be the natural 
quotient, let $\bar a_k = \pi_g \circ a_k= \pi_g \circ a_k'$ for $1 \le k\le 2 \rho -4$, and let $\bar 
b_j = \pi_g \circ b_j = \pi_g \circ b_j'$ and $\bar c_j = \pi_g \circ c_j = \pi_g \circ c_j'$ for 
$j=1,2$. Then all the $\bar a_k$'s , the $\bar b_j$'s, and the $\bar c_j$'s are loops in $\Sigma$, and 
the set $\{ [\bar a_k]; 1 \le k \le 2 \rho -4 \} \cup \{ [\bar b_1], [\bar b_2], [\bar c_1], [\bar 
c_2]\}$ is a basis for $H_1(\Sigma, \R)$. Here again, $F_\A: \A \to \A$ induces a diffeomorphism $G \in 
\Diff (\Sigma)$ via $\pi_g$, and we define $g \in \MM (\Sigma)$ as the isotopy class of $G$. Note that 
$G \circ \bar a_k = \bar a_k^{-1}$ for all $1 \le k \le 2 \rho -4$, $G \circ \bar b_1 = \bar b_2$, $G 
\circ \bar b_2 = \bar b_1$, $G \circ \bar c_1= \bar c_2$, $G \circ \bar c_2 = \bar c_1$, thus the 
spectrum of the action of $g$ on $H_1(\Sigma, \R)$ is $(-1, \dots, -1,-1,-1,1,1)$. Moreover, we have 
$g^2= \tau_c$.

\begin{figure}[htbp]
\centerline{
\setlength{\unitlength}{.4cm}
\begin{picture}(12,12)
\put(2,2){\includegraphics[width=3.2cm]{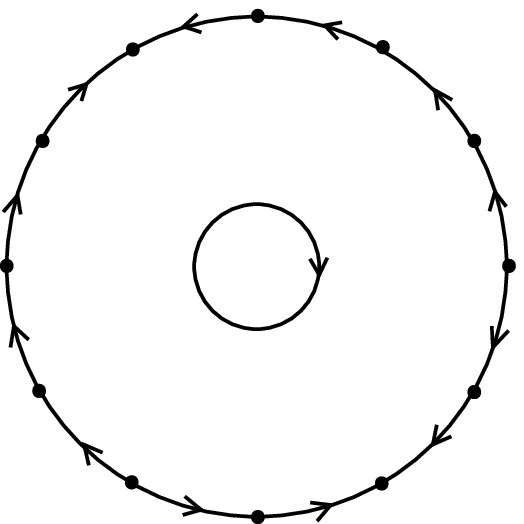}}
\put(10,7){\small $a_1$}
\put(9,9){\small $a_2$}
\put(7,10.3){\small $b_1$}
\put(4.2,10.4){\small $c_1$}
\put(2,9){\small $b_1'$}
\put(1,7){\small $c_1'$}
\put(1,4.5){\small $a_1'$}
\put(2,2.5){\small $a_2'$}
\put(4.5,1.1){\small $b_2$}
\put(6.7,1.3){\small $c_2$}
\put(9,2.2){\small $b_2'$}
\put(10,4.5){\small $c_2'$}
\put(6.3,5.9){\small $c$}
\end{picture}}
\centerline{{\bf Figure 4.2.} A surface of genus 3 with 1 boundary component.}
\end{figure}

\bigskip\noindent
We already know that $f^2=g^2 = \tau_c$. On the other hand, the spectrum of the action of $f$ on 
$H_1(\Sigma, \R)$ is different from the spectrum of the action of $g$, thus $f$ 
and $g$ are not conjugate. Suppose that there exists $f' \in \MM (\Sigma)$ conjugate to $f$ and commuting 
with $g$. Let $\theta: \MM (\Sigma) \to \MM (\Sigma_0, \QQ)$ be the corking of $\Sigma$. Then $\theta 
(f') (Q_1)= \theta (g) (Q_1)=Q_1$, and $\{ \theta(f'), \theta (g) \}$ generates a subgroup of $\MM 
(\Sigma_0, \QQ)$ isomorphic to $\Z /2 \Z \times \Z /2 \Z$. 
This is a finite subgroup of $\MM (\Sigma_0, \QQ)$ whose elements fix $Q_1$ and which is not 
cyclic, contradicting Lemma 4.1.
\qed

\bigskip\noindent
{\bf Corollary 4.3.} {\it Let $\Sigma$ be a surface of genus $\rho \ge 2$ and $q$ boundary components, 
$c_1, \dots, c_q$, where $q \ge 1$, and let $\PP = \{ P_1, \dots, P_n \}$ be a finite set of punctures 
in the interior of $\Sigma$. Then there exist $f,g \in \MM (\Sigma, \PP)à$ such that $f^2=g^2$, $f$ and 
$g$ are not conjugate, and none of the conjugates of $f$ commutes with $g$.}

\bigskip\noindent
{\bf Proof.} Obviously, we can assume that either $q \ge 2$ or $\PP \neq \emptyset$. 
(The case $q=1$ and $\PP= \emptyset$ is covered by the conclusion of Proposition 4.2.)
Let $\Sigma_1$ be 
a surface of genus $\rho$ and $1$ boundary component, $c^{(1)}$, and let $\Sigma_2$ be a surface of 
genus 0 and $q+1$ boundary components, $c_1, \dots, c_q, c^{(2)}$. The surface $\Sigma$ can be 
presented as the quotient $\Sigma= (\Sigma_1 \sqcup \Sigma_2) / \sim$, where $\sim$ is the equivalence 
relation which identifies $c^{(1)}(z)$ with $c^{(2)}(\bar z)= (c^{(2)})^{-1} (z)$ for all $z \in \S^1$ 
(see Figure 4.3). We 
denote by $\pi: \Sigma_1 \sqcup \Sigma_2 \to \Sigma$ the natural quotient, we assume that $\PP$ is 
contained in the interior of $\pi(\Sigma_2)$, and we set $c= \pi \circ c^{(1)} = \pi \circ (c^{(2)})^{-
1}$ and $\gamma= \langle c \rangle$. Then $\pi$ determines a homomorphism $\theta: \MM (\Sigma_1) 
\times \MM (\Sigma_2, \pi^{-1} (\PP)) \to \MM (\Sigma,\PP)$ whose image is $\Stab (\gamma)= \{h \in \MM 
(\Sigma,\PP); h( \gamma)= \gamma\}$. Furthermore, the restriction of $\theta$ to $\MM (\Sigma_1)$ is 
injective (see \cite{ParRol1}).

\begin{figure}[htbp]
\centerline{
\setlength{\unitlength}{.4cm}
\begin{picture}(23,12)
\put(1,1){\includegraphics[width=8.2cm]{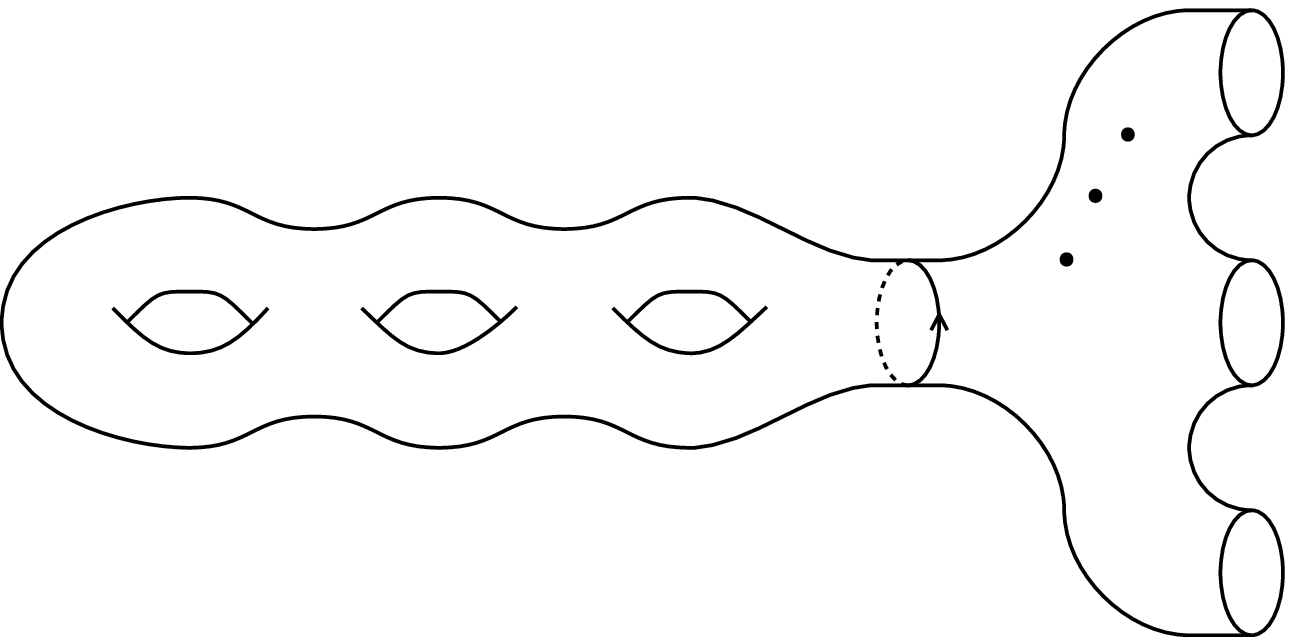}}
\put(16.2,5.9){\small $c$}
\put(17.5,6){\small $P_1$}
\put(18.8,7.5){\small $P_2$}
\put(18.5,9.5){\small $P_3$}
\put(7.5,2.5){$\Sigma_1$}
\put(16.9,2.5){$\Sigma_2$}
\end{picture}}
\centerline{{\bf Figure 4.3.} Cutting $\Sigma$ into two pieces.}
\end{figure}

\bigskip\noindent
By Proposition 4.2, there exist two periodic elements $f_1, g_1 \in \MM (\Sigma_1)$ such that $f_1^2= 
g_1^2 = \tau_{c^{(1)}}$, $f_1$ and $g_1$ are not conjugate in $\MM (\Sigma_1)$, and none of the 
conjugates of $f_1$ commutes with $g_1$. We set $f = \theta (f_1)$ and $g= \theta (g_1)$. Then $f^2= 
g^2 = \tau_c$ and $\Delta(f)= \Delta(g) = \{ \gamma\}$.

\bigskip\noindent
Suppose that there exists $h \in \MM (\Sigma,\PP)$ such that $h f h^{-1} = g$. We have 
\[
\{ \gamma\} = \Delta(g)= \Delta (hfh^{-1}) = h( \Delta (f)) = h (\{ \gamma\})\,,
\]
thus $h \in \Stab (\gamma)$. Let $(h_1, h_2) \in \MM (\Sigma_1) \times \MM (\Sigma_2, \pi^{-1} (\PP))$ 
such that $h = \theta (h_1,h_2)$. Then 
\[
\theta (h_1 f_1 h_1^{-1}) = \theta( (h_1,h_2) \cdot (f_1, \Id) \cdot (h_1^{-1}, h_2^{-1})) = hfh^{-1}
=g = \theta ((g_1, \Id)) = \theta(g_1)\,,
\]
thus $h_1 f_1 h_1^{-1} = g_1$ (since $\theta|_{\MM (\Sigma_1)}$ is injective): a contradiction.

\bigskip\noindent
Suppose that there exists $h \in \MM (\Sigma, \PP)$ such that $(hfh^{-1}) g = g (hfh^{-1})$. Set 
$\gamma' = h(\gamma)$ and choose an essential curve $c': \S^1 \to \Sigma \setminus \PP$ which 
represents $\gamma'$. We have $g^2= \tau_c$, $h f^2 h^{-1} = h \tau_c h^{-1} = \tau_{c'}$, and $g$ 
commutes with $h f h^{-1}$, thus $\tau_c$ and $\tau_{c'}$ commute, therefore $c'$ can be chosen so that 
$c' \cap c = \emptyset$ (see \cite{ParRol1}, for example). Both, $c$ and $c'$, separate the surface 
$\Sigma$ into two subsurfaces, one of genus $\rho$ with one boundary component, and one of genus $0$ 
with $q+1$ boundary components and containing the whole set $\PP$. Obviously, since $c' \cap c = 
\emptyset$, this is possible only if $c'$ and $c$ are isotopic, that is, $\gamma' = \gamma$. So, $h \in 
\Stab (\gamma)$. 
Let $(h_1,h_2) \in \MM (\Sigma_1) \times \MM (\Sigma_2, \pi^{-1} (\PP))$ such that $h= (h_1,h_2)$.
Then $(h_1 f_1 h_1^{-1}) g_1 = g_1 (h_1 f_1 h_1^{-1})$: a contradiction.
\qed

\bigskip\noindent
{\bf Remark.} If $q$ is large enough, then the elements $f$ and $g$ in the above corollary cannot be 
chosen to be periodic because of the following.

\bigskip\noindent
{\bf Lemma 4.4.} {\it Let $\Sigma$ be a surface of genus $\rho \ge 2$ and $q$ boundary components, 
$c_1, \dots, c_q$, where $q \ge 1$, and let $\PP = \{ P_1, \dots, P_n \}$ be a finite set of punctures 
in the interior of $\Sigma$. Let $K$ be the subgroup of $\MM (\Sigma, \PP)$ generated by $\{ 
\tau_{c_1}, \dots, \tau_{c_q}\}$. If $q > 2 \rho +2$, then all the periodic elements of $\MM (\Sigma, 
\PP)$ belong to $K$. In particular, if $q > 2 \rho +2$, and if $f,g \in \MM (\Sigma, \PP)$ are two 
periodic elements such that $f^m=g^m$ for some $m \ge 1$, then $f=g$.}

\bigskip\noindent
{\bf Proof.} Let $\theta: \MM (\Sigma, \PP) \to \MM (\Sigma_0, \PP \sqcup \QQ)$ be the corking of 
$(\Sigma, \PP)$. Let $f \in \MM (\Sigma, \PP)$ be a periodic element such that $\theta(f) \not\in K$. 
This means that $\theta(f)$ is of finite order $m \ge 2$. Let $F \in \Diff (\Sigma_0, \PP \sqcup \QQ)$ 
be a diffeomorphism of order $m$ which represents $\theta(f)$. Since $F(Q_i)=Q_i$ for all $Q_i \in 
\QQ$, by Lemma 2.27, we have
\[
\begin{array}{rll}
&\chi (\Sigma_0) + q(m-1) \le m \cdot \chi (\Sigma_0 /F)\\
\Rightarrow \quad & 2-2 \rho + q(m-1) \le 2m \quad &(\text{since } \chi (\Sigma_0/F) \le 2)\\
\Rightarrow \quad & 2 -2 \rho -q + (q-2)m \le 0\\
\Rightarrow \quad & 2 -2 \rho -q +2 (q-2) \le 0 \quad &(\text{since } m \ge 2)\\
\Rightarrow \quad & q \le 2 \rho +2\,.
\end{array}
\]
The uniqueness of a $m$-root of a given periodic element in $\MM (\Sigma, \PP)$, if $q > 2 \rho +2$, 
is simply due to the fact that $K \simeq \Z^q$ and such a root must lie in $K$.
\qed

\bigskip\noindent
We turn now to the good news of the section.

\bigskip\noindent
{\bf Theorem 4.5.} {\it Let $\Sigma$ be a surface of genus $\rho \ge 0$ and $q$ boundary components, 
$c_1, \dots, c_q$, where $q \ge 1$, and let $\PP = \{ P_1, \dots, P_n \}$ be a finite set of punctures 
in the interior of $\Sigma$. Let $f,g \in \MM (\Sigma, \PP)$ be two pseudo-Anosov elements. If $f^m = 
g^m$ for some $m \ge 1$, then $f=g$.}

\bigskip\noindent
{\bf Proof.}
Let $\theta: \MM (\Sigma, \PP) \to \MM (\Sigma_0, \PP \sqcup \QQ)$ be the corking of $(\Sigma, \PP)$. 
By Corollary 2.14, there are pseudo-Anosov representatives $F,G \in \Diff (\Sigma_0, \PP \sqcup \QQ)$ 
of $\theta(f), \theta(g)$, respectively, such that $F^m = G^m$. Let $\FF^s$ (resp. $\FF^u$) be the 
stable (resp. unstable) foliation of $F$, and let $\lambda$ be its dilatation coefficient. Then $\FF^s$ 
(resp. $\FF^u$) is also the stable (resp. unstable) foliation of $G$ and $\lambda$ is its dilatation 
coefficient. This implies that $F \circ G \circ F^{-1} \circ G^{-1} \in \Sym (\FF^s, \FF^u)$.

\bigskip\noindent
Fix some $Q_1 \in \QQ$. Let $\SS \SS_{Q_1}$ be the set of separatrices of $\FF^s$ at $Q_1$. We 
number $\LL_1, \LL_2, \dots, \LL_k$ the elements of $\SS\SS_{Q_1}$ so that the corresponding prongs 
are anticlockwise numbered. Let 
\linebreak
$\Norm_{Q_1} (\FF^s, \FF^u)$ be the subgroup of $\Norm (\FF^s, \FF^u)$ of 
elements $H \in \Norm (\FF^s, \FF^u)$ that fix $Q_1$. All the elements of $\Norm_{Q_1} 
(\FF^s, \FF^u)$ permute the set $\SS\SS_{Q_1}$, 
and this determines a well-defined homomorphism $\psi: \Norm_{Q_1} (\FF^s, \FF^u) \to \Sym 
(\SS\SS_{Q_1})$. Moreover, if $\omega \in \Sym (\SS\SS_{Q_1})$ denotes the cyclic permutation of $\LL_1, 
\LL_2, \dots, \LL_k$, then $\psi(H)$ is a power of $\omega$ for all $H \in \Norm_{Q_1} (\FF^s, \FF^u)$. In 
particular, if $H_1, H_2 \in \Norm_{Q_1} (\FF^s, \FF^u)$, then $\psi( H_1 \circ H_2 \circ H_1^{-1} \circ 
H_2^{-1}) =1$.

\bigskip\noindent
Note that $F,G \in \Norm_{Q_1} (\FF^s, \FF^u)$, thus, by the above observation, we have $\psi( F \circ G 
\circ F^{-1} \circ G^{-1}) =1$. Furthermore, we already know that $F \circ G \circ F^{-1} \circ G^{-1} 
\in \Sym( \FF^s, \FF^u)$, thus, by Lemma 2.9, $F \circ G \circ F^{-1} \circ G^{-1} = \Id$. So, 
$\theta( fgf^{-1} g^{-1})=1$.

\bigskip\noindent
Let $K$ be the subgroup of $\MM (\Sigma, \PP)$ generated by $\tau_{c_1}, \dots, \tau_{c_q}$. By 
Proposition 2.1, there exists $u \in K$ such that $fgf^{-1} = ug$. Since $u$ is central, it follows 
that
\[
g= g^m g g^{-m} = f^m g f^{-m} = u^m g\,,
\]
thus $u^m=1$. Since $K$ is torsion free, we deduce that $u=1$, that is, $f$ and $g$ commute.

\bigskip\noindent
Set $h= gf^{-1}$. Then
\[
g^m= (hf)^m = h^m f^m = h^m g^m\,,
\]
thus $h^m=1$, therefore $h=1$ (since $\MM (\Sigma, \PP)$ is torsion free), that is, $f=g$.
\qed


\section{Closed surfaces}

The good news of this section are the following.

\bigskip\noindent
{\bf Proposition 5.1.} {\it Let $\Sigma$ be a closed surface of genus $\rho \ge 2$, and let $\PP$ be a 
finite set of punctures in $\Sigma$. Let $f \in \MM (\Sigma, \PP)$ be a pseudo-Anosov element. Then the 
set of $g \in \MM (\Sigma, \PP)$ satisfying $g^m=f^m$ for some $m \ge 1$ is finite.}

\bigskip\noindent
{\bf Proof.} We take a pseudo-Anosov representative $F \in \Diff (\Sigma, \PP)$ of $f$, we denote by 
$\FF^s$ (resp. $\FF^u$) the stable (resp. unstable) foliation of $F$, we denote by $\lambda$ its 
dilatation coefficient, and we set $\tilde \Gamma_0 = \Sym (\FF^s, \FF^u)$. The subgroup $\tilde 
\Gamma_0$ is finite by Corollary 2.10 and, by Corollary 2.12, it is sent isomorphically onto a subgroup 
$\Gamma_0$ of $\MM(\Sigma, \PP)$ under the natural homomorphism $\Diff (\Sigma, \PP) \to \MM (\Sigma, 
\PP)$.

\bigskip\noindent
Set $\UU= \{ g \in \MM (\Sigma, \PP); g^m=f^m \text{ for some } m \ge 1 \}$. Let $g \in \UU$, and let 
$m \ge 1$ such that $g^m= f^m$. Clearly, this last equality implies that $g$ is a pseudo-Anosov element 
of $\MM (\Sigma, \PP)$. By Corollary 2.14, there exists a pseudo-Anosov representative 
$G \in \Diff (\Sigma, \PP)$ of $g$ such that $G^m=F^m$. Clearly, $\FF^s$ (resp. $\FF^u$) is the stable 
(resp. unstable) foliation of $G$, and $\lambda$ is its dilatation coefficient, hence $F^{-1} G \in 
\tilde \Gamma_0 = \Sym (\FF^s, \FF^u)$, therefore $f^{-1} g \in \Gamma_0$. This shows that $\UU 
\subset f \cdot \Gamma_0= \{ fh; h \in \Gamma_0\}$, hence $\UU$ is finite.
\qed

\bigskip\noindent
However, the situation can be pretty bad as we turn now to show.

\bigskip\noindent
{\bf Theorem 5.2.} {\it 
\begin{enumerate}
\item
Let $\Sigma$ be a closed surface of genus $\rho \ge 2$. Then there exist two pseudo-Anosov elements 
$f,g \in \MM (\Sigma)$ such that $f^{2(\rho +1)} = g^{2(\rho+1)}$, $f$ is not conjugate to $g$, and 
none of the conjugates of $f$ commutes with $g$.
\item
Let $\Sigma$ be a closed surface of genus $\rho \ge 4$, with $\rho \equiv 0\ (\mod\, 4)$. Then there 
exist two pseudo-Anosov elements $f,g \in \MM (\Sigma)$ such that $f^2=g^2$, $f$ is not conjugate to 
$g$, and none of the conjugates of $f$ commutes with $g$.
\end{enumerate}}

\bigskip\noindent
{\bf Proof of Theorem 5.2.1.} We consider a closed surface $\Sigma$ of genus $\rho \ge 2$. This admits 
a cellular decomposition (see Figure 5.1) with $2 (\rho+1)$ vertices, $P_1, P_2, \dots, P_{2(\rho 
+1)}$, with $4 (\rho+1)$ arrows, $a_1, a_2, \dots, a_{\rho+1}$, $b_1, b_2, \dots, b_{\rho+1}$, $c_1, 
c_2, \dots, c_{\rho+1}$, $d_1, d_2, \dots, d_{\rho+1}$, and with 4 faces, $A_1, A_2, A_3, A_4$, where 
the sources and the targets of the arrows are
\[
\begin{array}{lll}
\source (a_i)= \source (d_i)= P_{2i-1}\,, &\quad& \target (a_i)= \target (d_i)= P_{2i}\,,\\
\source (b_i)= \source (c_i) = P_{2i}\,, &\quad& \target (b_i) = \target (c_i)= P_{2i+1}\,,
\end{array}
\]
for all $1 \le i \le \rho+1$, and where the boundaries of the faces are
\[
\begin{array}{lll}
\partial A_1= a_1 b_1 a_2 b_2 \dots a_{\rho+1} b_{\rho+1}\,, &\quad& \partial A_2= c_{\rho+1}^{-1} 
a_{\rho+1}^{-1} \dots c_2^{-1} a_2^{-1} c_1^{-1} a_1^{-1}\,,\\
\partial A_3= b_{\rho+1}^{-1} d_{\rho+1}^{-1} \dots b_2^{-1} d_2^{-1} b_1^{-1} d_1^{-1}\,, &\quad& 
\partial A_4= d_1 c_1 d_2 c_2 \dots d_{\rho+1} c_{\rho+1}\,.
\end{array}
\]
Moreover, we endow $\Sigma$ with the hyperbolic metric so that each $A_j$ is isometric to a 
right-angled regular $2(\rho+1)$-gone in the hyperbolic plane.

\begin{figure}[htbp]
\centerline{
\setlength{\unitlength}{.4cm}
\begin{picture}(24,24)
\put(2,2){\includegraphics[width=8cm]{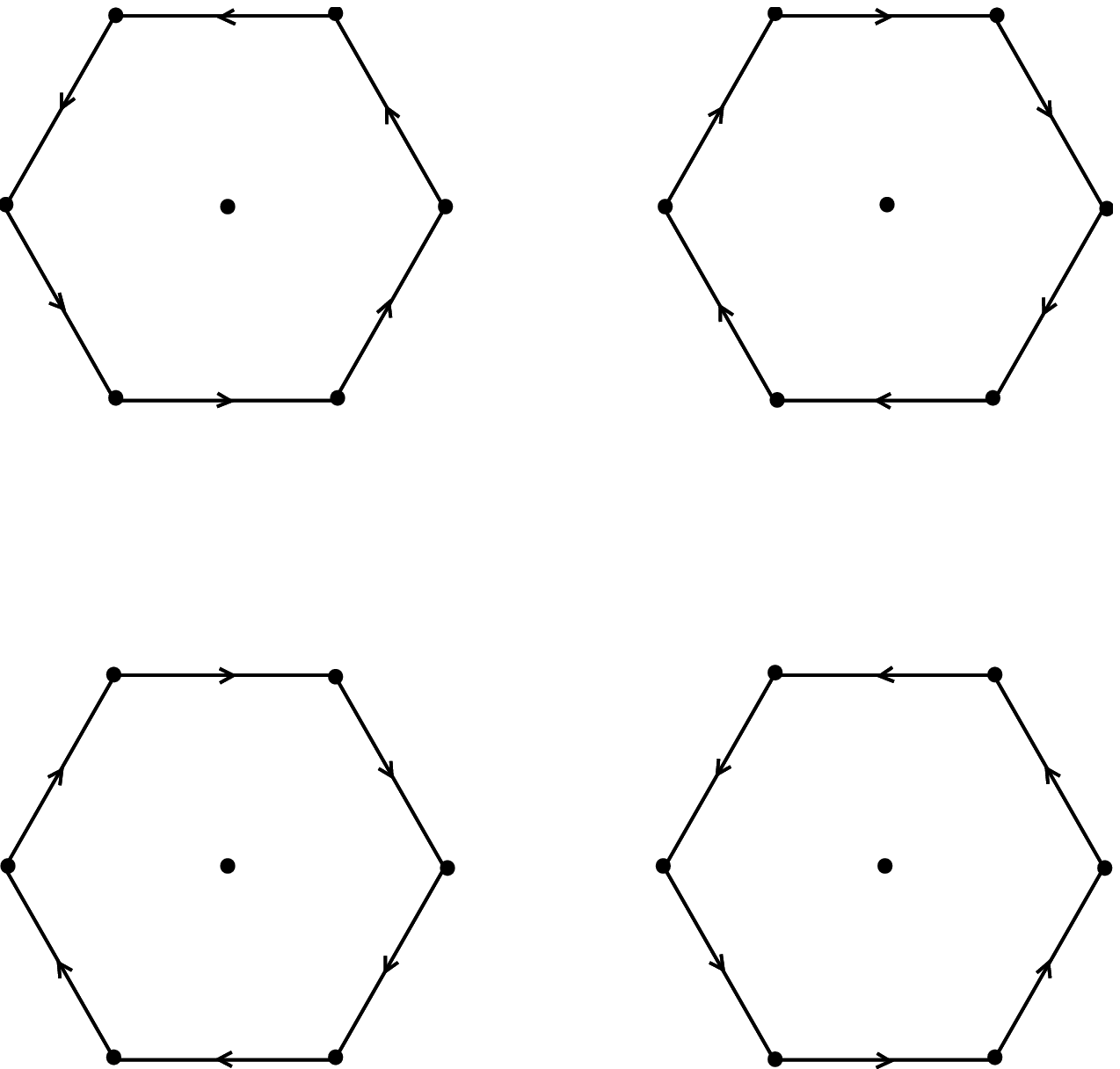}}
\put(8.5,13.5){\small $P_1$}
\put(10.5,17){\small $P_2$}
\put(8.5,21.2){\small $P_3$}
\put(3,21.4){\small $P_4$}
\put(0.8,17.6){\small $P_5$}
\put(3,13.5){\small $P_6$}
\put(9.2,15.5){\small $a_1$}
\put(9.4,19.3){\small $b_1$}
\put(5.5,21.5){\small $a_2$}
\put(2,19.5){\small $b_2$}
\put(2,15.5){\small $a_3$}
\put(5.5,13.2){\small $b_3$}
\put(5.5,18.2){\small $Q_1$}
\put(22.3,17.2){\small $P_1$}
\put(20.3,13.4){\small $P_2$}
\put(15.1,13.2){\small $P_3$}
\put(12.8,16.8){\small $P_4$}
\put(15,21.3){\small $P_5$}
\put(20,21.3){\small $P_6$}
\put(21.2,15.2){\small $a_1$}
\put(17.5,13.3){\small $c_1$}
\put(14,15.3){\small $a_2$}
\put(13.8,19.4){\small $c_2$}
\put(17.5,21.5){\small $a_3$}
\put(21,19.6){\small $c_3$}
\put(17.5,18.1){\small $Q_2$}
\put(10.3,5){\small $P_1$}
\put(8,1.3){\small $P_2$}
\put(3.5,1.3){\small $P_3$}
\put(1,5.2){\small $P_4$}
\put(3,9.4){\small $P_5$}
\put(8,9.4){\small $P_6$}
\put(9.3,3.6){\small $d_1$}
\put(6,1.2){\small $b_1$}
\put(2,3.5){\small $d_2$}
\put(1.9,7.5){\small $b_2$}
\put(5.5,9.7){\small $d_3$}
\put(9.2,7.5){\small $b_3$}
\put(5.7,6.2){\small $Q_3$}
\put(20,1.2){\small $P_1$}
\put(22.2,5){\small $P_2$}
\put(20.3,9){\small $P_3$}
\put(14.9,9.5){\small $P_4$}
\put(12.7,5.4){\small $P_5$}
\put(15.5,1.2){\small $P_6$}
\put(21.2,3.5){\small $d_1$}
\put(21.2,7.4){\small $c_1$}
\put(17.5,9.5){\small $d_2$}
\put(13.8,7.7){\small $c_2$}
\put(13.7,3.5){\small $d_3$}
\put(17.5,1.5){\small $c_3$}
\put(17.5,6.3){\small $Q_4$}
\end{picture}}
\centerline{{\bf Figure 5.1.} A cellular decomposition of the genus 2 oriented surface.}
\end{figure}

\bigskip\noindent
Consider the isometries $R,S,T \in \Isom (\Sigma)$ defined as follows. The transformation $R$ is an 
order $\rho+1$ isometry such that $R(A_j)=A_j$ for all $1 \le j \le 4$, and whose restriction to each 
$A_j$ is a rotation of angle ${ 2 \pi \over \rho +1}$ (resp. $-{ 2 \pi \over \rho +1}$)
centered at $Q_j$ for $j=1,4$ (resp. for $j=2,3$), where $Q_j$ denotes the center 
of $A_j$. Note that $R$ permutes cyclically $a_1, a_2, \dots, a_{\rho+1}$, as well as $b_1, b_2, \dots, 
b_{\rho+1}$, and $c_1, c_2, \dots, c_{\rho+1}$, and $d_1, d_2, \dots, d_{\rho+1}$. The isometry $S$ is 
the (unique) involution which sends $a_1$ onto $a_1^{-1}$, $d_1$ onto $d_1^{-1}$, $A_1$ onto $A_2$, and 
$A_4$ onto $A_3$. It fixes the centers of $a_1$ and $d_1$, and may fix two other points, depending on 
the evenness of $\rho+1$. The isometry $T$ is the (unique) involution which sends $a_i$ onto $d_i$, $b_i$ 
onto $c_i$, for all $1 \le i\le \rho+1$, $A_1$ onto $A_4$, and $A_2$ onto $A_3$. It fixes all the 
$P_j$'s. Let $\Gamma$ be the subgroup of $\Isom (\Sigma)$ generated by $\{ R,S,T \}$. It is easily 
checked that $\Gamma$ has the presentation
\[
\Gamma= \langle R,S,T \ |\ R^{\rho+1}= S^2= T^2 =1,\ S\circ R\circ S=R^{-1},\ T\circ R=R\circ T,\ T \circ 
S = S \circ T \rangle\,,
\]
and is isomorphic to $\DD_{2(\rho+1)} \times \CC_2$, where $\DD_{2(\rho+1)}$ denotes the dihedral group 
of order $2(\rho+1)$ and $\CC_2= \{ \pm 1\}$.

\bigskip\noindent
Let $\Sigma_0= \Sigma / \Gamma$ and let $\pi: \Sigma \to \Sigma_0$ be the natural quotient. Then 
$\Sigma_0$ is a sphere and $\pi$ is a ramified covering with 4 ramification points, $\bar P_0, \bar P_1, 
\bar P_2, \bar P_3$. Up to permutation, we can assume that $\pi^{-1} (\bar P_0)= \{ Q_1, Q_2, Q_3, Q_4 
\}$, $\pi^{-1} (\bar P_1)= \{ P_1, P_2, \dots, P_{\rho+1} \}$, $\pi^{-1} (\bar P_2)$ is the set of the 
centers of the $a_i$'s and of the $d_i$'s, and $\pi^{-1} (\bar P_3)$ is the set of the centers of the 
$b_i$'s and of the $c_i$'s. In particular, we have $r_\pi (\bar P_0)= \rho+1$ and $r_\pi (\bar P_1)= 
r_\pi (\bar P_2)= r_\pi (\bar P_3) = 2$.

\bigskip\noindent
Set $\PP_0= \{ \bar P_0, \bar P_1, \bar P_2, \bar P_3 \}$ and take a pseudo-Anosov diffeomorphism $F_0 
\in \Diff (\Sigma_0, \PP_0)$. Such a diffeomorphism can be easily constructed as follows. Choose a 
matrix $M \in SL_2 ( \Z)$ having two distinct real eigenvalues $\lambda$ and ${1 \over \lambda}$, with 
$| \lambda | >1$, and denote by $F_M$ the Anosov diffeomorphism of the torus $\T^2= \R^2/ \Z^2$ induced 
by $M$. Let $F_{-\Id}$ be the involution of $\T^2$ induced by $-\Id \in SL_2 ( \Z)$. Then $F_{-\Id}$ 
commutes with $F_M$, thus $F_M$ induces a diffeomorphism $F_0$ on $\Sigma_0= \T^2 / F_{-\Id}$. 
Moreover, the projection $\T^2 \to \T^2 / F_{-\Id} = \Sigma_0$ is a 2-folds ramified covering of 
$\Sigma_0$ with 4 ramification points that we can assume to be $\bar P_0, \bar P_1, \bar P_2, \bar 
P_3$, thus $F_0$ is a pseudo-Anosov diffeomorphism of $(\Sigma_0, \PP_0)$.

\bigskip\noindent
Let $\FF_0^s$ (resp. $\FF_0^u$) denote the stable (resp. unstable) foliation of $F_0$, and let $\FF^s$ 
(resp. $\FF^u$) denote the lift of $\FF_0^s$ (resp. $\FF_0^u$) on $\Sigma$. Upon replacing $F_0$ by some 
power, we can assume that $F_0$ lifts to a pseudo-Anosov diffeomorphism $F: \Sigma \to \Sigma$ whose 
stable (resp. unstable) foliation is $\FF^s$ (resp. $\FF^u$), and which 
fixes its separatrices (see Subsections 2.5 and 2.6).

\bigskip\noindent
It is easily checked from Proposition 2.2 that $\Ind (\FF_0^s, \FF_0^u: \bar P_i)= 1$ for all $0 \le 
i\le 4$, and that the singularities of $(\FF_0^s, \FF_0^u)$ are precisely $\bar P_0, \bar P_1, \bar 
P_2, \bar P_3$. On the other hand, we already know that $r_\pi (\bar P_0)= \rho +1$ and $r_\pi (\bar 
P_i)= 2$ for all $1 \le i\le 3$. So, $\bar P_0$ is a pivot. By Proposition 2.18, it follows that $\Gamma= 
\Sym (\FF^s, \FF^u)$.

\bigskip\noindent
Set $F_1 = F \circ R$ and $F_2= F \circ S$, and denote by $f_1 \in \MM (\Sigma)$ (resp. $f_2 \in \MM 
(\Sigma)$) the isotopy class of $F_1$ (resp. $F_2$). By Lemma 2.19, $R$ and $S$ commute with $F$, thus 
$F_1^{2(\rho+1)} = F^{2(\rho +1)} \circ R^{ 2(\rho+1)} = F^{2(\rho+1)} = F^{2(\rho+1)} \circ 
S^{2(\rho+1)} = F_2^{2(\rho+1)}$, therefore $f_1^{2(\rho+1)} = f_2^{2(\rho+1)}$.

\bigskip\noindent
Suppose that $f_1$ is conjugate to $f_2$. By Proposition 2.20, there exists $H \in \Norm (\FF^s, \FF^u)$ such 
that $H \circ R \circ H^{-1} = S$. But this equality is impossible because $R$ is of order $\rho+1$ while 
$S$ is of order 2.

\bigskip\noindent
Suppose that there exists a conjugate of $f_1$ which commutes with $f_2$. By Proposition 2.21, 
there exists $H \in \NNorm (\FF^s, \FF^u)$ such that $H \circ R \circ H^{-1}$ commutes with $S$. 
Recall that $\Gamma$ is isomorphic to $\DD_{2 (\rho +1)} \times \CC_2$. Moreover, $H \circ R \circ H^{-
1}$ is an element of $\Sym (\FF^s, \FF^u) \simeq \Gamma$ (by Lemma 2.7), and it is of order $\rho +1$ 
(because $R$ is of order $\rho +1$). A direct inspection on the group $\Gamma \simeq \DD_{2 (\rho+1)} 
\times \CC_2$ shows that the only elements of $\Gamma$ of order $\rho+1$ are:
\begin{itemize}
\item
$R^k$, with $k$ and $\rho+1$ coprime,
\item
$R^k \circ T$, with $\rho+1$ even, and $k$ and $\rho+1$ coprime,
\item
$R^k \circ T$, with $\rho+1 \equiv 2\ (\text{mod}\, 4)$, and $k$ and ${\rho+1 \over 2}$ coprime.
\end{itemize}
It is easily checked that none of these elements commutes with $S$: a contradiction.
\qed

\bigskip\noindent
The following lemma gives a simple way for constructing pseudo-Anosov braids. It is also a preliminary 
to the proof of Theorem 5.2.2.

\bigskip\noindent
{\bf Lemma 5.3.} {\it Let $\PP= \{ P_1, \dots, P_n\}$ be a set of $n$ punctures in the sphere $\S^2$, 
where $n$ is a prime number, and let $f \in \MM (\S^2, \PP)$ be a non-periodic element which permutes 
cyclically $P_1, P_2, \dots, P_n$. Then $f$ is a pseudo-Anosov element.}

\bigskip\noindent
{\bf Proof.} Let $f \in \MM (\S^2, \PP)$ which is neither periodic, nor pseudo-Anosov (so, is 
reducible), and assume that it acts cyclically on $P_1, P_2, \dots, P_n$. The hypothesis that $f$ is 
reducible means that $\Delta (f) \neq \emptyset$. Set $\Delta (f) = \{ \delta_0, \delta_1, \dots, 
\delta_p\}$. We choose essential curves $d_0, d_1, \dots, d_p$ such that $\langle d_i \rangle = 
\delta_i$ for all $0 \le i\le p$, and $d_i \cap d_j = \emptyset$ for all $0 \le i \neq j \le p$. We 
denote by $\Sigma_\Delta$ the surface obtained from $\S^2$ cutting along the $d_i$'s, and by 
$\pi_\Delta : \Sigma_\Delta \to \Sigma$ the natural quotient. Let $\Sigma_{\Delta\,1}, \dots, 
\Sigma_{\Delta\,l}$ be the connected components of $\Sigma_\Delta$, $\PP_\Delta = \pi_\Delta^{-
1}(\PP)$, $\PP_{\Delta\,i}= \PP_\Delta \cap \Sigma_{\Delta\,i}$, and $\PP_i= \pi_\Delta 
(\PP_{\Delta\,i})$ for $1 \le i\le l$. For all $1 \le i\le l$ there exists $1 \le j\le l$ such that 
$f(\Sigma_{\Delta\,i})$ is $\Sigma_{\Delta\,j}$ up to isotopy, and, therefore, $f(\PP_i) = \PP_j$. 
Suppose $\PP_i \neq \emptyset$ for $1 \le i \le k$ and $\PP_i = \emptyset$ for $k+1 \le i\le l$. The 
fact that $\S^2$ is a sphere implies that $k \ge 2$ and that $| \PP_i| \ge 2$ for some $1 \le i\le k$ 
(say $|\PP_1| \ge 2$), and the fact that $f$ acts transitively on $\PP$ implies that $|\PP_i| = 
|\PP_j|$ for all $1 \le i,j \le k$. It follows that $n= k \cdot |\PP_1|$, contradicting the hypothesis that 
$n$ is a prime number.
\qed

\bigskip\noindent
{\bf Proof of Theorem 5.2.2.} We assume that $\rho= 4l$, where $l \ge 1$, and we consider a closed 
oriented surface $\Sigma$ of genus $\rho$. We represent $\Sigma$ embedded in $\R^3= \C \times \R$ as 
follows (see Figure 5.2). Let $\gamma_k$ be the segment in $\C= \C \times \{ 0\} \subset \R^3$ which 
joins $0$ to $K \cdot \exp ({2ik\pi \over \rho})$, where $K$ is a fixed real number greater than $2/ 
\sin ({\pi \over \rho})$, and let $S_k$ be the circle in $\C$ centered at $(K+1) \cdot \exp ({2ik \pi 
\over \rho})$ of radius 1, for all $0 \le k\le 2 \rho -1$. Set
\[
\Upsilon = \left( \cup_{k=0}^{2\rho-1} \gamma_k \right) \cup \left( \cup_{k=0}^{2\rho-1} S_k 
\right)\,.
\]
Note that $\Upsilon$ is a connected graph embedded in $\C \subset \R^3$. We define $\Sigma$ as the 
boundary of a regular neighborhood of $\Upsilon$.

\begin{figure}[htbp]
\centerline{
\setlength{\unitlength}{.4cm}
\begin{picture}(20,12)
\put(1,1){\includegraphics[width=7.2 cm]{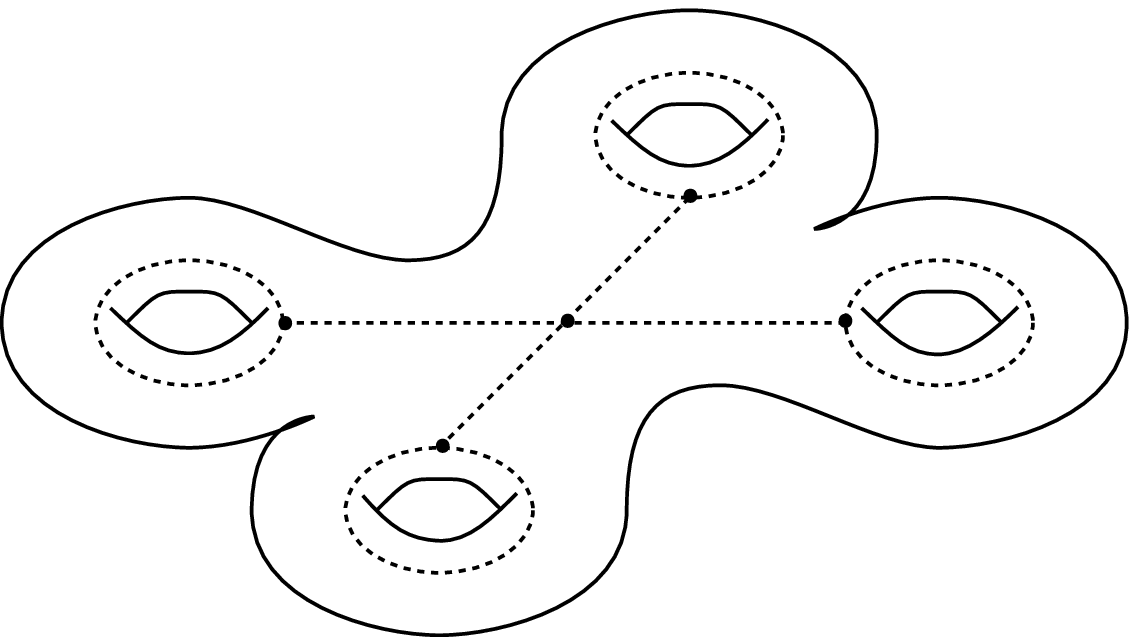}}
\put(17.8,6){\small $S_1$}
\put(11.5,10.3){\small $S_2$}
\put(1.5,5.5){\small $S_3$}
\put(7.5,1.3){\small $S_4$}
\put(12,6.3){\small $\gamma_1$}
\put(10,7){\small $\gamma_2$}
\put(7.5,6.3){\small $\gamma_3$}
\put(7.8,5){\small $\gamma_4$}
\end{picture}}
\centerline{{\bf Figure 5.2.} A surface of genus 4.}
\end{figure}

\bigskip\noindent
Consider the isometries $R,S \in \Isom (\Sigma)$ defined as follows. Let $\{ \vec e_1, \vec e_2, \vec 
e_3 \}$ be the standard basis of $\R^3$. Then $R$ is the rotation of angle ${ 2 \pi \over \rho}$ around 
the axis $\R \vec e_3$, and $S$ is the half-turn around the axis $\R \vec e_1$. Let $\Gamma$ be the 
group generated by $\{ R,S \}$. It is easily checked that $\Gamma$ has the presentation
\[
\Gamma= \langle R,S\ |\ R^\rho= S^2=1\,, S\circ R\circ S=R^{-1} \rangle\,,
\]
and is isomorphic to the dihedral group $\DD_{2 \rho}$. Note that $R^k\circ S$ is a half-turn with 
6 fixed points on $\Sigma$ if $k$ is even, $R^k \circ S$ has 2 fixed points on $\Sigma$ if 
$k$ is odd, and $R^k$ has 2 fixed points on $\Sigma$ if $1 \le k \le \rho-1$.
In particular, the elements of $\Gamma$ having precisely 6 fixed points are those of the 
form $R^k \circ S$ with $k$ even.

\bigskip\noindent
Let $\Sigma_0= \Sigma/\Gamma$, and let $\pi: \Sigma \to \Sigma_0$ be the natural quotient. Then 
$\Sigma_0$ is a sphere and $\pi$ is a ramified covering with 5 ramification points, $\bar P_0, \bar 
P_1, \bar P_2, \bar P_3, \bar P_4$. Up to a permutation, we céan assume that $\pi^{-1} (\bar P_0)= 
\Sigma \cap \R \vec e_3$ (which has 2 points) and, consequently, that $r_\pi(\bar P_0)= \rho$. Note 
that $r_\pi(\bar P_i)=2$ for the remaining indices $1 \le i \le 4$.

\bigskip\noindent
Set $\PP_0= \{ \bar P_0, \bar P_1, \bar P_2, \bar P_3, \bar P_4 \}$. Let $f_0 \in \MM (\Sigma_0, 
\PP_0)$ be a non-periodic element which permutes cyclically $\bar P_0, \bar P_1, \dots, \bar P_4$. 
By Lemma 5.3,
$f_0$ is a pseudo-Anosov element. Let $F_0 \in \Diff (\Sigma_0, \PP_0)$ be a pseudo-Anosov 
diffeomorphism which represents $f_0$, and let $\FF_0^s$ (resp. $\FF_0^u$) be the stable (resp. 
unstable) foliation of $F_0$. Since $F_0$ permutes cyclically $\bar P_0, \bar P_1, \dots, \bar P_4$, we 
have $\Ind (\FF_0^s, \FF_0^u: \bar P_i)= \Ind (\FF_0^s, \FF_0^u: \bar P_0)$ for all $1 \le i\le 4$. On 
the other hand, by Proposition 2.2, one of these indices must be 1, thus $\Ind (\FF_0^s, \FF_0^u: \bar 
P_i) = 1$ for all $0 \le i\le 4$. Let $\SS_0$ denote the set of singularities of $(\FF_0^s, \FF_0^u)$. 
Then, again by Proposition 2.2, there exists $\bar P_5 \in \Sigma_0 \setminus \PP_0$ such that $\SS_0 = 
\PP_0 \cup \{ \bar P_5 \}$ and $\Ind( \FF_0^s, \FF_0^u: \bar P_5)=3$.

\bigskip\noindent
Let $\FF^s$ (resp. $\FF^u$) be the lift of $\FF_0^s$ (resp. $\FF_0^u$) on $\Sigma$. By Lemma 2.16, there 
exists some $k \ge 1$ such that $F_0^k$ lifts to a pseudo-Anosov diffeomorphism $F: \Sigma \to \Sigma$ 
whose stable (resp. unstable) foliation is $\FF^s$ (resp. $\FF^u$). Obviously, $k$ can be chosen so 
that $F$ fixes its separatrices. Recall that 
$r_\pi(\bar P_0)=\rho$ and $\Ind (\FF_0^s, \FF_0^u: \bar P_0)=1$, $r_\pi (\bar P_i) =2$ and $\Ind 
(\FF_0^s, \FF_0^u: \bar P_i)=1$ for all $1 \le i\le 4$, and $r_\pi(\bar P_5)= 1$ and $\Ind (\FF_0^s, 
\FF_0^u: \bar P_5)=3$, thus $\bar P_0$ is a pivot. By Lemma 2.18, it follows that $\Gamma= \Sym 
(\FF^s, \FF^u)$.

\bigskip\noindent
Set $F_1= F \circ S$ and $F_2= F \circ R \circ S$, and denote by $f_1 \in \MM(\Sigma)$ (resp. $f_2 \in 
\MM (\Sigma)$) the isotopy class of $F_1$ (resp. $F_2$). By Lemma 2.19, $R$ and $S$ commute with $F$, 
thus $F_1^2= F^2= F_2^2$, therefore $f_1^2=f_2^2$.

\bigskip\noindent
Suppose that $f_1$ and $f_2$ are conjugate in $\MM(\Sigma)$. By Proposition 2.20, there exists $H \in \Norm 
(\FF^s, \FF^u)$ such that $H \circ S \circ H^{-1} = R \circ S$. But this equality is impossible because $S$ 
has 6 fixed points on $\Sigma$ while $R \circ S$ has only 2.

\bigskip\noindent
Suppose that there exists a conjugate of $f_1$ which commutes with $f_2$. By Proposition 2.21, there exists $H 
\in  \NNorm (\FF^s, \FF^u)$ such that $H \circ S \circ H^{-1}$ commutes with $R \circ S$. Observe 
that $H \circ S \circ H^{-1}$ should be an order 2 element of $\Sym( \FF^s, \FF^u)= \Gamma$ which 
fixes 6 points on $\Sigma$,
thus, as pointed out before, it is of the form $R^k \circ S$ with $k$ even. 
On the other hand, the only elements of order 2 in $\Gamma \simeq \DD_{2 \rho}$ which commute with 
$R \circ S$ are $R \circ S$, $R^{{\rho-2 \over 2}} \circ S$, and $R^{{\rho \over 2}}$, and all of them 
have precisely 2 fixed points on $\Sigma$ (note that ${\rho -2 \over 2}$ is odd). This is a 
contradiction.
\qed


\section{Pure subgroups}

{\bf Theorem 6.1.} {\it Let $\Sigma$ be a surface of genus $\rho \ge 0$ and $q$ boundary components, 
$c_1, \dots, c_q$, and let $\PP= \{ P_1, \dots, P_n\}$ be a finite set of punctures in the interior of 
$\Sigma$. Let $\Gamma$ be a pure subgroup of $\MM (\Sigma, \PP)$. If $f,g \in \Gamma$ are such that 
$f^m = g^m$ for some $m \ge 1$, then $f=g$.}

\bigskip\noindent
{\bf Proof.} We consider two elements $f,g \in \Gamma$ such that $f^m= g^m$ for some $m \ge 1$. 
Obviously, we have $\Delta (f) = \Delta (g)$ and $f,g$ are either both periodic, or both pseudo-Anosov, 
or both reducible.

\bigskip\noindent
Assume that $f,g$ are both periodic. Let $K$ be the subgroup of $\MM (\Sigma, \PP)$ generated by $\{ 
\tau_{c_1}, \dots, \tau_{c_q}\}$. The only periodic pure elements of $\MM (\Sigma, \PP)$ are the 
elements of $K$, and $K \simeq \Z^q$, thus $f,g \in K$, and the equality $f^m=g^m$ clearly implies $f=g$.

\bigskip\noindent
Now, assume that $f,g$ are both pseudo-Anosov. Let $\theta: \MM (\Sigma, \PP) \to \MM (\Sigma_0, \PP 
\sqcup \QQ)$ be the corking of $(\Sigma, \PP)$. By Corollary 2.14, there exist pseudo-Anosov 
representatives $F,G \in \Diff (\Sigma_0, \PP \sqcup \QQ)$ of $\theta(f), \theta(g)$, respectively, 
such that $F^m=G^m$. Let $\FF^s$ (resp. $\FF^u$) be the stable (resp. unstable) foliation of $F$, and 
let $\lambda$ be its dilatation coefficient. Then $\FF^s$ (resp. $\FF^u$) is also the stable (resp. 
unstable) foliation of $G$, and $\lambda$ is its dilatation coefficient. It follows that $GF^{-1} \in 
\Sym (\FF^s, \FF^u)$, thus $GF^{-1}$ has finite order, therefore $gf^{-1}$ is periodic (since $GF^{-
1}$ represents $\theta (gf^{-1})$). As pointed out before, the periodic elements of $\Gamma$ are 
precisely the elements of $\Gamma \cap K$, thus $g f^{-1} \in K$. Set $h=gf^{-1}$. The subgroup $K$ is 
contained in the center of $\MM (\Sigma, \PP)$, thus
\[
g^m= (hf)^m = h^m f^m = h^m g^m\,,
\]
therefore $h^m=1$. We conclude that $h=1$ (since $K$ is torsion free), that is, $g=f$.

\bigskip\noindent
Now, assume that $f,g$ are both reducible. We use the same notations as in Subsection 2.7. So, we set 
$\Delta= \Delta (f)= \Delta (g)= \{ \delta_0, \delta_1, \dots, \delta_p \}$, and we choose essential 
curves $d_0, d_1, \dots, d_p$ of $(\Sigma, \PP)$ such that $\langle d_i \rangle = \delta_i$ for all $0 
\le i\le p$, and $d_i \cap d_j = \emptyset$ for all $0 \le i \neq j \le p$. We denote by 
$\Sigma_\Delta$ the surface obtained from $\Sigma$ cutting along the $d_i$'s. Thus, $\Sigma_\Delta$ is 
a non-necessarily connected surface whose boundary components are $c_1, \dots, c_q, d_0^{(1)}, 
d_1^{(1)}, \dots, d_p^{(1)}, d_0^{(2)}, d_1^{(2)}, \dots, d_p^{(2)}$ and $\Sigma= \Sigma_\Delta /\sim$, 
where $\sim$ is the equivalence relation which identifies $d_i^{(1)}(z)$ with $d_i^{(2)}(\bar z)
= (d_i^{(2)})^{-1}(z)$ for 
all $0 \le i\le p$ and all $z \in \S^1$. Let $\pi_\Delta: \Sigma_\Delta \to \Sigma$ denote the natural 
quotient. Then $d_i= \pi_\Delta \circ d_i^{(1)} = \pi_\Delta \circ (d_i^{(2)})^{-1}$. We denote by 
$\Sigma_{\Delta\,1}, \dots, \Sigma_{\Delta\,l}$ the connected components of $\Sigma_\Delta$, we set 
$\PP_\Delta= \pi_\Delta^{-1} (\PP)$, and we set $\PP_{\Delta\,k}= \PP_\Delta \cap \Sigma_{\Delta\,k}$ 
for all $1 \le k\le l$. We set $\MM (\Sigma_\Delta, \PP_\Delta)= \MM (\Sigma_{\Delta\,1}, 
\PP_{\Delta\,1}) \times \dots \times \MM (\Sigma_{\Delta\,l}, \PP_{\Delta\,l})$, and we denote by 
$\theta_\Delta: \MM (\Sigma_\Delta, \PP_\Delta) \to \MM (\Sigma, \PP)$ the homomorphism induced by 
$\pi_\Delta$.

\bigskip\noindent
Recall from Subsection 2.7 that, if $h \in \Stab (\Delta)$ and $h$ is pure, then it belongs to $\Im 
\theta_\Delta$. Moreover, if $h= \theta_\Delta(h_1, \dots, h_l)$, then each $h_k$ is a pure element of 
$\MM (\Sigma_{\Delta\, k}, \PP_{\Delta\,k})$. This observation is important in the present proof.

\bigskip\noindent
Let $\tilde K$ be the subgroup of $\MM (\Sigma_\Delta, \PP_\Delta)$ generated by $\{ \tau_{c_1}, \dots, 
\tau_{c_q}, \tau_{d_0^{(1)}}, \tau_{d_1^{(1)}}, \dots, \tau_{d_p^{(1)}}, 
\linebreak
\tau_{d_0^{(2)}}, 
\tau_{d_1^{(2)}}, \dots, \tau_{d_p^{(2)}} \}$, and let $\tilde K_k = \tilde K \cap \MM 
(\Sigma_{\Delta\,k}, \PP_{\Delta\,k})$ for $1 \le k\le l$. Then $\tilde K_k$ is the subgroup of $\MM 
(\Sigma_{\Delta\,k}, \PP_{\Delta\,k})$ generated by the Dehn twists along the boundary components of 
$\Sigma_{\Delta\,k}$, and we have $\tilde K= \tilde K_1 \times \dots \times \tilde K_l$.

\bigskip\noindent
Take $(f_1, \dots, f_l) \in \theta_\Delta^{-1} (f)$ and $(g_1, \dots, g_l) \in \theta_\Delta^{-1} (g)$. 
The equality $\Delta= \Delta(f)$ implies that each $f_k$ is either a pseudo-Anosov element of $\MM 
(\Sigma_{\Delta\,k}, \PP_{\Delta\,k})$, or 
a periodic element (namely, and element of $\tilde K_k$).
Similarly, each $g_k$ is either 
a pseudo-Anosov element of $\MM (\Sigma_{\Delta\,k}, \PP_{\Delta\,k})$, or an element of $\tilde K_k$. 
Furthermore, the equality $f^m=g^m$ implies that $f_k^m \equiv g_k^m\ (\mod\, \tilde K_k)$, thus either 
$f_k,g_k$ are both pseudo-Anosov elements of $\MM (\Sigma_{\Delta\,k}, \PP_{\Delta\,k})$, or they both 
belong to $\tilde K_k$. Let $h= g f^{-1}$, and, for $1 \le k \le l$, let $h_k = g_k f_k^{-1}$. Then $h 
\in \Stab(\Delta)$ and $\theta_\Delta (h_1, \dots, h_l) =h$. Moreover, since $h$ is pure, each $h_k$ is 
a pure element of $\MM(\Sigma_{\Delta\,k}, \PP_{\Delta\,k})$. Obviously, if $f_k,g_k \in \tilde K_k$, 
then $h_k = g_k f_k^{-1} \in \tilde K_k$.

\bigskip\noindent
Assume that $f_k, g_k$ are both pseudo-Anosov elements of $\MM (\Sigma_{\Delta\,k}, \PP_{\Delta\,k})$. 
Let $\mu_k: \MM (\Sigma_{\Delta\,k}, \PP_{\Delta\,k}) 
\linebreak
\to \MM (\Sigma_{0\,k}, \PP_{\Delta\,k} \sqcup 
\QQ_k)$ be the corking of $(\Sigma_{\Delta\,k}, \PP_{\Delta\,k})$. From the congruence $f_k^m \equiv 
g_k^m\ (\mod\, \tilde K_k)$ we deduce that $\mu_k (f_k)^m = \mu_k(g_k)^m$. By Corollary 2.14, there 
exist pseudo-Anosov representatives $F,G \in \Diff (\Sigma_{0\,k}, \PP_{\Delta\,k} \sqcup \QQ_k)$ of 
$\mu_k(f_k), \mu_k(g_k)$, respectively, such that $F^m=G^m$. Let $\FF^s$ (resp. $\FF^u$) be the stable 
(resp. unstable) foliation of $F$, and let $\lambda$ be its dilatation coefficient. Then $\FF^s$ (resp. 
$\FF^u$) is the stable (resp. unstable) foliation of $G$, and $\lambda$ is its dilatation coefficient. 
It follows that $G F^{-1} \in \Sym (\FF^s, \FF^u)$, thus $GF^{-1}$ has finite order, therefore $h_k$ 
is periodic (since $GF^{-1}$ is a representative of $\mu_k(h_k)$). Finally, $h_k$ is pure and the only 
pure periodic elements of $\MM (\Sigma_{\Delta\,k}, \PP_{\Delta\,k})$ are the elements of $\tilde K_k$, 
thus $h_k \in \tilde K_k$.

\bigskip\noindent
Let $K_0$ be the subgroup of $\MM(\Sigma, \PP)$ generated by $\{ \tau_{c_1}, \dots, \tau_{c_q}, 
\tau_{d_0}, \tau_{d_1}, \dots, \tau_{d_p} \}$. Note that $K_0$ is a free abelian group of rank $p+q+1$ 
and is contained in the centralizer of $f$ in $\MM(\Sigma, \PP)$. The fact that $h_k \in \tilde K_k$ 
for all $1 \le k\le l$ implies that $h \in K_0$, thus
\[
g^m= (hf)^m = h^m f^m= h^m g^m\,,
\]
therefore $h^m=1$. We conclude that $h=1$ (since $K_0$ is torsion free), that is, $f=g$.
\qed


\section{Proof of Theorem 3.6}

Let $\Sigma$ be a surface of genus 1 with $q$ boundary components, $c_1, 
\dots, c_q$, where $q \ge 1$, and let $f,g \in \MM (\Sigma)$ such that $f^m=g^m$ for some $m \ge 1$. In 
particular, $\Delta (f)= \Delta(f^m) = \Delta(g^m)= \Delta(g)$. Set $\Delta= \{ \delta_0, \delta_1, 
\dots, \delta_p \} = \Delta(f)= \Delta(g)$, and choose essential curves $d_0, d_1, \dots, d_p : \S^1 
\to \Sigma$ such that $\langle d_i \rangle = \delta_i$ for all $0 \le i \le p$, and $d_i \cap d_j = 
\emptyset$ for all $0 \le i \neq j \le p$. We denote by $\Sigma_\Delta$ the surface obtained from 
$\Sigma$ cutting along the $d_i$'s. So, $\Sigma_\Delta$ is a compact (but not necessarily connected) 
surface whose boundary components are $c_1, \dots, c_q$, $d_0^{(1)}, \dots, d_p^{(1)}$, $d_0^{(2)}, 
\dots, d_p^{(2)}$, and $\Sigma= \Sigma_\Delta /\sim$, where $\sim$ is the equivalence relation which 
identifies $d_i^{(1)}(z)$ with $d_i^{(2)}(\bar z) = (d_i^{(2)})^{-1} (z)$ 
for all $0 \le i\le p$ and all $z \in \S^1$. We 
denote by $\Sigma_{\Delta\,1}, \dots, \Sigma_{\Delta\,l}$ the connected components of $\Sigma_\Delta$.

\bigskip\noindent
We define a graph $\tilde \Gamma$ as follows. The set $\{ v_1, \dots, v_l\}$ of vertices is in 
bijection with the set of connected components of $\Sigma_\Delta$, the set $\{ e_0, e_1, \dots, e_p\}$ of edges 
is in bijection with $\Delta= \{ \delta_0, \delta_1, \dots, \delta_p\}$, and an edge $e_i$ connects the 
vertices $v_j$ and $v_k$ if $d_i^{(1)}$ is a boundary component of $\Sigma_{\Delta\,j}$ and $d_i^{(2)}$ 
is a boundary component of $\Sigma_{\Delta\,k}$. In other words, $\tilde \Gamma$ is the dual of the 
decomposition of $\Sigma$ determined by $d_0, d_1, \dots, d_p$. For convenience, we will use the graph 
$\Gamma$ which extends $\tilde \Gamma$ as follows. We add $q$ vertices $v_1', \dots, v_q'$ and $q$ 
edges $e_1', \dots, e_q'$ to $\tilde \Gamma$, and we set that $e_i'$ connects $v_i'$ with $v_j$ if 
$c_i$ is a boundary component of $\Sigma_{\Delta\,j}$. So, $v_1', \dots, v_q'$ are valence 1 vertices 
of $\Gamma$. If $\Sigma_{\Delta\,j}$ is of genus $1$, then $v_j$ may have valence 1, 2, or more. 
Otherwise, $\Sigma_{\Delta\,j}$ is of genus $0$ and $v_j$ is of valence $\ge 3$ in $\Gamma$.

\bigskip\noindent
We consider four different possible configurations:

\bigskip\noindent
{\bf Case 1 :}
one of the components of $\Sigma_{\Delta}$, say $\Sigma_{\Delta\,1}$, is of genus 1 (see Figure 7.1).
In this case $\Gamma$ is a tree.

\begin{figure}[htbp]
\centerline{
\setlength{\unitlength}{.4cm}
\begin{picture}(21,11)
\put(1,2){\includegraphics[width=7.6 cm]{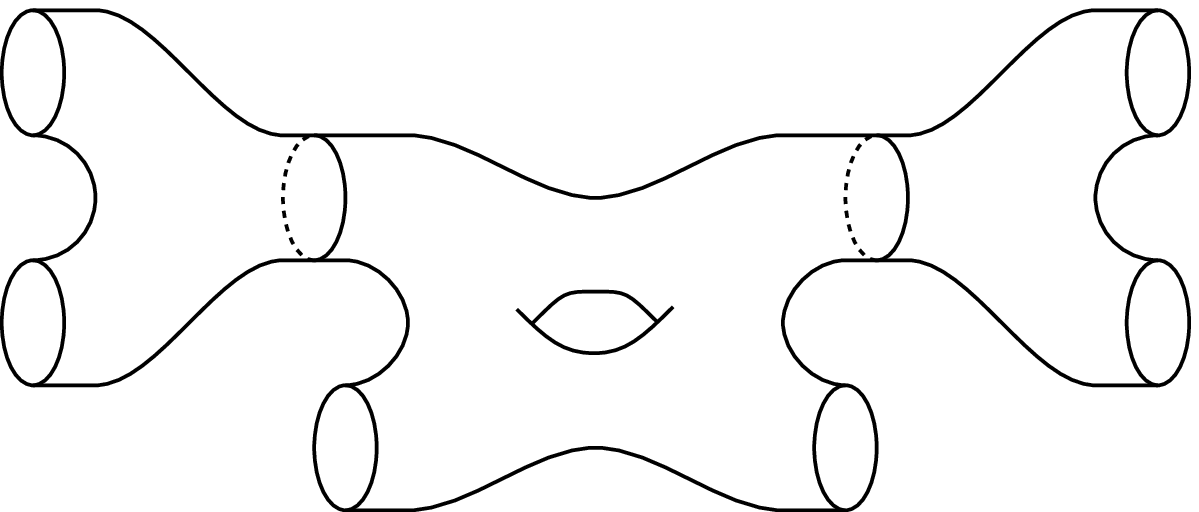}}
\put(5.7,8.5){\small $d_0$}
\put(14.4,8.5){\small $d_1$}
\put(10,1){$\Sigma_{\Delta\,1}$}
\put(3,3){$\Sigma_{\Delta\,2}$}
\put(16.5,3){$\Sigma_{\Delta\,3}$}
\end{picture}}
\centerline{{\bf Figure 7.1.} The first case.}
\end{figure}

\bigskip\noindent
If all the components of $\Sigma_\Delta$ are of genus $0$, then $\pi_1(\Gamma)$ is isomorphic to $\Z$ 
and $\Gamma$ has a unique reduced cycle. The remaining cases depend on the length of this reduced cycle (length 
1, or 2, or more). 

\bigskip\noindent
{\bf Case 2 :}
there exist $j \in \{ 0,1, \dots, p\}$ (say $j=0$) and a component of $\Sigma_\Delta$ (say 
$\Sigma_{\Delta\, 1}$) such that $d_0^{(1)}$ and $d_0^{(2)}$ are both boundary components of 
$\Sigma_{\Delta\, 1}$ (see Figure 7.2).
In this case $\Gamma$ has a unique reduced cycle which is of length 1. The vertex of this cycle is 
$v_1$ and the edge is $e_0$.

\begin{figure}[htbp]
\centerline{
\setlength{\unitlength}{.4cm}
\begin{picture}(21,11)
\put(1,2){\includegraphics[width=7.6 cm]{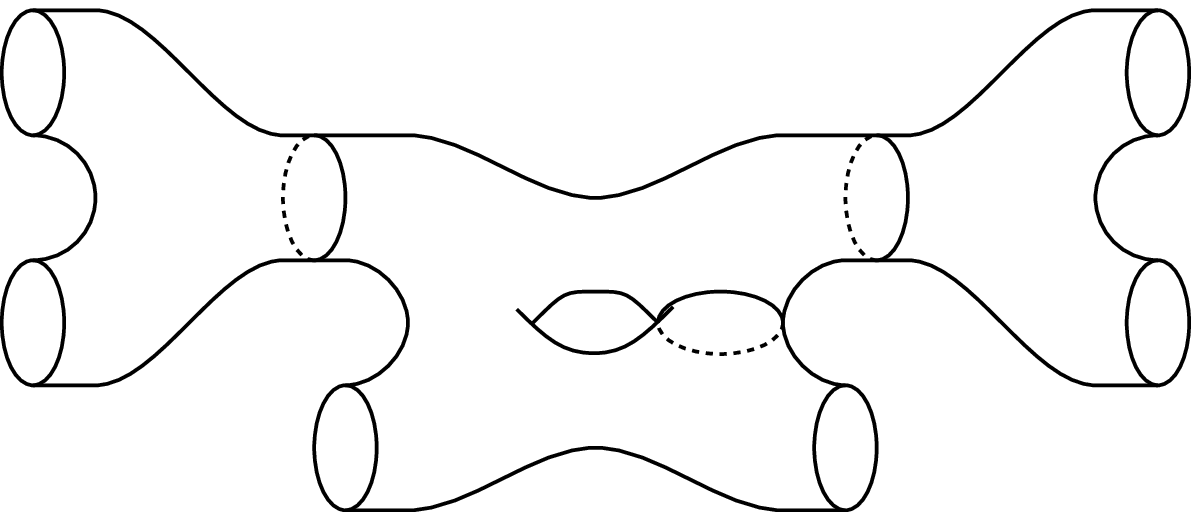}}
\put(13.8,4.9){\small $d_0$}
\put(5.7,8.5){\small $d_1$}
\put(14.4,8.5){\small $d_2$}
\put(10,1){$\Sigma_{\Delta\,1}$}
\put(3,3){$\Sigma_{\Delta\,2}$}
\put(16.5,3){$\Sigma_{\Delta\,3}$}
\end{picture}}
\centerline{{\bf Figure 7.2.} The second case.}
\end{figure}

\bigskip\noindent
{\bf Case 3 :}
there exist $j_1, j_2 \in \{ 0,1, \dots, p\}$ (say $j_1=0$ and $j_2=1$) and two components of 
$\Sigma_\Delta$ (say $\Sigma_{\Delta\,1}$ and $\Sigma_{\Delta\,2}$) such that $d_0^{(1)}$ and 
$d_1^{(1)}$ are boundary components of $\Sigma_{\Delta\, 1}$, and $d_0^{(2)}$ and $d_1^{(2)}$ are 
boundary components of $\Sigma_{\Delta\,2}$ (see Figure 7.3).
In this case $\Gamma$ has a unique reduced cycle which is of length 2. The vertices of this cycle are $v_1, v_2$ and the edges are $e_0, e_1$.

\begin{figure}[htbp]
\centerline{
\setlength{\unitlength}{.4cm}
\begin{picture}(21,11)
\put(1,2){\includegraphics[width=7.6 cm]{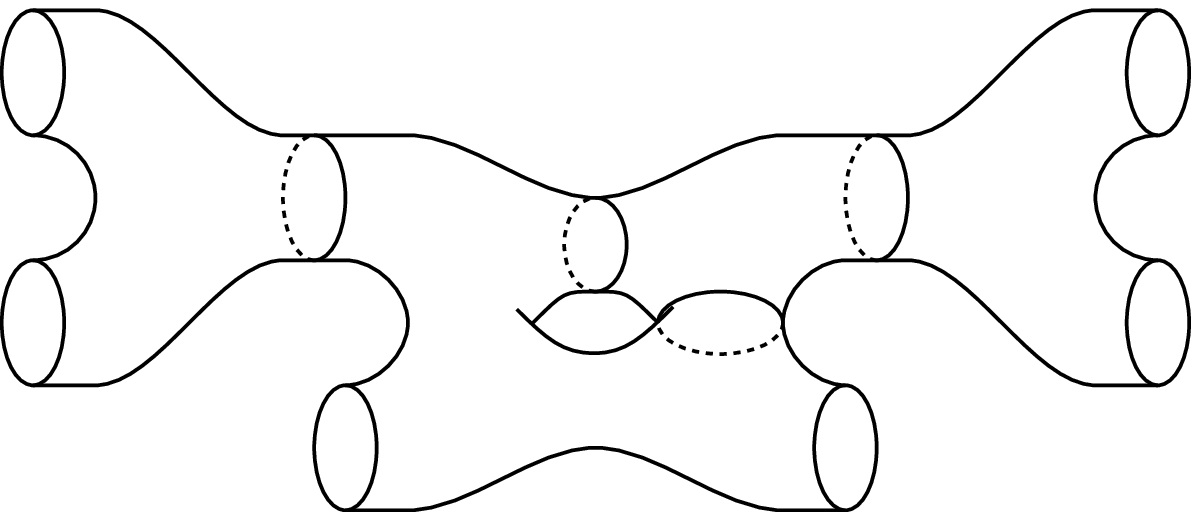}}
\put(13.8,4.9){\small $d_0$}
\put(10,7.5){\small $d_1$}
\put(5.7,8.5){\small $d_2$}
\put(14.4,8.5){\small $d_3$}
\put(10,1){$\Sigma_{\Delta\,1}$}
\put(12.2,6.4){$\Sigma_{\Delta\,2}$}
\put(3,3){$\Sigma_{\Delta\,3}$}
\put(16.5,3){$\Sigma_{\Delta\,4}$}
\end{picture}}
\centerline{{\bf Figure 7.3.} The third case.}
\end{figure}

\bigskip\noindent
{\bf Case 4 :}
none of the above (see Figure 7.4).
In this case $\Gamma$ has a unique reduced cycle which is of length $\ge 3$.

\begin{figure}[htbp]
\centerline{
\setlength{\unitlength}{.4cm}
\begin{picture}(21,11)
\put(1,2){\includegraphics[width=7.6 cm]{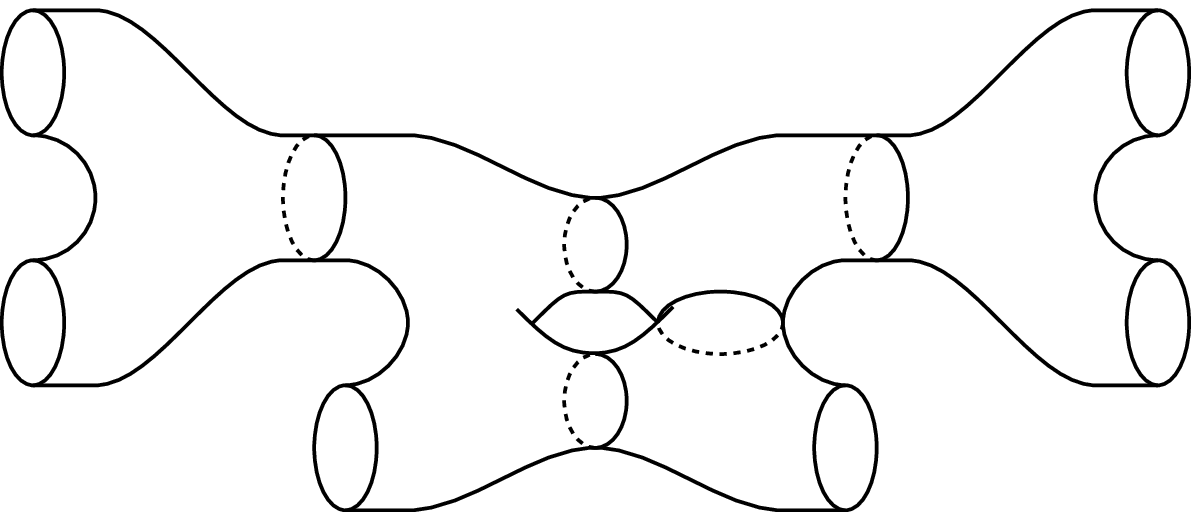}}
\put(10,2){\small $d_0$}
\put(13.8,4.9){\small $d_1$}
\put(10,7.5){\small $d_2$}
\put(5.7,8.5){\small $d_4$}
\put(14.4,8.5){\small $d_3$}
\end{picture}}
\centerline{{\bf Figure 7.4.} The fourth case.}
\end{figure}

\bigskip\noindent
The following lemmas are preliminaries. Lemma 7.2 is especially needed to treat Case 1, Lemma 7.3 is 
needed to treat Case 2, and Lemma 7.4 is needed to treat Case 3. Lemma 7.1 will be used in all the 
cases.

\bigskip\noindent
{\bf Lemma 7.1.} {\it Let $\Sigma$ be a surface of genus $0$ with $q$ boundary components, $c_1, \dots, 
c_q$, where $q \ge 3$. Let $K$ be the subgroup of $\MM (\Sigma)$ generated by $\{ \tau_{c_1}, \dots, 
\tau_{c_q}\}$. Let $f,g \in \MM(\Sigma)$ be two pseudo-Anosov elements. If $f^m \equiv g^m\ (\mod\, K)$ 
for some $m \ge 1$, then $f \equiv g\ (\mod\, K)$.}

\bigskip\noindent
{\bf Lemma 7.2.} {\it Let $\Sigma$ be a surface of genus 1 with 
$q$ boundary components, $c_1, \dots, c_q$, where $q \ge 1$.
Let $f,g \in \MM( 
\Sigma)$ be two periodic elements. If $f^m = g^m$ for some $m \ge 1$, then $f$ and $g$ are conjugate.}

\bigskip\noindent
{\bf Lemma 7.3.} {\it Let $\Sigma$ be a surface of genus 1 with $q$ boundary components, 
$c_1, \dots, c_q$, where $q \ge 1$,
and let $d: \S^1 \to \Sigma$ be a non-separating closed curve (see Figure 7.5). Set $\delta= \langle d 
\rangle$. Let $f,g \in \MM (\Sigma)$ such that $\Delta(f)= \Delta(g)= \{ \delta\}$. If $f^m=g^m$ for 
some $m \ge 1$, then $f$ and $g$ are conjugate.}

\begin{figure}[htbp]
\centerline{
\setlength{\unitlength}{.4cm}
\begin{picture}(11,8)
\put(1,1){\includegraphics[width=3.6 cm]{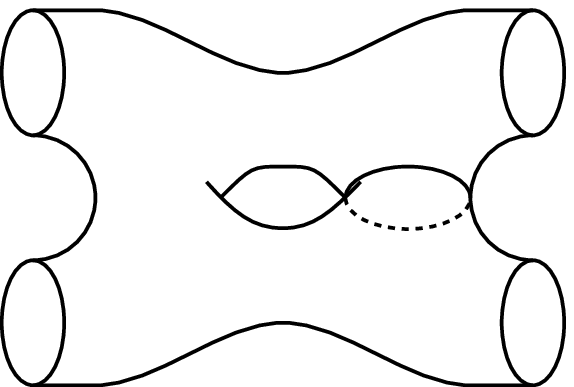}}
\put(8.8,3.8){\small $d$}
\end{picture}}
\centerline{{\bf Figure 7.5.} A non-separating curve in $\Sigma$.}
\end{figure}

\bigskip\noindent
{\bf Lemma 7.4.} {\it Let $\Sigma$ be a surface of genus 1 with $q$ boundary components, $c_1, \dots, 
c_q$, where $q \ge 2$. Let $d_1, d_2: \S^1 \to \Sigma$ be two non-separating curves such that $d_1$ is 
not isotopic to $d_2^{\pm 1}$, and $d_1 \cap d_2 = \emptyset$ (see Figure 7.6). Set $\delta_i= \langle 
d_i \rangle$ for $i \in \{ 1,2 \}$, and $\Delta= \{ \delta_1, \delta_2\}$. Let $f,g \in \MM(\Sigma)$ 
such that $\Delta(f)= \Delta(g)=\Delta$. If $f^m=g^m$ for some $m \ge 1$, then $f$ and $g$ are 
conjugate.}

\begin{figure}[htbp]
\centerline{
\setlength{\unitlength}{.4cm}
\begin{picture}(11,8)
\put(1,1){\includegraphics[width=3.6 cm]{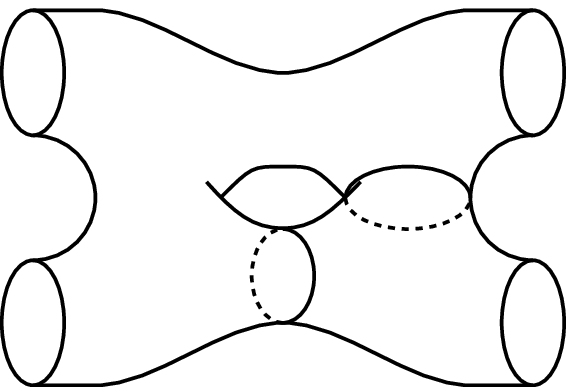}}
\put(5.2,1){\small $d_1$}
\put(8.8,3.8){\small $d_2$}
\end{picture}}
\centerline{{\bf Figure 7.6.} Two non-separating curves in $\Sigma$.}
\end{figure}

\bigskip\noindent
Now, we assume Lemmas 7.1--7.4 (they will be proved later), and turn to finish the proof of Theorem 3.6.

\bigskip\noindent
{\bf Proof of Theorem 3.6, Case 1.} We assume that $\Sigma_{\Delta\,1}$ is of genus 1 (see Figure 7.1).
For $2 \le i \le l$, the surface $\Sigma_{\Delta\,i}$ is of genus 0 and has at least 3 boundary 
components. Each $\Sigma_{\Delta\,i}$ embeds in $\Sigma$ ($1 \le i\le l$), and this embedding 
determines a homomorphism $\theta_{\Delta\,i}: \MM (\Sigma_{\Delta\,i}) \to \MM (\Sigma)$ which is 
injective (see \cite{ParRol1}). Let $\pi: \Sigma_\Delta \to \Sigma$ denote the quotient. We set 
$\MM(\Sigma_\Delta)= \MM( \Sigma_{\Delta\,1}) \times \MM( \Sigma_{\Delta\,2}) \times \dots \times \MM 
(\Sigma_{\Delta\,l})$, and we denote by $\theta_\Delta: \MM (\Sigma_\Delta) \to \MM(\Sigma)$ the 
homomorphism induced by $\pi$.

\bigskip\noindent
We denote by $\tilde K$ the subgroup of $\MM (\Sigma_\Delta)$ generated by $\{ \tau_{c_i}; 1 \le i \le 
q\} \cup \{ \tau_{d_i^{(1)}}, \tau_{d_i^{(2)}}; 0 \le i\le p\}$. This is a free abelian group of rank 
$q+2p+2$. Set $\tilde K_i= \tilde K \cap \MM (\Sigma_{\Delta\,i})$ for $1 \le i \le l$. Then, 
for all $1 \le i \le l$, $\tilde K_i$ is generated by the Dehn twists along the boundary components
of $\Sigma_{\Delta\,i}$, and
$\tilde K = \tilde K_1 \times \tilde K_2 \times \dots \times \tilde K_l$.

\bigskip\noindent
Since $\Delta (f) = \Delta$, the mapping class $f$ induces an isomorphism $f_\ast: \Gamma \to \Gamma$. 
Recall that $\Gamma$ is a tree and its leafs (i.e. the valence 1 vertices) are $v_1', 
\dots, v_q'$, and possibly $v_1$. The curves $c_1, \dots c_q$ are boundary components of $\Sigma$, thus 
there are fixed by $f$, therefore $f_\ast(v_i') = v_i'$ for all $1 \le i\le q$. It is easily proved that 
that this implies that $f_\ast$ is the identity on $\Gamma$.

\bigskip\noindent
The above shows that $f(\delta_k) = \delta_k$ and $f$ preserves the orientation of $d_k$ for all $0 \le 
k\le p$, and $\Sigma_{\Delta\,j}$ is invariant under $f$ up to isotopy for all $1 \le j \le l$. It 
follows that there exists $(f_1, f_2, \dots, f_l) \in \MM (\Sigma_\Delta)$ such that $f=\theta_\Delta 
(f_1, f_2, \dots, f_l)$. Furthermore, we have $\Delta(f_i) = \emptyset$ for all $1 \le i\le l$, since 
$\Delta(f) = \Delta = \{ \delta_0, \delta_1, \dots, \delta_p \}$.
Similarly, there exists $(g_1,g_2, \dots, g_l) \in \MM (\Sigma_\Delta)$ such 
that $g = \theta_\Delta (g_1,g_2, \dots, g_l)$ and $\Delta(g_i)= \emptyset$ for all $1 \le i\le l$. 

\bigskip\noindent
The equality $f^m= g^m$ implies that $f_i^m \equiv g_i^m\ (\mod\, \tilde K_i)$ for all $1 \le i\le 
l$. Let $i \in \{2, \dots, l\}$. Since $\Delta (f_i) = \Delta(g_i) = \emptyset$, the elements $f_i,g_i$ 
are either periodic or pseudo-Anosov. Obviously, the only periodic elements of 
$\MM(\Sigma_{\Delta\,i})$ are the elements of $\tilde K_i$. On the other hand, if both, $f_i$ and 
$g_i$, are pseudo-Anosov, then, by Lemma 7.1, $f_i \equiv g_i\ (\mod\, \tilde K_i)$. So, there exists 
$\tilde k \in \tilde K$ such that
\[
\theta_\Delta (1, g_2, \dots, g_l) = \theta_\Delta( 1, f_2, \dots, f_l) \cdot \theta_\Delta (\tilde k)\,.
\]
Then
\[
g^m \cdot \theta_{\Delta\,1} (g_1)^{-m} = \theta_\Delta (1, g_2, \dots, g_l)^m= \theta_\Delta (1, f_2, 
\dots, f_l)^m \cdot \theta_\Delta (\tilde k)^m 
= f^m \cdot \theta_{\Delta\,1} (f_1)^{-m} \cdot \theta_\Delta (\tilde k)^m\,,
\]
thus
\[
\theta_{\Delta\,1} (g_1)^m = \theta_{\Delta\,1} (f_1)^m \cdot \theta_\Delta (\tilde k)^{-m}\,.
\]
We know that $g_1^m \equiv f_1^m\ (\mod\, \tilde K_1)$, thus $\theta_\Delta (\tilde k)^{-m} \in 
\theta_\Delta (\tilde K_1)$. Moreover, $\theta_\Delta (\tilde K)$ is a free abelian group and 
$\theta_\Delta (\tilde K_1)$ is a direct factor of $\theta_\Delta (\tilde K)$, hence $\theta_\Delta 
(\tilde k) \in \theta_\Delta (\tilde K_1)$. So, we can assume $\tilde k \in \tilde K_1$. Set $f_1'=f_1 
\tilde k$ and $f_i'=g_i$ for $2 \le i\le l$. Then $f= \theta_\Delta (f_1', f_2', \dots, f_l')$, 
$\Delta(f_1')=\emptyset$, and $\theta_{\Delta\,1} (f_1')^m = \theta_{\Delta\,1} (g_1)^m$. Since 
$\theta_{\Delta\,1}: \MM (\Sigma_{\Delta\,1}) \to \MM (\Sigma)$ is injective, the later equality 
implies that $f_1'^m=g_1^m$. If $f_1'$ and $g_1$ are both pseudo-Anosov, then $f_1'=g_1$ by Theorem 4.5. 
If $f_1'$ and $g_1$ are both periodic, then they are conjugate by Lemma 7.2. So, $f$ and $g$ are 
conjugate.

\bigskip\noindent
{\bf Proof of Theorem 3.6, Case 2.} We assume that $d_0^{(1)}$ and $d_0^{(2)}$ are boundary components 
of $\Sigma_{\Delta\,1}$ (see Figure 7.2).
Set $\bar \Sigma_{\Delta\,1} = \Sigma_{\Delta\,1} /\sim$, where $\sim$ is the equivalence relation 
which identifies $d_0^{(1)}(z)$ with $d_0^{(2)}(\bar z) = (d_0^{(2)})^{-1}(z)$ 
for all $z \in \S^1$. Note that $\bar 
\Sigma_{\Delta\,1}$ is a surface of genus 1 which embeds in $\Sigma$, and this embedding determines a 
homomorphism $\bar \theta_{\Delta\,1}: \MM (\bar \Sigma_{\Delta\,1}) \to \MM(\Sigma)$ which 
is injective (see \cite{ParRol1}). 
For $2 \le i\le l$, the surface $\Sigma_{\Delta\,i}$ is of genus 0 and has at least 3 boundary 
components. It embeds in $\Sigma$, and this embedding determines an injective homomorphism 
$\theta_{\Delta\,i}: \MM (\Sigma_{\Delta\,i}) \to \MM (\Sigma)$. Set $\bar \Sigma_\Delta= \bar 
\Sigma_{\Delta\,1} \sqcup \Sigma_{\Delta\,2} \sqcup \dots \sqcup \Sigma_{\Delta\,l}$, and denote by 
$\pi: \bar \Sigma_\Delta \to \Sigma$ the natural quotient. We set $\MM (\bar \Sigma_\Delta)= \MM 
(\bar \Sigma_{\Delta\,1}) \times \MM (\Sigma_{\Delta\,2}) \times \dots \times \MM 
(\Sigma_{\Delta\,l})$, and we denote by $\bar \theta_\Delta: \MM (\bar \Sigma_\Delta) \to \MM (\Sigma)$ 
the homomorphism induced by $\pi$.

\bigskip\noindent
Here again, $f$ induces an isomorphism $f_\ast: \Gamma 
\to \Gamma$. Recall that $\Gamma$ has a unique reduced cycle which is of length 1. The vertex of this 
cycle is $v_1$ and the edge is $e_0$. In particular, we have $f_\ast (v_1) = v_1$ and $f_\ast (e_0) = 
e_0^{\pm 1}$. Let $\Gamma' = \Gamma \setminus \{e_0\}$. Then $f$ induces an isomorphism $f_\ast: 
\Gamma' \to \Gamma'$, the graph $\Gamma'$ is a tree, the leafs of $\Gamma'$ are $v_1', \dots, v_q'$, and 
possibly $v_1$, and $f_\ast(v_i') = v_i'$ for all $1 \le i\le q$, thus $f_\ast$ is the identity on 
$\Gamma'$.

\bigskip\noindent
The above shows that $f(\delta_k) = \delta_k$ for all $0 \le k\le p$, $f$ preserves the orientation of 
$d_k$ for all $1 \le k\le p$ (not necessarily for $k=0$), and $\Sigma_{\Delta\,j}$ is invariant under 
$f$ up to isotopy for all $1 \le j \le l$. It follows that there exists $(f_1, f_2, \dots, f_l) \in \MM 
(\bar \Sigma_\Delta)$ such that $f= \bar \theta_\Delta (f_1, f_2, \dots, f_l)$. Furthermore, we have $\Delta 
(f_1)= \{ \delta_0\}$, and $\Delta(f_i)= \emptyset$ for all $2 \le i\le l$, since $\Delta(f) = \Delta = 
\{ \delta_0, \delta_1, \dots, \delta_l\}$.
Similarly, there exists $(g_1,g_2, \dots, g_l) \in \MM (\bar \Sigma_\Delta)$ such that $\bar 
\theta_\Delta (g_1, g_2, \dots, g_l) = g$, $\Delta(g_1)= \{ \delta_0 \}$, and $\Delta(g_i)= \emptyset$ 
for all $2 \le i \le l$. Using the same arguments as in Case 1, it is easily shown that $(f_1, f_2, 
\dots, f_l)$ can be chosen so that $f_i=g_i$ for all $2 \le i\le l$, and $f_1^m=g_1^m$. By Lemma 7.3, 
it follows that $f_1$ and $g_1$ are conjugate in $\MM (\bar \Sigma_{\Delta\,1})$, thus $f$ and $g$ are 
conjugate.

\bigskip\noindent
{\bf Proof of Theorem 3.6, Case 3.} We assume that $d_0^{(1)}$ and $d_1^{(1)}$ are boundary components 
of $\Sigma_{\Delta\,1}$, and 
$d_0^{(2)}$ and $d_1^{(2)}$ are boundary components of $\Sigma_{\Delta\,2}$ (see Figure 7.3).
Set $\bar \Sigma_{\Delta\,1} = (\Sigma_{\Delta\,1} \sqcup \Sigma_{\Delta\,2}) /\sim$, where $\sim$ is 
the equivalence relation which identifies $d_i^{(1)}(z)$ with $d_i^{(2)}(\bar z)
=(d_i^{(2)})^{-1}(z)$ for all $i \in \{ 
0,1\}$ and all $z \in \S^1$. Note that $\bar \Sigma_{\Delta\,1}$ is a surface of genus 1 which embeds 
in $\Sigma$, and this embedding determines a homomorphism $\bar \theta_{\Delta\,1}: \MM (\bar 
\Sigma_{\Delta\,1}) \to \MM (\Sigma)$ which is injective. For $3 \le i \le l$, the surface 
$\Sigma_{\Delta\,i}$ is of genus 0 and has at least 3 boundary components. It embeds in $\Sigma$, and 
this embedding determines an injective homomorphism $\theta_{\Delta\,i}: \MM (\Sigma_{\Delta\,i}) \to 
\MM (\Sigma)$. Set $\bar \Sigma_\Delta = \bar \Sigma_{\Delta\,1} \sqcup \Sigma_{\Delta\,3} \sqcup \dots 
\sqcup \Sigma_{\Delta\,l}$, and denote by $\pi: \bar \Sigma_\Delta \to \Sigma$ he natural quotient. 
We set $\MM (\bar \Sigma_\Delta)= \MM (\bar \Sigma_{\Delta\,1}) \times \MM (\Sigma_{\Delta\,3}) 
\times \dots \times \MM (\Sigma_{\Delta\,l})$, and we denote by $\bar \theta_\Delta: \MM (\bar 
\Sigma_\Delta) \to \MM (\Sigma)$ the homomorphism induced by $\pi$.

\bigskip\noindent
Again, $f$ induces an isomorphism $f_\ast: \Gamma \to \Gamma$. Recall that $\Gamma$ has a	 unique 
reduced cycle which is of length 2. The vertices of this cycle are $v_1, v_2$ and the edges are $e_0, 
e_1$. The graph $\Gamma' = \Gamma \setminus \{ e_0, e_1\}$ is the disjoint union of two trees, 
$\Gamma_1', \Gamma_2'$, where $v_1$ (resp. $v_2$) is assumed to be a vertex of $\Gamma_1'$ (resp. 
$\Gamma_2'$). The vertex $v_1$ has valence $\ge 3$, thus $\Gamma_1'$ has more than one vertex and, 
equivalently, more than one leaf. We have $f_\ast (\Gamma')= \Gamma'$ so either $f_\ast (\Gamma_1') = 
\Gamma_1'$, or $f_\ast (\Gamma_1') = \Gamma_2'$. We actually have $f_\ast(\Gamma_1') = \Gamma_1'$ 
because $\Gamma_1'$ has at least one leaf of the form $v_j'$ and $f_\ast(v_j')= v_j'$. All the leafs 
of $\Gamma_1'$ but possibly one, $v_1$, are fixed by $f_\ast$, thus $f_\ast$ is the identity of 
$\Gamma_1'$. Similarly, $f_\ast$ is the identity on $\Gamma_2'$. In particular, we have $f_\ast (v_1) = 
v_1$ and $f_\ast (v_2) = v_2$.

\bigskip\noindent
The above shows that $f(\{ \delta_0, \delta_1 \}) = \{ \delta_0, \delta_1 \}$, $f (\delta_k) = \delta_k$ 
and $f$ preserves the orientation of $d_k$ for all $2 \le k \le p$, and $\Sigma_{\Delta\,j}$ is 
invariant under $f$ up to isotopy for all $1 \le j \le l$. It follows that there exists $(f_1, f_3, 
\dots, f_l) \in \MM (\bar \Sigma_\Delta)$ such that $\bar \theta_\Delta (f_1, f_3, \dots, f_l) = f$. 
Furthermore, we have $\Delta(f_1) = \{ \delta_0, \delta_1 \}$, and $\Delta (f_i)= \emptyset$ for all $3 
\le i \le l$, because $\Delta (f) = \Delta = \{ \delta_0, \delta_1, \dots, \delta_p \}$. 
Similarly, there exists $(g_1, g_3, \dots, g_l) \in \MM (\bar \Sigma_\Delta)$ such that $\bar 
\theta_\Delta (g_1, g_3, \dots, g_l)= g$, $\Delta(g_1)= \{ \delta_0, \delta_1 \}$, and $\Delta(g_i)= 
\emptyset$ for all $3 \le i\le l$. Using again the same arguments as in Case 1, it is easily shown that 
$(f_1, f_3, \dots, f_l)$ can be chosen so that $f_i=g_i$ for all $3 \le i \le l$, and $f_1^m=g_1^m$. By 
Lemma 7.4, it follows that $f_1$ and $g_1$ are conjugate in $\MM (\bar \Sigma_{\Delta\,1})$, thus $f$ 
and $g$ are conjugate.

\bigskip\noindent
{\bf Proof of Theorem 3.6, Case 4.} (See Figure 7.4.)
For $1 \le i \le l$, the surface $\Sigma_{\Delta\,i}$ is of genus $0$ and has at least 3 boundary 
components. It embeds in $\Sigma$, and this embedding determines an injective homomorphism 
$\theta_{\Delta\,i}: \MM (\Sigma_{\Delta\,i}) \to \MM (\Sigma)$ (see \cite{ParRol1}). Let $\pi: 
\Sigma_\Delta \to \Sigma$ denote the natural quotient. We set $\MM (\Sigma_\Delta) = 
\MM(\Sigma_{\Delta\,1}) \times \MM(\Sigma_{\Delta\,2}) \times \dots \times \MM(\Sigma_{\Delta\,l})$, 
and we denote by $\theta_\Delta: \MM (\Sigma_\Delta) \to \MM (\Sigma)$ the homomorphism induced by 
$\pi$.

\bigskip\noindent
Like in the previous cases, $f$ induces an isomorphism $f_\ast: \Gamma \to \Gamma$. Recall that $\Gamma$ 
has a unique reduced cycle which is of length $\ge 3$. We assume that the vertices of this cycle are 
$v_1, \dots, v_r$ and the edges are $e_0, e_1, \dots, e_{r-1}$. The graph $\Gamma'= \Gamma \setminus \{ 
e_0, e_1, \dots, e_{r-1} \}$ is the disjoint union of $r$ trees, $\Gamma_1', \dots, \Gamma_r'$, where 
$v_i$ is assumed to be a vertex of $\Gamma_i'$. Since $v_i$ has valence $\ge 3$, the tree $\Gamma_i'$ 
has more than one vertex and, equivalently, more than one leaf. We have $f_\ast(\Gamma')= \Gamma'$ and 
$\Gamma_i'$ has a leaf of the form $v_j'$, thus $f_\ast(\Gamma_i') = \Gamma_i'$. Furthermore, all the 
leafs but possibly one, $v_i$, are fixed under $f_\ast$, thus $f_\ast$ is the identity on $\Gamma_i'$. 
This also implies that $f_\ast (v_i)=v_i$. Finally, we have $f_\ast (e_i) = e_i$ for all $0 \le i \le 
r-1$ because the reduced cycle of $\Gamma$ is invariant under $f_\ast$, all its vertices are fixed by 
$f_\ast$, and its length is $\ge 3$. We conclude that $f_\ast$ is the identity on the whole graph 
$\Gamma$.

\bigskip\noindent
The above shows that $f(\delta_k) = \delta_k$ and $f$ preserves the orientation of $d_k$ for all $0 \le 
k\le p$, and $\Sigma_{\Delta\,j}$ is invariant under $f$ up to isotopy for all $1 \le j\le l$. It 
follows that there exists $(f_1, f_2, \dots, f_l) \in \MM (\Sigma_\Delta)$ such that $f = 
\theta_\Delta (f_1, f_2, \dots, f_l)$. Furthermore, we have $\Delta(f_i)= \emptyset$ for all $1 \le 
i\le l$ because $\Delta (f)= \Delta = \{ \delta_0, \delta_1, \dots, \delta_p\}$.
Similarly, there exists $(g_1, g_2, \dots, g_l) \in \MM (\Sigma_\Delta)$ such that 
$\theta_\Delta (g_1, g_2, \dots, g_l) =g$ and $\Delta(g_i)= \emptyset$ for all $1 \le i\le l$. Using 
the same arguments as in Case 1, it is easily shown that $(f_1, f_2, \dots, f_l)$ can be chosen so that 
$f_i=g_i$ for all $2 \le i\le l$, and $f_1^m=g_1^m$. We conclude that $f_1=g_1$ by Proposition 3.2, hence $f=g$.
\qed


\subsection{Proof of Lemma 7.1.}

{\bf Lemma 7.5.} {\it Let $H: \S^2 \to \S^2$ be a diffeomorphism of the sphere $\S^2$ of finite order $m 
\ge 2$. Then $H$ has at most 2 fixed points.}

\bigskip\noindent
{\bf Proof.}
Consider the ramified covering $\pi: \S^2 \to \S^2/H$ and denote by $Q_1, \dots, Q_l$ the ramification 
points of $\pi$. By Lemma 2.27, we have
\[
\chi (\S^2) + \sum_{i=1}^l (m-o(Q_i)) = m \cdot \chi (\S^2/H)\,.
\]
Suppose that $H$ has at least 3 fixed points. Then $l \ge 3$ and we can suppose $o(Q_1)= o(Q_2) = 
o(Q_3)=1$. From the above equality follows 
\[
3m-1 = \chi(\S^2) + 3(m-1) \le m \cdot \chi(\S^2/H) \le 2m\,,
\]
and this inequality holds only if $m \le 1$: a contradiction.
\qed

\bigskip\noindent
{\bf Proof of Lemma 7.1.}  Let $\Sigma$ be a surface of genus $0$ with $q$ boundary components, $c_1, \dots, 
c_q$, where $q \ge 3$. Let $K$ be the subgroup of $\MM (\Sigma)$ generated by $\{ \tau_{c_1}, \dots, 
\tau_{c_q}\}$. Let $f,g \in \MM(\Sigma)$ be two pseudo-Anosov elements such that $f^m \equiv g^m\ (\mod\, K)$ 
for some $m \ge 1$.

\bigskip\noindent
Consider the corking $\theta: \MM(\Sigma) \to \MM (\Sigma_0, \QQ)$ of $\Sigma$ as defined in Subsection 
2.1. By Proposition 2.1, the kernel of $\theta$ is precisely $K$, thus the congruence $f^m \equiv g^m\ 
(\mod\,K)$ is equivalent to the equality $\theta(f)^m = \theta(g)^m$. Take a pseudo-Anosov 
representative $F$ of $\theta(f)$ and denote by $\FF^s$ (resp. $\FF^u$) the stable (resp. unstable) 
foliation of $F$. By Corollary 2.14, there exists a pseudo-Anosov representative $G \in \Diff 
(\Sigma_0, \QQ)$ of $\theta(g)$ such that $F^m=G^m$. Set $H=F^{-1}G$. Then $H \in \Sym(\FF^s, 
\FF^u)$, a finite group (see Corollary 2.10). Observe furthermore that $H(Q)=Q$ for all $Q \in \QQ$ and 
$|\QQ| \ge 3$, thus, by Lemma 7.5, $H=\Id$ and $F=G$. We conclude that $\theta(f)=\theta(g)$, that is, $f 
\equiv g\ (\mod\,K)$.
\qed


\subsection{Proof of Lemma 7.2}

The following lemmas 7.6--7.10 are preliminaries to the proof of Lemma 7.2.

\bigskip\noindent
Let $\Sigma$ be a surface, and let $\PP$ be a 
finite set of punctures in the interior of $\Sigma$. Then $\PDiff (\Sigma, \PP)$ denotes the group of 
$F \in \Diff (\Sigma,\PP)$ which pointwise fixed the elements of $\PP$, and $\PP\MM (\Sigma, \PP)$ 
denotes the group of isotopy classes of elements of $\PDiff (\Sigma, \PP)$. Note that $\PP \MM (\Sigma, 
\PP)$ is a normal subgroup of $\MM (\Sigma, \PP)$ and we have the exact sequence
\[
1 \to \PP \MM (\Sigma, \PP) \to \MM (\Sigma, \PP) \to \Sym (\PP) \to 1\,.
\]

\bigskip\noindent
{\bf Lemma 7.6.} {\it Let $\PP= \{ P_1, P_2, P_3, P_4 \}$ be a set of 4 punctures in the torus $\T^2$. 
Let $f,g \in \PP \MM (\T^2, \PP)$ be two elements of order 2. Then $f$ and $g$ are conjugate in $\PP 
\MM (\T^2, \PP)$.}

\bigskip\noindent
{\bf Proof.}
Take a diffeomorphism $F \in \PDiff (\T^2, \PP)$ of order 2 which represents $f$, and consider the 
ramified covering $\pi: \T^2 \to \T^2/F$. Let $Q_1, \dots, Q_l$ be the ramification points of $\pi$. 
Since $F$ is of order 2, we have $o(Q_i)=|\pi^{-1} (Q_i)| =1$ for all $1 \le i\le l$. Moreover, we can and do 
assume that $\pi^{-1} (Q_i)= \{ P_i\}$ for $1 \le i\le 4$. In particular, $l \ge 4$. By Lemma 2.27, we have
\[
\chi(\T^2) +l = 2 \cdot \chi( \T^2/ F)\,.
\]
thus $\chi (\T^2/F)=2$ (that is, $\T^2/F$ is a sphere), and $l=4$.

\bigskip\noindent
We choose a cellular decomposition of $\T^2/F$ whose vertices are $Q_1, Q_2, Q_3, Q_4$, and having 4 
arrows, $a_1, a_2, a_3, a_4$, and two faces, $A_1, A_2$ (see Figure 7.7). The source of $a_i$ is assumed 
to be $Q_i$ and its target to be $Q_{i+1}$ (the indices are taken mod 4). Moreover, we assume $\partial 
A_1= a_1a_2a_3a_4$ and $\partial A_2= a_4^{-1} a_3^{-1} a_2^{-1} a_1^{-1}$. For $1 \le i\le 4$, the set $\pi^{-
1}(a_i)$ is the union of two edges, $b_i$ and $b_i'$, such that $b_i \cap b_i' = \{ P_i, P_{i+1}\}$. 
For $j=1,2$, the set $\pi^{-1}(A_j)$ is the union of two faces, $B_j$ and $B_j'$, such that $B_j \cap 
B_j'= \{ P_1, P_2, P_3, P_4\}$. Furthermore, the set $\{P_1, P_2, P_3, P_4\} \cup \{b_1,b_1', \dots, b_4,b_4'\} 
\cup \{B_1,B_1',B_2,B_2'\}$ determines a cellular decomposition of $\T^2$. Up to a permutation of 
either some $\{b_i, b_i'\}$, or some $\{B_j, B_j'\}$, such a cellular decomposition is unique. So, $F$ is unique 
up to a conjugation in $\PDiff (\T^2, \PP)$.
\qed

\begin{figure}[htbp]
\centerline{
\setlength{\unitlength}{.4cm}
\begin{picture}(25,12)
\put(2,2){\includegraphics[width=8.8 cm]{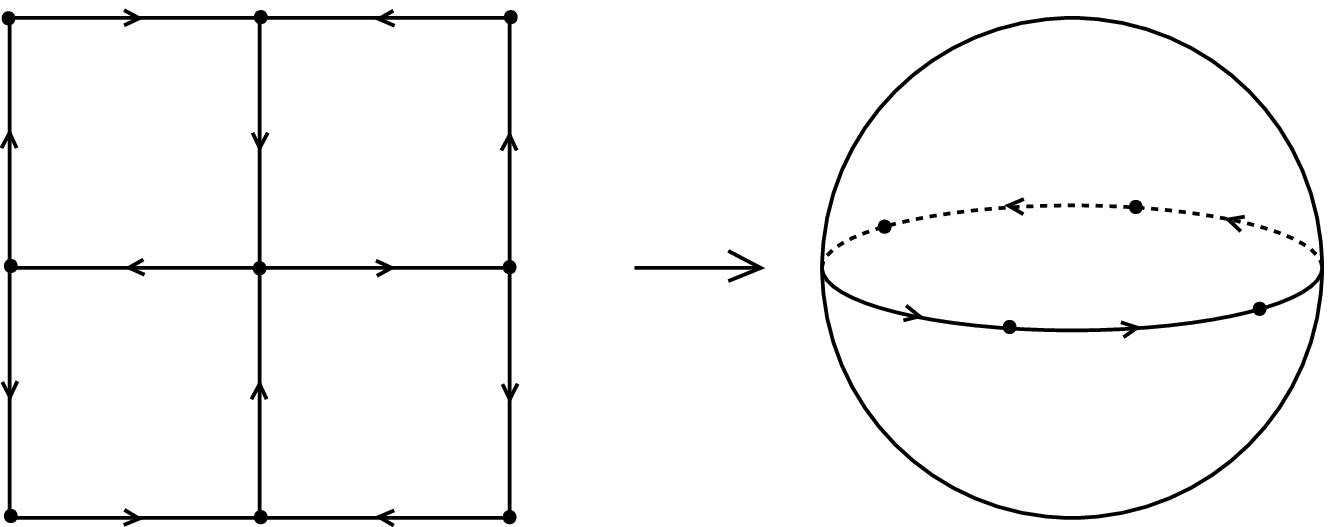}}
\put(1.2,1.2){\small $P_3$}
\put(5.9,1.2){\small $P_4$}
\put(10.4,1.2){\small $P_3$}
\put(1,6){\small $P_2$}
\put(5.2,5.4){\small $P_1$}
\put(10.6,6.1){\small $P_2$}
\put(1,10.5){\small $P_3$}
\put(5.9,10.8){\small $P_4$}
\put(10.4,10.8){\small $P_3$}
\put(3.7,1.2){\small $b_3'$}
\put(7.9,1.2){\small $b_3$}
\put(1,3.9){\small $b_2'$}
\put(5.2,3.9){\small $b_4'$}
\put(10.6,3.9){\small $b_2'$}
\put(3.9,6.7){\small $b_1'$}
\put(7.7,6.7){\small $b_1$}
\put(1,7.9){\small $b_2$}
\put(5.3,7.9){\small $b_4$}
\put(10.6,7.9){\small $b_2$}
\put(3.7,10.8){\small $b_3'$}
\put(7.9,10.8){\small $b_3$}
\put(3.7,4){\small $B_1'$}
\put(7.7,4){\small $B_2$}
\put(3.7,8){\small $B_2'$}
\put(7.7,8){\small $B_1$}
\put(13.1,6.6){$\pi$}
\put(18.5,4.4){\small $Q_1$}
\put(22.3,4.6){\small $Q_2$}
\put(16.2,7.3){\small $Q_4$}
\put(20.5,7.6){\small $Q_3$}
\put(16.5,4.8){\small $a_4$}
\put(20.5,4.5){\small $a_1$}
\put(18.5,7.6){\small $a_3$}
\put(22.3,7.4){\small $a_2$}
\put(19.5,2.5){\small $A_2$}
\put(19.5,9.3){\small $A_1$}
\end{picture}}
\centerline{{\bf Figure 7.7.} Ramified covering associated to an order 2 diffeomorphism of $\T^2$.}
\end{figure}

\bigskip\noindent
{\bf Lemma 7.7.} {\it Let $\PP= \{P_1, P_2, P_3\}$ be a set of 3 punctures in the torus $\T^2$. Let 
$f,g \in \PP \MM (\T^2, \PP)$ be two elements of order 2. Then $f$ and $g$ are conjugate.}

\bigskip\noindent
{\bf Proof.}
Choose an extra point $P_4 \in \T^2 \setminus \PP$. Let $F \in \PDiff (\T^2, \PP)$ be a diffeomorphism 
of order 2 which represents $f$, and consider the ramified covering $\pi: \T^2 \to \T^2/F$. Let $Q_1, 
\dots, Q_l$ be the ramification points of $\pi$. Here again, we have $o(Q_i)=| \pi^{-1} (Q_i)| = 1$ for all $1 
\le i \le l$, and we can assume that $\pi^{-1} (Q_i)= \{ P_i\}$ for $1 \le i \le 3$. In particular, $l 
\ge 3$. By Lemma 2.27, we have
\[
\chi (\T^2) +l= 2 \cdot \chi( \T^2/F)\,,
\]
thus $\chi (\T^2/F)=2$ (that is, $\T^2/F$ is a sphere), and $l=4$. Obviously, we can assume that $\pi^{-1} 
(Q_4)= \{P_4\}$.

\bigskip\noindent
Let $\varphi: \PP\MM (\T^2, \PP \cup \{ P_4 \}) \to \PP\MM (\T^2, \PP)$ be the natural epimorphism. The 
above shows that there exists $\tilde f \in \PP\MM (\T^2, \PP \sqcup \{P_4\})$ of order 2 such that $f= 
\varphi(\tilde f)$. Similarly, there exists $\tilde g \in \PP\MM( \T^2, \PP \sqcup \{P_4\})$ of order 2 
such that $g= \varphi (\tilde g)$. By Lemma 7.6, $\tilde f$ and $\tilde g$ are conjugate, thus $f$ and $g$ 
are conjugate, too.
\qed

\bigskip\noindent
{\bf Lemma 7.8.} {\it Let $\PP= \{ P_1, P_2, P_3\}$ be a set of 3 punctures in the torus $\T^2$. Let $f,g 
\in \PP\MM (\T^2, \PP)$ be two elements of order 3. Then $f$ and $g$ are conjugate.}

\bigskip\noindent
{\bf Proof.}
The proof is similar to the one of Lemma 7.6. We choose a diffeomorphism $F \in \PDiff 
(\T^2, \PP)$ of order 3 which represents $f$, and we consider the ramified covering $\pi: \T^2 \to 
\T^2/F$. Let $Q_1, \dots, Q_l$ be the ramification points of $\pi$. We have $o(Q_i)=|\pi^{-1}(Q_i)|=1$ for all $1 
\le i \le l$, and we can assume that $\pi^{-1} (Q_i)= \{P_i\}$ for $1 \le i\le 3$. 
In particular, $l \ge 3$.
By Lemma 2.27, we have
\[
\chi( \T^2) +2l = 3 \cdot \chi (\T^2/F)\,,
\]
thus $\chi(\T^2/F)=2$ (that is, $\T^2/F$ is a sphere), and $l=3$.

\bigskip\noindent
We choose a cellular decomposition of $\T^2/F$, whose vertices are $Q_1,Q_2,Q_3$, and having 3 arrows, 
$a_1,a_2,a_3$, and 2 faces, $A_1,A_2$ (see Figure 7.8). The source of $a_i$ is assumed to be $Q_i$ and 
its target to be $Q_{i+1}$. Moreover, we assume that $\partial A_1= a_1a_2a_3$ and $\partial A_2= a_3^{-1} 
a_2^{-1} a_1^{-1}$. For $1 \le i\le 3$, the set $\pi^{-1}(a_i)$ is the union of three edges, $b_i, b_i', 
b_i''$, such that $b_i \cap b_i' = b_i \cap b_i'' = b_i' \cap b_i'' = \{ P_i, P_{i+1}\}$. For $j=1,2$, 
the set $\pi^{-1}(A_j)$ is the union of three faces, $B_j, B_j', B_j''$, such that $B_j \cap B_j'= B_j \cap 
B_j''= B_j' \cap B_j'' = \{ P_1, P_2, P_3\}$. Furthermore, $\{ P_1, P_2, P_3 \} \cup \{b_1,b_1',b_1'', 
b_2,b_2',b_2'', b_3,b_3',b_3''\} \cup \{ B_1,B_1',B_1'', B_2,B_2',B_2''\}$ determines a cellular 
decomposition of $\T^2$. Up to a permutation of either some $\{ b_i,b_i',b_i'' \}$, or some 
$\{B_j,B_j',B_j''\}$, such a decomposition is unique. So, $F$ is unique up to a conjugation in $\PDiff 
(\T^2, \PP)$.
\qed

\begin{figure}[htbp]
\centerline{
\setlength{\unitlength}{.4cm}
\begin{picture}(29,14)
\put(2,1.5){\includegraphics[width=10.4 cm]{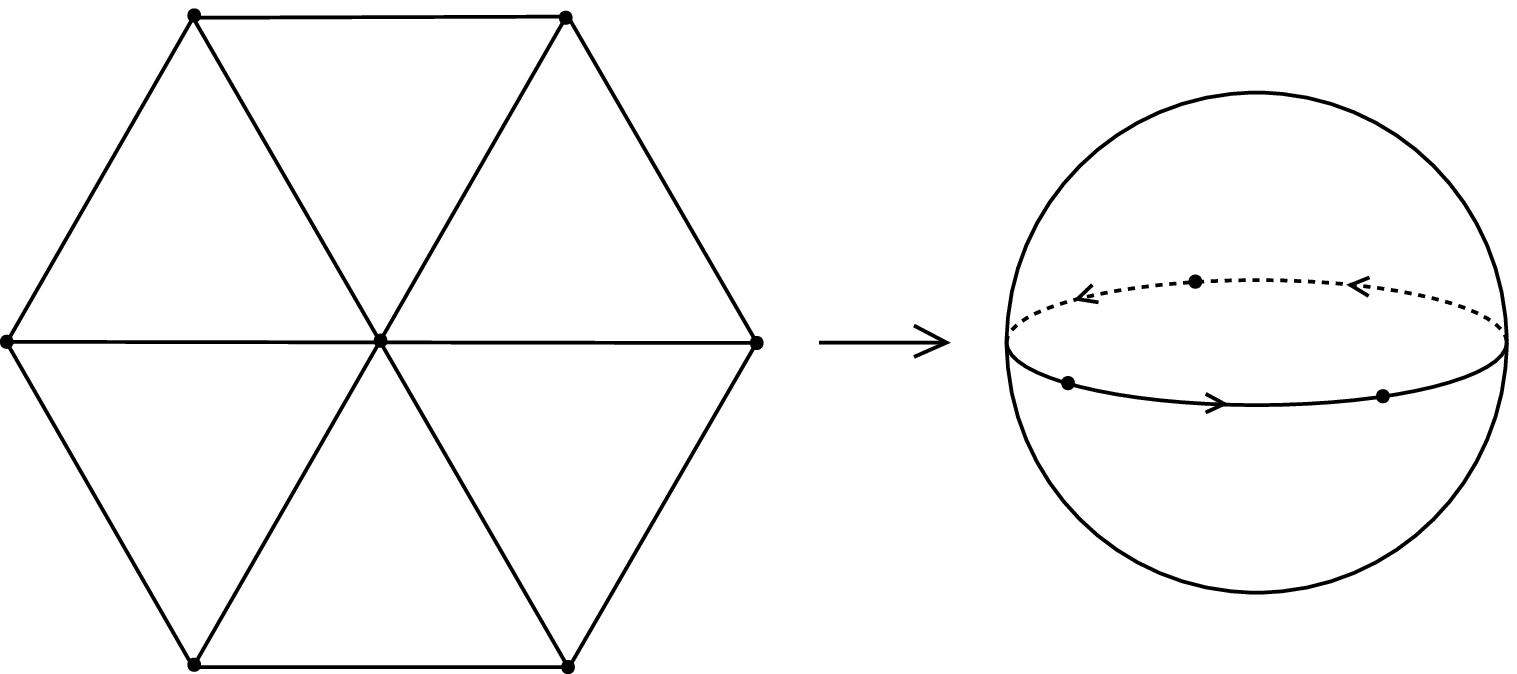}}
\put(4.4,0.8){\small $P_3$}
\put(11.8,0.8){\small $P_2$}
\put(0.9,7){\small $P_2$}
\put(8,8.4){\small $P_1$}
\put(15,6.9){\small $P_3$}
\put(4.3,13){\small $P_3$}
\put(11.5,13){\small $P_2$}
\put(7.7,0.9){\small $b_2''$}
\put(2.8,3.5){\small $b_2'$}
\put(6.9,3.7){\small $b_3''$}
\put(10.1,4.5){\small $b_1$}
\put(13.3,3.5){\small $b_2$}
\put(4.7,6.2){\small $b_1''$}
\put(11.1,7.4){\small $b_3$}
\put(2.6,9.7){\small $b_2$}
\put(6.1,9.1){\small $b_3'$}
\put(8.9,9.8){\small $b_1'$}
\put(13.5,9.8){\small $b_2'$}
\put(7.9,13.1){\small $b_2''$}
\put(4.5,4.5){\small $B_1''$}
\put(7.9,3){\small $B_2''$}
\put(11.3,4.9){\small $B_1$}
\put(4.5,8.7){\small $B_2'$}
\put(7.8,11){\small $B_1'$}
\put(11.5,8.9){\small $B_2$}
\put(16.7,7.4){$\pi$}
\put(20.1,5.6){\small $Q_1$}
\put(25.5,5.2){\small $Q_2$}
\put(22.1,8.6){\small $Q_3$}
\put(23,5.3){\small $a_1$}
\put(20.3,8.4){\small $a_3$}
\put(25,8.5){\small $a_2$}
\put(23.3,3.5){\small $A_2$}
\put(23.3,10.2){\small $A_1$}
\end{picture}}
\centerline{{\bf Figure 7.8.} Ramified covering associated to an order 3 diffeomorphism of $\T^2$.}
\end{figure}

\bigskip\noindent
{\bf Lemma 7.9.} {\it Let $\Sigma$ be a surface of genus 1 with 3 boundary components, $c_1,c_2,c_3$. 
Then there exist $f_0,g_0, \in \MM(\Sigma)$ such that $g_0^3= f_0^2 = \tau_{c_1} \tau_{c_2} 
\tau_{c_3}$.}

\bigskip\noindent
{\bf Proof.}
Consider the curves $a_1,a_2,a_3,b$ pictured in Figure 7.9, and set $f_0= \tau_{a_1} \tau_b \tau_{a_2} 
\tau_b \tau_{a_3} \tau_b$ and $g_0= \tau_{a_1} \tau_{a_2} \tau_{a_3} \tau_b$. Then $g_0^3= f_0^2 = 
\tau_{c_1} \tau_{c_2} \tau_{c_3}$. These equalities can be easily checked inspecting the action of $g_0^3$ 
(resp. $f_0^2$) on the isotopy classes of the curves $a_1, a_2, a_3, b$, and on the isotopy classes of 
some chosen arcs that join different components of $\partial \Sigma$. This can be found, for instance, 
in \cite{LabPar1}.
\qed

\begin{figure}[htbp]
\centerline{
\setlength{\unitlength}{.4cm}
\begin{picture}(11,8)
\put(1,1){\includegraphics[width=3.6 cm]{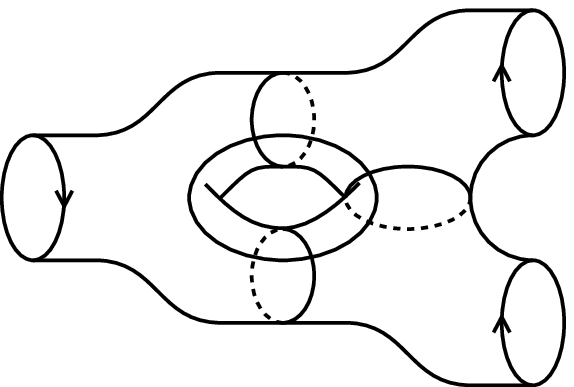}}
\put(5.2,1.3){\small $a_1$}
\put(8.7,3.8){\small $a_2$}
\put(5.2,6.3){\small $a_3$}
\put(2.2,3.5){\small $c_3$}
\put(8,1.5){\small $c_1$}
\put(8,5.5){\small $c_2$}
\put(3.5,3.5){\small $b$}
\end{picture}}
\centerline{{\bf Figure 7.9.} A surface of genus 1 with 3 boundary components.}
\end{figure}

\bigskip\noindent
{\bf Lemma 7.10.} {\it Let $\Sigma$ be a surface of genus 1 with 4 boundary components, $c_1,c_2,c_3,c_4$. 
Then there exists $f_0 \in \MM(\Sigma)$ such that $f_0^2= \tau_{c_1} \tau_{c_2} \tau_{c_3} 
\tau_{c_4}$.}

\bigskip\noindent
{\bf Proof.}
Consider the curves $a_1, a_2, a_3, a_4, b$ pictured in Figure 7.10, and set $f_0= \tau_{a_1} \tau_{a_3} 
\tau_b \tau_{a_2} \tau_{a_4} \tau_b$. Then $f_0^2= \tau_{c_1} \tau_{c_2} \tau_{c_3} \tau_{c_4}$. Here 
again, this equality can be easily checked inspecting the action of $f_0^2$ on the isotopy classes of 
the curves $a_1, a_2, a_3, a_4, b$, and on the isotopy classes of some chosen arcs that join different 
components of $\partial \Sigma$.
\qed

\begin{figure}[htbp]
\centerline{
\setlength{\unitlength}{.4cm}
\begin{picture}(11,8)
\put(1,1){\includegraphics[width=3.6 cm]{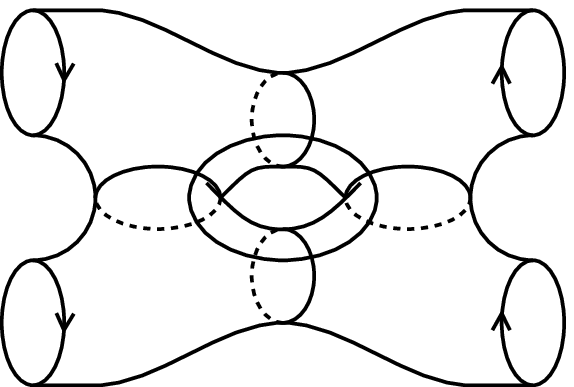}}
\put(5.2,1.2){\small $a_1$}
\put(8.7,3.7){\small $a_2$}
\put(5.2,6.3){\small $a_3$}
\put(1.5,4){\small $a_4$}
\put(8,1.5){\small $c_1$}
\put(8,5.5){\small $c_2$}
\put(2.2,5.5){\small $c_3$}
\put(2.2,1.5){\small $c_4$}
\put(6.5,2.5){\small $b$}
\end{picture}}
\centerline{{\bf Figure 7.10.} A surface of genus 1 with 4 boundary components.}
\end{figure}

\bigskip\noindent
{\bf Proof of Lemma 7.2.} Let $\Sigma$ be a surface of genus 1 with 
$q$ boundary components, $c_1, \dots, c_q$, where $q \ge 1$, and
let $f,g \in \MM( 
\Sigma)$ be two periodic elements such that $f^m = g^m$ for some $m \ge 1$.

\bigskip\noindent
If $q=1$, then $\MM (\Sigma) \simeq \BB_3$, and if $q=2$, then $\MM (\Sigma)= \BB_4 \times \Z$ (see 
\cite{LabPar1}, \cite{Gerva1}). In both cases the conclusion of Lemma 7.2 follows from \cite{Gonza1}. 
So, we can assume $q \ge 3$.

\bigskip\noindent
Let $\theta: \MM (\Sigma) \to \MM (\Sigma_0, \QQ)$ be the corking of $\Sigma$, as defined in Subsection 
2.1. Note that $\Sigma_0$ is a torus and $|\QQ|=q\ge 3$. Since $f$ and $g$ are periodic, both, 
$\theta(f)$ and $\theta(g)$, are finite order elements. Moreover, they both belong to $\PP\MM (\Sigma_0, 
\QQ)$. We denote by $o(f)$ (resp. $o(g)$) the order of $\theta(f)$ (resp. $\theta(g)$). Let $F \in 
\PDiff (\Sigma_0, \QQ)$ be a diffeomorphism of order $o(f)$ which represents $\theta(f)$, and consider 
the ramified covering $\pi: \Sigma_0 \to \Sigma_0/F$. The fact that $\QQ$ is pointwise fixed under $F$ 
implies, by Lemma 2.27, that
\[
\chi(\Sigma_0) + q(o(f)-1) \le o(f) \cdot \chi(\Sigma_0/F)\,.
\]
Obviously, this inequality holds only if either $o(f)=1$, or $o(f)=2$ and $q \in \{ 3,4\}$, or 
$o(f)=3$ and $q=3$. Similarly, we have either $o(g)=1$, or $o(g)=2$ and $q \in \{3,4\}$, or $o(g)=3$ 
and $q=3$.

\bigskip\noindent
Set $O(q)= \{1,2,3\}$ if $q=3$, $O(q)= \{1,2\}$ if $q=4$, and $O(q)= \{1\}$ if $q \ge 5$. By Lemmas 7.9 and 7.10, 
for all $r \in O(q)$, there exists $u_r \in \MM (\Sigma)$ such that $(u_r)^r = \tau_{c_1} 
\tau_{c_2} \dots \tau_{c_q}$. By Lemmas 7.6, 7.7 and 7.8, there exist $x_1, \dots, x_q \in \Z$ such that $f$ 
is conjugate to $u_{o(f)} \cdot \tau_{c_1}^{x_1} \dots \tau_{c_q}^{x_q}$. Thus $f^{6m}$ is conjugate to 
$\tau_{c_1}^{6mx_1+{6m \over o(f)}} \dots \tau_{c_q}^{6mx_q+ {6m\over o(f)}}$. This implies that
\[
f^{6m} = \tau_{c_1}^{6mx_1+{6m \over o(f)}} \dots \tau_{c_q}^{6mx_q+ {6m\over o(f)}}\,,
\]
since $\tau_{c_1}^{6mx_1+{6m \over o(f)}} \dots \tau_{c_q}^{6mx_q+ {6m\over o(f)}}$ is central in $\MM 
(\Sigma)$. Similarly, there exist $y_1, \dots, y_q \in \Z$ such that $g$ is conjugate to $u_{o(g)} 
\cdot \tau_{c_1}^{y_1} \dots \tau_{c_q}^{y_q}$, and
\[
g^{6m} = \tau_{c_1}^{6my_1+{6m \over o(g)}} \dots \tau_{c_q}^{6my_q+{6m \over o(g)}}\,.
\]
Now, the equality $f^{6m} = g^{6m}$ implies that $o(f)=o(g)$ and $x_i=y_i$ for all $1 \le i\le q$, thus $f$ 
and $g$ are conjugate.
\qed


\subsection{Proof of Lemma 7.3}

The following lemmas 7.11--7.13 are preliminaries to the proof of Lemma 7.3. Lemma 7.11 can be easily 
proved using similar arguments as in the proofs of Lemmas 7.5 and 7.8, so its proof is left to the 
reader. Lemma 7.12 is quite standard, so its proof is also left to the reader.

\bigskip\noindent
{\bf Lemma 7.11.} {\it Let $\PP= \{ P_1, P_2, P_3, P_4\}$ be a set of 4 punctures in the sphere $\S^2$. 
Let $f,g \in \MM (\S^2, \PP)$ be two elements of order 2 such that $f(P_1)= g(P_1)=P_1$, $f(P_2)= 
g(P_2)= P_2$, $f(P_3)= g(P_3)=P_4$, and $f(P_4)=g(P_4)=P_3$. Then there exists $h \in \PP \MM (\S^2, 
\PP)$ such that $g=hfh^{-1}$.}
\qed

\bigskip\noindent
{\bf Lemma 7.12.} {\it Let $\Sigma$ be a surface of genus 0 with 2 boundary components, $c_1, c_2$, and 
let $\PP= \{ P_1, P_2 \}$ be a set of 2 punctures in the interior of $\Sigma$. Then there exists $f_0 
\in \MM(\Sigma, \PP)$ such that $f_0(P_1)=P_2$, $f_0(P_2)=P_1$, and $f_0^2=\tau_{c_1} \tau_{c_2}$.}
\qed

\bigskip\noindent
{\bf Lemma 7.13.} {\it Let $\Sigma$ be a surface of genus 0 with $q$ boundary components, $c_1, \dots, 
c_q$, where $q \ge 1$, and let $\PP= \{P_1, P_2 \}$ be a set of two punctures in the interior of 
$\Sigma$. Let $f,g \in \MM(\Sigma, \PP)$ be two periodic elements. If $f^m=g^m$ for some $m \ge 1$, then 
$f$ and $g$ are conjugate.}

\bigskip\noindent
{\bf Proof.}
If $q=1$, then $\MM(\Sigma, \PP) \simeq \Z$. The conclusion of Lemma 7.13 obviously holds in this case, so 
we can assume $q \ge 2$.

\bigskip\noindent
Consider the corking $\theta: \MM(\Sigma, \PP) \to \MM (\Sigma_0, \PP \sqcup \QQ)$ of $(\Sigma, \PP)$. 
Then $\Sigma_0$ is a sphere and $|\QQ|=q$. Since $f$ and $g$ are periodic, $\theta(f)$ and $\theta(g)$ 
have finite order. Let $o(f)$ (resp. $o(g)$) denote the order of $\theta(f)$ (resp. $\theta(g)$). Choose 
a diffeomorphism $F \in \Diff (\Sigma_0, \PP \sqcup \QQ)$ of order $o(f)$ which represents $\theta(f)$, 
and consider the ramified covering $\pi: \Sigma_0 \to \Sigma_0/F$. Note that $\QQ$ is pointwise fixed 
under $F$, thus, by Lemma 2.27,
\[
\chi(\Sigma_0) + q(o(f)-1) \le o(f) \cdot \chi( \Sigma_0/F)\,.
\]
This equality holds only if either $o(f)=1$, or $q=2$. If $q=2$, then $o(f)=2$ because, in that 
case, $|\PP|=2$ and $\PP$ must be an orbit of $F$. Similarly, we have either $o(g)=1$, or $o(g)=2$ and 
$q=2$.

\bigskip\noindent
Set $O(q)= \{1,2\}$ if $q=2$, and $O(q)=\{1\}$ if $q \ge 3$. By Lemma 7.12, for all $r \in O(q)$, there 
exists $u_r \in \MM (\Sigma, \PP)$ such that $(u_r)^r= \tau_{c_1} \dots \tau_{c_q}$. By Lemma 7.11, there 
exist $x_1, \dots, x_q \in \Z$ such that $f$ is conjugate to $u_{o(f)} \cdot \tau_{c_1}^{x_1} \dots 
\tau_{c_q}^{x_q}$
(recall that, if $o(f)=2$, then $q=2$, $\PP \sqcup \QQ$ has 4 elements, and $\theta (f)$ permutes
the elements of $\PP$ and fixes the elements of $\QQ$),
hence $f^{2m}= \tau_{c_1}^{2mx_1+{2m\over o(f)}} \dots 
\tau_{c_q}^{2mx_q+{2m\over o(f)}}$. Similarly, there exist $y_1, \dots, y_q \in \Z$ such that $g$ is 
conjugate to $u_{o(g)} \cdot \tau_{c_1}^{y_1} \dots \tau_{c_q}^{y_q}$, and $g^{2m} = 
\tau_{c_1}^{2my_1 + {2m \over o(g)}} \dots \tau_{c_q}^{2my_q + {2m \over o(g)}}$. Finally, the equality 
$f^{2m} = g^{2m}$ implies that $o(f)=o(g)$ and $x_i=y_i$ for all $1 \le i \le q$, thus $f$ and $g$ are 
conjugate.
\qed

\bigskip\noindent
{\bf Proof of Lemma 7.3.} Let $\Sigma$ be a surface of genus 1 with $q$ boundary components, 
$c_1, \dots, c_q$, where $q \ge 1$,
and let $d: \S^1 \to \Sigma$ be a non-separating closed curve (see Figure 7.5). 
We denote by $\Sigma_\Delta$ the surface obtained from $\Sigma$ cutting along $d$. So, $\Sigma_\Delta$ 
is a surface of genus 0 with $q+2$ boundary components, $c_1, \dots, c_q, d^{(1)}, d^{(2)}$, and 
$\Sigma= \Sigma_\Delta / \sim$, where $\sim$ is the equivalence relation which identifies $d^{(1)}(z)$ 
with $d^{(2)}(\bar z) = (d^{(2)})^{-1}(z)$ for all $z \in \S^1$. Take two copies 
$\D_1$ and $\D_2$ of the standard disk 
$\D$, and denote by $e_i: \S^1 \to \partial \D_i$ the boundary component of $\D_i$ ($i=1,2$). Define 
the surface $\Sigma_0= (\Sigma_\Delta \sqcup \D_1 \sqcup \D_2)/\sim$, where $\sim$ is the equivalence 
relation which identifies $d^{(i)}(z)$ with $e_i(\bar z) = e_i^{-1}(z)$ for all $z \in \S^1$ and all $i \in \{ 1,2 
\}$. Note that $\Sigma_0$ is a surface of genus 0 with $q$ boundary components, $c_1, \dots, c_q$. We 
fix a point $Q_i$ in the interior of $\D_i$ for all $i \in \{1,2\}$.

\bigskip\noindent
We view $\delta= \langle d \rangle$ as a 0-simplex of $\CC (\Sigma)$, and we set $\Stab (\delta)= \{ h 
\in \MM (\Sigma); h(\delta)= \delta \}$. Let $h \in \Stab(\delta)$. Then we can choose a representative 
$H \in \Diff(\Sigma)$ of $h$ such that either $H \circ d=d$ (if $h$ 
preserves the orientation of $d$), or $H \circ d= d^{-1}$ (if $h$ reverses 
the orientation of $d$). Such a diffeomorphism can be left to a diffeomorphism $H_\Delta: 
\Sigma_\Delta \to \Sigma_\Delta$ (which does not necessarily belong to $\Diff(\Sigma_\Delta)$), and 
$H_\Delta$ can be extended to a diffeomorphism $H_0 \in \Diff (\Sigma_0, \{ Q_1, Q_2\})$, which is 
unique up to isotopy. Now, the following can be easily proved using classical techniques that can be 
found, for example, in \cite{Ivano2}.

\bigskip\noindent
{\bf Claim 1.} {\it The mapping $h \mapsto H_0$ determines an epimorphism $\mu: \Stab(\delta) \to \MM 
(\Sigma_0, \{Q_1, Q_2\})$ whose kernel is the cyclic group $\langle \tau_d \rangle$ generated by 
$\tau_d$.}

\bigskip\noindent
Let $f,g \in \MM(\Sigma)$ such that $\Delta(f)= \Delta(g)= \{ \delta \}$, and such that $f^m=g^m$ for 
some $m \ge 1$. In particular, $f,g \in \Stab(\delta)$. We have $\Delta( \mu(f))= \Delta (\mu(g))= 
\emptyset$ and $\mu(f)^m= \mu(g)^m$, thus either $\mu(f)$ and $\mu(g)$ are both pseudo-Anosov, or 
they are both periodic.  

\bigskip\noindent
Suppose that $\mu(f)$ and $\mu(g)$ are pseudo-Anosov. Then, by Theorem 4.5, $\mu(f)=\mu(g)$. 
By Claim 1, it follows that there exists $x \in \Z$ such that $g= f \tau_d^x$. Note that $\tau_d$ lies 
in the center of $\Stab (\delta)$, thus $f^m= g^m=f^m \tau_d^{mx}$, therefore $x=0$ and $f=g$.

\bigskip\noindent
Suppose that $\mu(f)$ and $\mu(g)$ are periodic. By Lemma 7.13, there are conjugate. So, there exists 
$h \in \Stab (\delta)$ such that $\mu(g)= \mu(hfh^{-1})$. Since $\mu(f)$ is periodic, there exist $r 
\ge 1$ and $y_1, \dots, y_q \in \Z$ such that $\mu(f)^r = \tau_{c_1}^{y_1} \dots \tau_{c_q}^{y_q}$. 
Furthermore, by Claim 1, there exists $t \in \Z$ such that $f^r= \tau_{c_1}^{y_1} \dots 
\tau_{c_q}^{y_q} \tau_d^t$. Note that $\tau_{c_1}, \dots, \tau_{c_q}, \tau_d$ lie in the center of 
$\Stab(\delta)$, thus $f^r$ also lies in the center of $\Stab(\delta)$, therefore $h f^r h^{-1} = f^r$. 
On the other hand, by Claim 1, there exists $x \in \Z$ such that $g= h f h^{-1} \cdot \tau_d^x$. So,
\[
f^{rm}= g^{rm} = (h f h^{-1})^{rm} \tau_d^{rmx} = f^{rm} \tau_d^{rmx}\,,
\]
thus $x=0$ and $g = hfh^{-1}$.
\qed


\subsection{Proof of Lemma 7.4}

{\bf Proof of Lemma 7.4.}  Let $\Sigma$ be a surface of genus 1 with $q$ boundary components, $c_1, \dots, 
c_q$, where $q \ge 2$. Let $d_1, d_2: \S^1 \to \Sigma$ be two non-separating curves such that $d_1$ is 
not isotopic to $d_2^{\pm 1}$, and $d_1 \cap d_2 = \emptyset$ (see Figure 7.6). Set $\delta_i= \langle 
d_i \rangle$ for $i \in \{ 1,2 \}$, and $\Delta= \{ \delta_1, \delta_2\}$. 
Let $\Sigma_\Delta$ denote the surface obtained from $\Sigma$ cutting along $d_1$ and 
$d_2$. So, $\Sigma_\Delta$ is the union of two surfaces of genus 0, $\Sigma_{\Delta\,1}$ and 
$\Sigma_{\Delta\,2}$. The boundary components of $\Sigma_{\Delta\,1}$ can be assumed to be $c_1, \dots, 
c_p, d_1^{(1)}, d_2^{(1)}$, and the boundary components of $\Sigma_{\Delta\,2}$ to be $c_{p+1}, \dots, 
c_q, d_1^{(2)}, d_2^{(2)}$, where $1 \le p \le q-1$. We have $\Sigma= \Sigma_\Delta/\sim$, where $\sim$ 
is the equivalence relation which identifies $d_i^{(1)}(z)$ with $d_i^{(2)}(\bar z) = 
(d_i^{(2)})^{-1}(z)$ for all $i \in \{ 
1,2\}$ and all $z \in \S^1$.

\bigskip\noindent
For each $i \in \{ 1,2 \}$, we take two copies $\D_1^{(i)}$ and $\D_2^{(i)}$ of the standard disk, and 
we denote by $e_j^{(i)} : \S^1 \to \partial \D_j^{(i)}$ the boundary component of $\D_j^{(i)}$. We 
define the surface $\Sigma_{0\,i} = (\Sigma_{\Delta\,i} \sqcup \D_1^{(i)} \sqcup \D_2^{(i)}) /\sim$, 
where $\sim$ is the equivalence relation which identifies $d_j^{(i)}(z)$ with $e_j^{(i)}(\bar z) = 
(e_j^{(i)})^{-1}(z)$ for 
all $j \in \{1,2\}$ and all $z \in \S^1$. Note that $\Sigma_{0\,1}$ (resp. $\Sigma_{0\,2}$) is a 
surface of genus 0 with $p$ (resp. $q-p$) boundary components. We fix a point $Q_j^{(i)}$ in the 
interior of $\D_j^{(i)}$ for all $i,j \in \{1,2\}$.

\bigskip\noindent
Set $\Stab(\Delta)= \{ h \in \MM (\Sigma); h (\Delta)= \Delta \}$. Let $h \in \Stab(\Delta)$. We can 
choose a representative $H \in \Diff (\Sigma)$ of $h$ such that either $H \circ d_1=d_1$ and $H \circ 
d_2=d_2$, or $H \circ d_1 = d_2$ and $H \circ d_2 = d_1$. Then $H$ determines diffeomorphisms 
$H_{\Delta \, i}: \Sigma_{\Delta\, i} \to \Sigma_{\Delta\, i}$, $i \in \{1,2\}$, and these 
diffeomorphisms extend to diffeomorphisms $H_{0\,i} \in \Diff( \Sigma_{0\, i}, \{ Q_1^{(i)}, Q_2^{(i)} 
\})$, $i \in \{1,2\}$. Now, 
in the same manner as Claim 1 in the proof of Lemma 7.3, the following can be easily proved using
classical techniques.

\bigskip\noindent
{\bf Claim 1.} {\it Write $\Gamma_i= \MM (\Sigma_{0\, i}, \{ Q_1^{(i)}, Q_2^{(i)} \})$ and denote by 
$\varphi_i: \Gamma_i \to \Sym_2$ the natural epimorphism. Then the mapping $h \mapsto (H_{0\,1}, 
H_{0\,2})$ determines a homomorphism $\mu: \Stab(\Delta) \to \Gamma_1 \times \Gamma_2$. The kernel of 
$\mu$ is the subgroup $\langle \tau_{d_1}, \tau_{d_2} \rangle$ generated by $\{ \tau_{d_1}, \tau_{d_2} 
\}$, which is isomorphic to $\Z^2$. The image of $\mu$ is $\{ (h_1,h_2) \in \Gamma_1 \times \Gamma_2; 
\varphi_1 (h_1)= \varphi_2 (h_2) \}$.}

\bigskip\noindent
Let $f,g \in \MM (\Sigma)$ such that $\Delta (f)= \Delta (g)= \Delta$ and $f^m = g^m$ for some $m \ge 
1$. Note that $f,g \in \Stab( \Delta)$. Set $(f_1, f_2)= \mu (f)$ and $(g_1, g_2)= \mu(g)$. Then 
$f_i^m= g_i^m$ and $\Delta(f_i)= \Delta(g_i)= \emptyset$, for $i=1,2$. In particular, $f_i$ and $g_i$ 
are either both pseudo-Anosov, or both periodic.

\bigskip\noindent
{\bf Claim 2.} {\it For $i \in \{ 1,2 \}$, there exist $\tilde h_i \in \mu^{-1} (\Gamma_i)$ and $r_i \ge 
1$ such that $g_i= \mu( \tilde h_i) f_i \mu( \tilde h_i)^{-1}$, $( \varphi_i \circ \mu) (\tilde h_i) 
=1$, and $\tilde h_i$ commutes with $\tilde f_i^{r_i}$ for all $\tilde f_i \in \mu^{-1} (f_i)$.}

\bigskip\noindent
We assume $i=1$. The case $i=2$ can be treated in the same way.
If $f_1$ and $g_1$ are both pseudo-Anosov, then, by 
Theorem 4.5, $f_1=g_1$. In this case, $\tilde h_1 =1$ and $r_1=1$ verify the conclusion of Claim 2. So, 
we can assume that $f_1$ and $g_1$ are periodic. By Lemma 7.11, there exists $h_1 \in \Gamma_1$ such 
that $\varphi_1 (h_1)= 1$ and $g_1 = h_1 f_1 h_1^{-1}$. On the other hand, since $f_1$ is periodic, 
there exist $r_1 \ge 1$ and $y_1, \dots, y_p \in \Z$ such that $f_1^{r_1} = \tau_{c_1}^{y_1} \dots 
\tau_{c_p}^{y_p}$. It follows that, for any $\tilde f_1 \in \mu^{-1} (f_1)$, there exist $t_1,t_2 \in 
\Z$ such that $\tilde f_1^{r_1} = \tau_{c_1}^{y_1} \dots \tau_{c_p}^{y_p} \tau_{d_1}^{t_1} 
\tau_{d_2}^{t_2}$. Clearly, such an element lies in the center of $\{ \tilde h \in \mu^{-1} (\Gamma_1); 
(\varphi_1 \circ \mu) (\tilde h) = 1 \}$. In particular, it commutes with any $\tilde h_1 \in \mu^{-1} 
(h_1)$.

\bigskip\noindent
Now, set $h= \tilde h_1 \tilde h_2$ and $r=r_1r_2$. Then $h f^r h^{-1} = f^r$ and $\mu(hfh^{-1}) = 
\mu(g)$. In particular, there exist $x_1,x_2 \in \Z$ such that $hfh^{-1} = g \tau_{d_1}^{x_1} 
\tau_{d_2}^{x_2}$. Furthermore, since $g \in \Stab(\Delta)$, either $g \tau_{d_1} = \tau_{d_1} g$ and 
$g \tau_{d_2} = \tau_{d_2} g$ (if $g (\delta_1)= \delta_1$ and $g(\delta_2)= \delta_2$), or $g 
\tau_{d_1} = \tau_{d_2} g$ and $g \tau_{d_2} = \tau_{d_1} g$ (if $g (\delta_1)= \delta_2$ and 
$g(\delta_2)= \delta_1$).

\bigskip\noindent
Suppose that $g \tau_{d_1} = \tau_{d_1} g$ and $g \tau_{d_2} = \tau_{d_2} g$. Then
\[
g^{rm}= f^{rm} = (hfh^{-1})^{rm} = g^{rm} \tau_{d_1}^{rmx_1} \tau_{d_2}^{rmx_2}\,,
\]
thus $x_1=x_2=0$ and $g= hfh^{-1}$.

\bigskip\noindent
Suppose that $g \tau_{d_1} = \tau_{d_2} g$ and $g \tau_{d_2} = \tau_{d_1} g$. Then
\[
(hfh^{-1})^2 = g \tau_{d_1}^{x_1} \tau_{d_2}^{x_2} g \tau_{d_1}^{x_1} \tau_{d_2}^{x_2} = g^2 
\tau_{d_1}^{x_1+x_2} \tau_{d_2}^{x_1+x_2}\,,
\]
thus
\[
g^{2rm}= f^{2rm} = (hfh^{-1})^{2rm} = g^{2rm} \tau_{d_1}^{rm(x_1+x_2)} \tau_{d_2}^{rm(x_1+x_2)}\,,
\]
therefore $x_1+x_2=0$. We conclude that
\[
\tau_{d_1}^{-x_1} h f h^{-1} \tau_{d_1}^{x_1} = \tau_{d_1}^{-x_1} g \tau_{d_2}^{-x_2} = g \tau_{d_2}^{-
x_1-x_2} = g\,.
\]
\qed


\bigskip\bigskip\noindent
{\bf Christian Bonatti},

\smallskip\noindent
Institut de Math\'ematiques de Bourgogne, UMR 5584 du CNRS, Universit\'e de Bourgogne, B.P. 
47870, 21078 Dijon cedex, France

\smallskip\noindent
E-mail: {\tt bonatti@u-bourgogne.fr} 

\bigskip\noindent
{\bf Luis Paris},

\smallskip\noindent 
Institut de Math\'ematiques de Bourgogne, UMR 5584 du CNRS, Universit\'e de Bourgogne, B.P. 
47870, 21078 Dijon cedex, France

\smallskip\noindent
E-mail: {\tt lparis@u-bourgogne.fr}

\end{document}